\documentclass[12pt]{article}
\usepackage{fullpage}%varioref,multicol

\usepackage{bm}
\usepackage[dvipsnames]{xcolor}
\usepackage{amsmath,amssymb}
\usepackage{amsthm}
\usepackage{wrapfig}
\usepackage{epsfig}
\usepackage{graphicx}
\usepackage{makeidx}
\usepackage{multicol}
\usepackage{tikz}
\usepackage{algorithm}
\usepackage{algpseudocode}
\usepackage{soul}

\usepackage{array,multirow,booktabs,bbm} % for fancier tables

\usepackage{algorithmicx}
\usepackage{algpseudocode}
\makeatletter
\algrenewcommand\ALG@beginalgorithmic{\sffamily\small}
\makeatother
\usepackage{epstopdf}
\usepackage{mathrsfs,mathtools}
\usepackage{xcolor}
\usepackage{enumitem}
\newcommand{\N}{{\cal N}}

\graphicspath{{./Figs/}}
\newcommand{\bs}[1]{\boldsymbol{#1}}
\newcommand{\bmu}{\bs{\mu}}

\newcommand{\bc}{\bs{c}}
\newcommand{\bu}{\bs{u}}
\newcommand{\br}{\bs{r}}
\newcommand{\bv}{\bs{v}}

\newcommand{\bx}{\bs{x}}

\newcommand{\calP}{\mathcal{P}}
\newcommand{\calD}{\mathcal{D}}

\newcommand{\calF}{\mathcal{F}}

\newcommand{\calN}{\mathcal{N}}

\DeclareMathOperator*{\argmax}{argmax}
\DeclareMathOperator*{\argmin}{argmin}

\theoremstyle{remark}

\begin{document}
\graphicspath{{Figs/}}

\title{A hyper-reduced MAC scheme for the parametric Stokes and Navier-Stokes equations
\footnote{%\newline ${   } \quad \ {  }$ 
Y. Chen was partially supported by National Science Foundation grant DMS-1719698 and by the UMass Dartmouth MUST program established by Dr. Ramprasad Balasubramanian via sponsorship from the Office of Naval Research. L. Ji acknowledges the support from NSFC (grant Nos. 12071288 and 21773165). Z. Wang was partially supported by National Science Foundation DMS-1913073.
}
}

\author{
Yanlai Chen\footnote{Department of Mathematics, University of Massachusetts Dartmouth, 285 Old Westport Road, North Dartmouth, MA 02747, USA. Email: {\tt{yanlai.chen@umassd.edu}}.}, 
\,
Lijie Ji \footnote{School of Mathematical Sciences, Shanghai Jiao Tong University, Shanghai 200240, China. Email: {\tt{sjtujidreamer@sjtu.edu.cn}}.},
\,  
Zhu Wang\footnote{
Department of Mathematics, University of South Carolina, Columbia, SC 29208, United
States. Email: {\tt{wangzhu@math.sc.edu}}.} 
}

\date{\empty}

\maketitle

\begin{abstract}
The need for accelerating the repeated solving of certain parametrized systems motivates the development of more efficient reduced order methods. The classical reduced basis method is popular due to an offline-online decomposition and a mathematically rigorous {\em a posterior} error estimator which guides a greedy algorithm offline. For nonlinear and nonaffine problems, hyper reduction techniques have been introduced to make this decomposition efficient. However, they may be tricky to implement and often degrade the online computation efficiency.

To avoid this degradation, reduced residual reduced over-collocation (R2-ROC) was invented integrating empirical interpolation techniques on the solution snapshots and well-chosen residuals, the collocation philosophy, and the simplicity of evaluating the hyper-reduced well-chosen residuals. In this paper, we introduce an adaptive enrichment strategy for R2-ROC rendering it capable of handling parametric fluid flow problems. Built on top of an underlying Marker and Cell (MAC) scheme, a novel hyper-reduced MAC scheme is therefore presented and tested on Stokes  and Navier-Stokes equations demonstrating its high efficiency, stability and accuracy.

\end{abstract}

\section{Introduction}
Developing fast algorithms for parametrized systems is a {critical} topic for researchers in many fields, especially when the problems require repeated numerical simulations. If the dimension of parameter space is high, the computation cost {tends to} be unacceptable. The reduced basis method (RBM) \cite{Quarteroni2015,HesthavenRozzaStammBook} was developed for parametrized PDEs (pPDEs), which is computation efficient and can {produce} an accurate approximation of the original model. Many existing literature and mathematical analysis have proved its effectiveness and cost saving {capability} \cite{BinevCohenDahmenDevorePetrovaWojtaszczyk, CohenDeVore2015}.
The key of RBM to realize cost saving is the offline-online decomposition procedure \cite{Rozza2008,Haasdonk2017Review}. That is, the algorithm calculates several high fidelity solutions through the full order model to construct the reduced basis space and build the reduced order model during the offline process. 
The online process only needs solving the reduced system {whose assembly and resolution are} independent of the {degrees of freedom of the} full order model. The efficiency of this offline-online procedure is trivial for the linear and affine system. 
However, for a nonaffine and nonlinear system, {\em Empirical Interpolation Method} (EIM) 
or its discrete version (DEIM) \cite{Barrault2004,grepl2007efficient, ChaturantabutSorensen2010,PeherstorferButnaruWillcoxBungartz2014} have to be applied to achieve online efficiency.

To illustrate the EIM, we consider the following time dependent linear equation with a nonaffine parameter dependence:
\[ \frac{\partial u}{\partial t}-\calP(u, t, \bx; \bmu) = f, ~\bx \in \Omega, ~t\in [0, T].\]
{Adopting a simple backward Euler scheme, we have} the following linear system {to solve for advancing from $t_k$ to $t_{k+1}$}:
\[ \left(I + \tau A(X^\N; \bmu)\right) \bu^{k+1} = \bu^k +\tau f^\N,\]
where $\bmu$ is the parameter, $\bu^{k} = u^\N(t_k, X^\N)$ represents the discretization of $u$ at time $t_k$, $f^\N$ is the discretization of $f$, $X^\N$ represents all the discretized nodes, $\N$ is the {degrees of freedom} of FDM, $I$ is the identity matrix, $\tau$ is the time step size, and $A$ is the discrete counterpart of operator $-\calP$. 
EIM generates an approximation of the original nonaffine operator $\calP(u, t, \bx; \bmu)$ by a linear combination of $\bmu$-independent {operators}
\[ \calP(u, t, \bx; \bmu) \approx \sum_{q=1}^{Q}\theta_q(\bmu)\calP_q(u, t, \bx).\]
where $\{\theta_q(\bmu)\}_{q=1}^{Q}$ are parameter dependent coefficients, 
$\calP_q(u, t, \bx)\coloneqq \calP(u, t, \bx; \bmu^q)$, 
and $\{\bmu^q\}_{q=1}^{Q}$ are determined through the greedy framework of EIM. Using this affine approximation and letting $V\bc^{k+1}$ denote the reduced order approximation of $\bu^{k+1}$, the reduced system can be written as
\[(V^TV +\sum_{q=1}^Q \theta_q(\bmu)\tau V^TA_qV)\bc^{k+1} = V^TV\bc^k +\tau V^Tf^\N, \]
where $A_q$ is the discretized matrix of the parameter independent operator $-\calP_q$. This reduced form is obtained from Galerkin projection, one of the commonly used variational approaches. Petrov-Galerkin projection \cite{BennerGugercinWillcox2015,CarlbergBouMoslehFarhat2011} {can be adopted as well}.
Obviously, this reduced solver  is not dependent on the dimension of the full order model. However, its 
complexity is linearly dependent on the number of EIM terms $Q$. If $Q$ is too large, assembling this affine approximation will lead to an efficiency degradation. {This may happen when, e.g.} $\Omega$ is a parametrized domain \cite{chen2012certified,BenaceurEhrlacherErnMeunier2018,stabile2020efficient}.

%  2021.3.9
To overcome the degradation from EIM technique, we adopt a (over-)collocation approach as we did in 
\cite{ChenGottlieb2013,ChenJiNarayanXu2020,ChenSigalJiMaday2021}.
Methods in \cite{ChenJiNarayanXu2020,ChenSigalJiMaday2021} are named as {\em over-collocation} method because they choose extra collocation points compared to the collocation methods in \cite{ChenGottlieb2013}. 
{The reduced residual reduced over-collocation (R2-ROC) \cite{ChenSigalJiMaday2021} integrates empirical interpolation techniques on the solution snapshots and well-chosen residuals, the collocation philosophy, and the simplicity of evaluating the hyper-reduced well-chosen residuals for a stable surrogate solver for many pPDEs while keeping the online computational complexity independent of $Q$.} 
{However R2-ROC, applied directly to the more complicated unsteady fluid flow problems, leads to surrogate solutions of oscillating accuracy especially when the RB dimension is low. Motivated by the facts that the instability / lack of robustness is usually ameliorated by adapting or enriching the test space in Finite Element Method and that the classical RBM adopts global basis functions, we design an adaptive strategy to enrich the set of collocation points during the Offline process. The idea is to control the accuracy oscillation in the surrogate solutions by monitoring a {\em robustness indicator}  and making sure that it is essentially non-oscillatory. We systematically select more collocation points, from a well-chosen residual, when a new basis is determined until the robustness indicator falls below a prescribed tolerance. 
We note that, when the robustness indicator stays below this tolerance, the adaptively enriched R2-ROC  degrades to the original R2-ROC algorithm. As our numerical results will show, the original R2-ROC suffices for the steady Stokes and Navier-Stokes equations. That is, a bare minimum of $2N-1$ points for $N$-dimensional RB spaces provide a stable approximation for those cases. However, the new adaptive enrichment is essential in generating reliable surrogate solutions for the unsteady Navier-Stokes equations.}

{We end the introduction by noting that} the collocation idea has also been explored by many other researchers. For instance, authors in \cite{Astrid2008MissingPoint,astrid2004fast} chose collocations based on the missing point interpolation and reduced the amount of original equations that were then solved by the Galerkin method. In \cite{ryckelynck2005priori,ryckelynck2009hyper}, the reduced basis space was constructed through the proper orthogonal decomposition (POD) method, and collocation points were selected by an adaptive algorithm in the context of finite element method. It was also applied to nonlinear dynamical systems with randomly sampled collocation points \cite{bos2004accelerating}. More recently, the reduced quadrature scheme has been developed and the calculation of residuals has been accelerated with sampled integration points. This approach uses less integration points, together with new weights, to generate full residual approximations. Authors 
in \cite{farhat2014dimensional} developed an Energy-conserving Sampling and Weighting hyper reduction (ECSW) method to approximate a vector function using only a reduced mesh, which is determined by solving a non-negative least squares problem. Empirical quadrature procedure (EQP) in \cite{yano2019lp} and Empirical cubature method in \cite{rocha2020adaptive} were all designed with reducing the integration points to accelerate calculations of nonlinear residuals. 
Moreover, model reduction for fluid dynamical system has been investigated by many research groups. Authors in \cite{wang2012proper,iliescu2014variational} used the POD method with eddy viscosity closure models to generate reduced basis space of structurally dominated turbulent flows. Rozza and Stabile \cite{stabile2018finite} also used the POD to construct the reduced basis space of the incompressible Navier-Stokes equations with moderate Reynolds number.
 Ito and Ravindran presented the reduced-order modeling to approximate viscous incompressible flows \cite{ito1998reduced}, where the reduced basis is constructed based on a Lagrange approach. 
Authors in \cite{fick2017reduced} generated the basis space through a greedy framework with a constrained Galerkin projection. 
{Our proposed R2-ROC with the adaptive enrichment differentiate itself via its} greedy framework {for simultaneous basis and} collocation point generation through empirical interpolation procedures \cite{MagicPt_2009,MadayMulaTurinici_GEIMSIAM} {and a novel adaptive mechanism for additional collocation points necessary for unsteady Navier-Stokes equations}.

The rest of paper is organized as follows. In Section \ref{sec:fluids}, we present the full order numerical schemes for the Stokes and Navier-Stokes equations. In Section \ref{sec:R2-ROC-Alg}, we introduce the adaptively enriched R2-ROC method through a generalized time dependent pPDE. In Section \ref{num:final}, numerical results are presented for Stokes equations and incompressible Navier-Stokes equations. Finally, concluding remarks are drawn in Section \ref{sec:conclusion}.

\section{The full order model (FOM) for fluid dynamics}
\label{sec:fluids}

In this section, we introduce three parametric fluid dynamics equations, namely the Stokes equations, steady incompressible Navier-Stokes (NS) equations, and the unsteady incompressible NS equations. As they are well-known, we simply list the equations and present the full order numerical schemes for each of them.

\subsection{Stokes equations}
\label{sec:stokes}
We consider the Stokes problem \cite{MadayMulaPateraYano2015,mu2017simple,rozza2007stability} in a two-dimensional domain $\Omega =[0,1] \times [0,1]$:
 \begin{align} 
 \begin{split}            
- \frac{1}{\text{Re}}(u_{xx} + u_{yy}) +\partial_x p &= 0, ~ \text{in} ~~\Omega,  \\
 - \frac{1}{\text{Re}}(v_{xx} + v_{yy})+ \partial_y p &=  0, ~ \text{in} ~~\Omega,  \\
  u_x + v_y &= 0, ~ \text{in} ~~ \Omega.   
  \end{split}
  \label{eq:stokes}
 \end{align}
The first two are ($x$- and $y$-coordinate) momentum equations and the last one is the so-called continuity equation. Here, $p$ denotes the pressure and $u, v$ represent the velocity components in the $x$- and $y$-directions, respectively. 
We consider the Reynolds number $\text{Re}$  a parameter, and introduce a second parameter in a parametric boundary condition on the top side of the domain while setting the rest of the boundary homogeneous:
\begin{align}
\begin{split} 
u &= 1 + \sin\left(\frac{\pi}{2}\nu x\right), ~\text{on} ~~\Gamma_{\text{top}}=\left\{(x,y): y=1\right\}\\
u &= 0, ~\text{on} ~~\partial \Omega \backslash \Gamma_{\text{top}}\\  
v& = 0, ~ \text{on} ~~\partial \Omega.  
\end{split}
\end{align}
 \begin{figure}[!htb]
\centering
\includegraphics[width=0.32\textwidth]{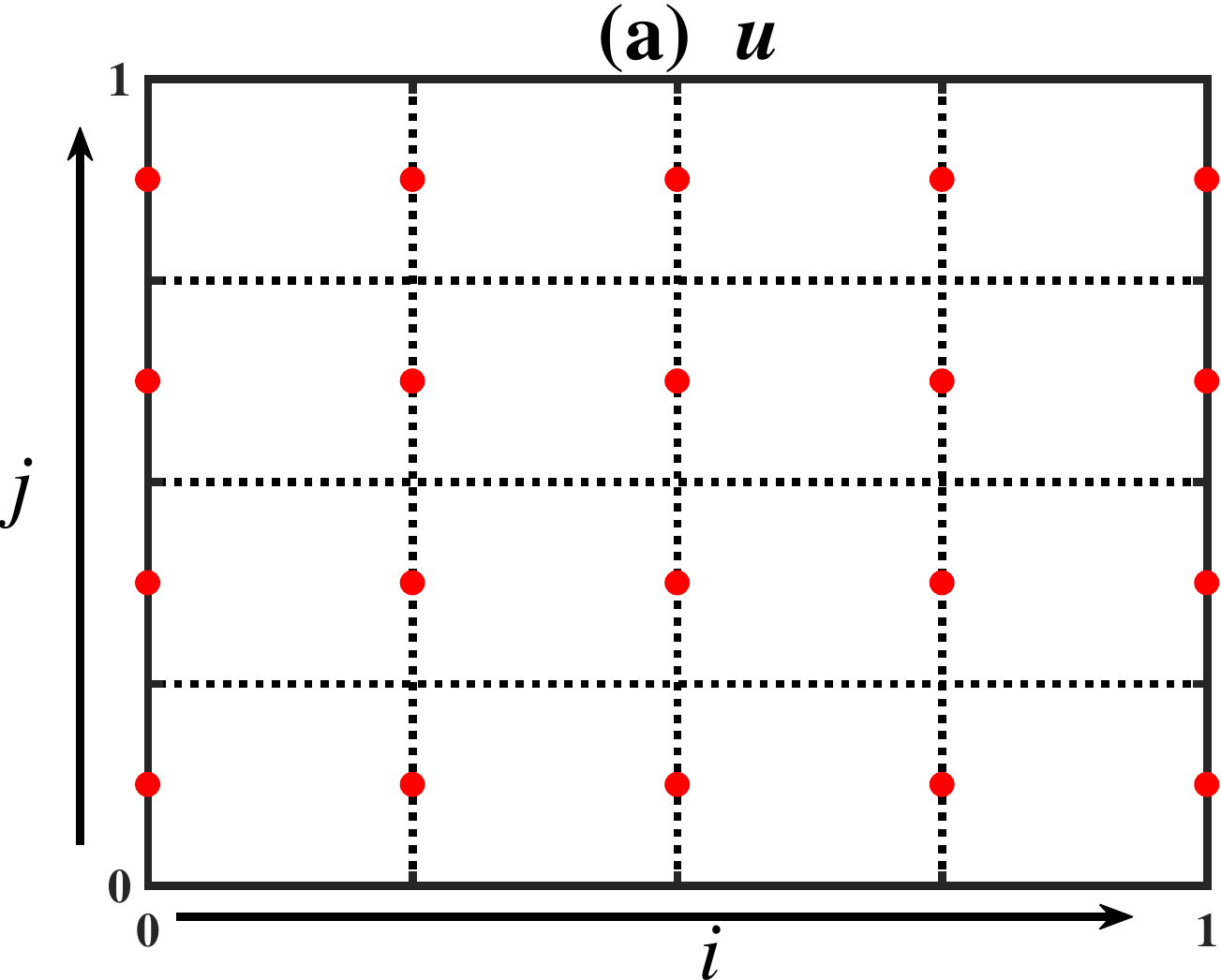}
\includegraphics[width=0.32\textwidth]{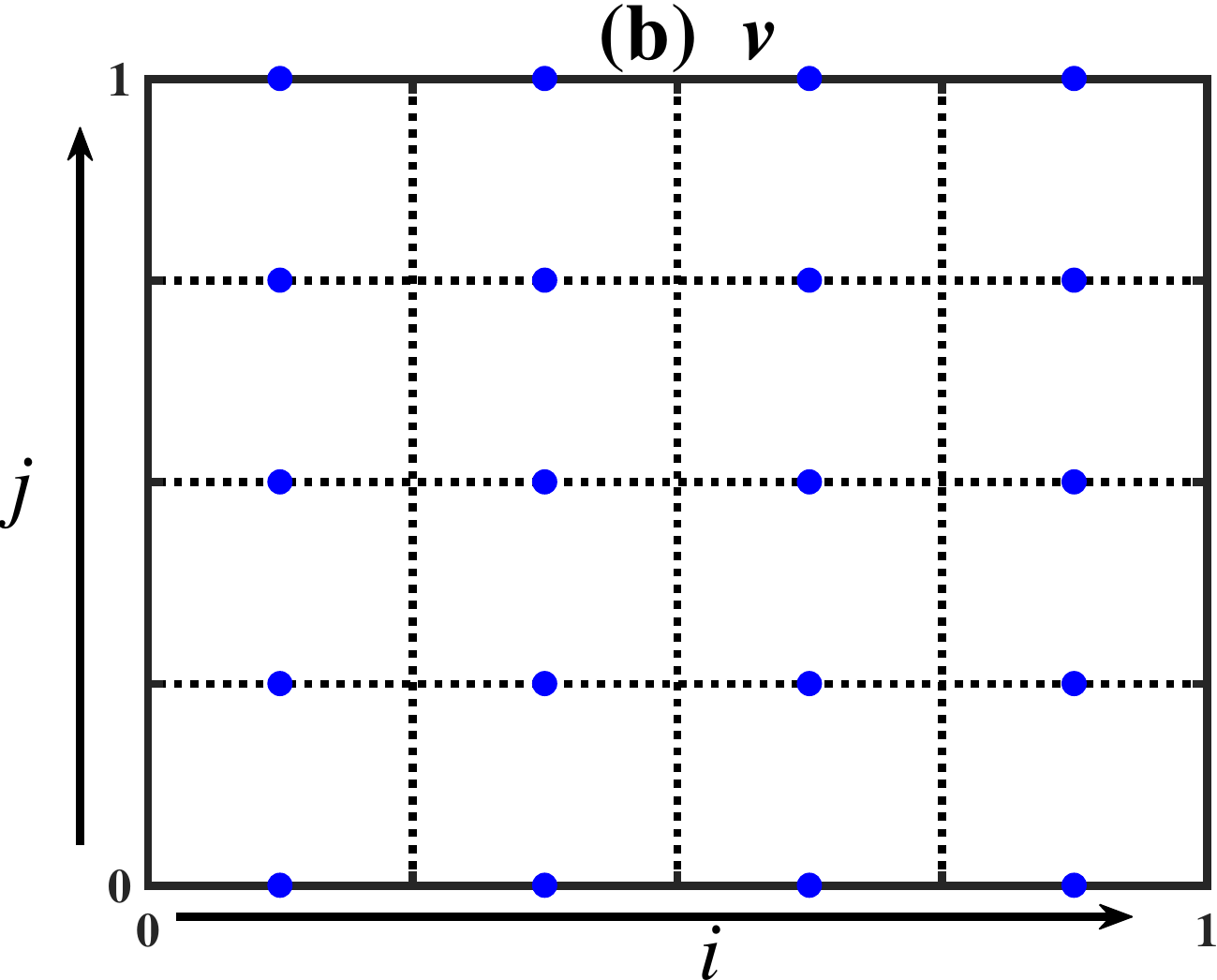}
\includegraphics[width=0.32\textwidth]{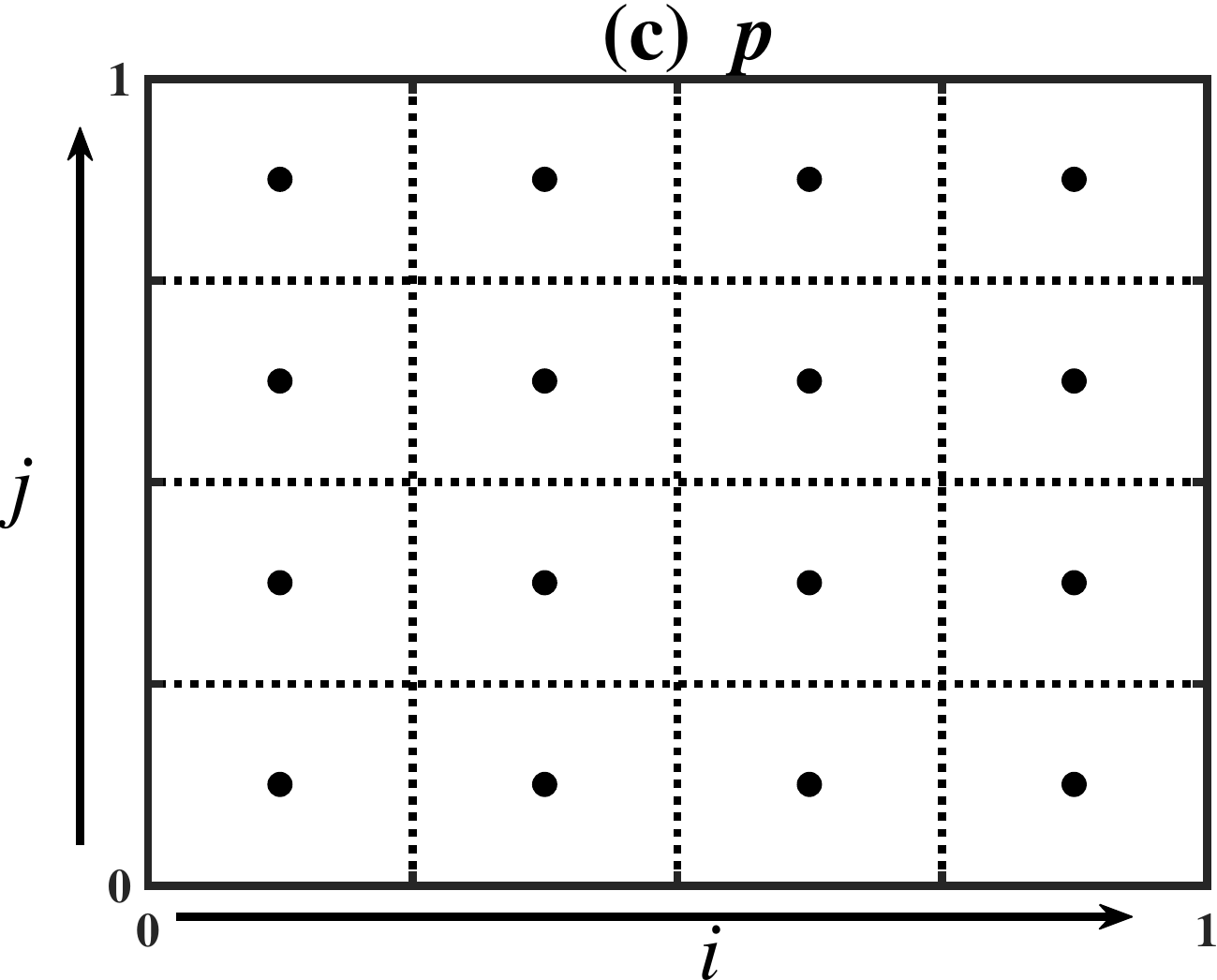}
\caption{Location and indices of $(u, v, p)$ variables {in the MAC scheme.} \cite{MAC_chenlong}.}
\label{fig:mac_uvp} 
\end{figure}
This paper concerns how to solve \eqref{eq:stokes} efficiently when $\bmu = (\text{Re}, \nu)$ changes value in a prescribed domain. Next, we first present the finite difference scheme, Marker and Cell (MAC) \cite{harlow1965numerical,MAC_chenlong,nicolaides1992analysis}, for solving \eqref{eq:stokes} with any given $\bmu$ value. Its idea is to place the unknows $(u, v, p)$ at different locations in a staggered fashion, see Figure \ref{fig:mac_uvp}, and then discretize the $x$-coordinate momentum equation  at vertical edges, $y$-coordinate momentum equation at horizontal edges, and the continuity equation at cell centers. Namely, we have
\begin{align}   
\begin{split}          
- \frac{1}{\text{Re}}\left (\frac{u^\N_{i-1,j} -2u_{i,j}^\N +u^\N_{i+1,j}}{h_x^2} +\frac{u^\N_{i,j-1} -2u^\N_{i,j} + u^\N_{i,j+1}}{h_y^2}\right)+ \frac{p^\N_{i+1,j} - p^\N_{i,j}}{h_x} = 0,\\
- \frac{1}{\text{Re}} \left(\frac{v^\N_{i-1,j} -2v^\N_{i,j} +v^\N_{i+1,j}}{h_x^2} +
\frac{v^\N_{i,j-1} -2v^\N_{i,j} +v^\N_{i,j+1}}{h_y^2}\right) + \frac{p^\N_{i,j+1} - p^\N_{i,j}}{h_y} = 0,\\
\frac{u^\N_{i,j}-u^\N_{i-1,j}}{h_x} + \frac{v^\N_{i,j} - v^\N_{i,j-1}}{h_y} = 0.
 \end{split}
  \label{eq:stokes:scheme}
\end{align}
where $\{u, v, p\}_{i,j}^\N$ denote each variable's value at the $(i,j)$-location of the corresponding grid, see Figure \ref{fig:mac_uvp}. 
The superscript $^\N$, equals to the total number of degrees of freedom, dictates that the quantity is a numerical approximation as opposed to the exact solution. 
Concatenating all unknowns of the velocities and pressure and denoting the resulting vector by $\bu (\bmu) = [u^\N;v^\N;p^\N]$, we can rewrite the numerical scheme~\eqref{eq:stokes:scheme} in a matrix form as:
\begin{align}
A_{\rm S}^\N(\bmu) \bu (\bmu)= f^\N(\bmu),
\label{eq:stokes:matrix}
\end{align}
where $f^\N$ encodes the boundary treatment of $u, v, p$. We end by a remark that, in a multi-query or real-time setting, a system \eqref{eq:stokes:matrix} of size $\N \times \N$ instantiated by a particular $\bmu$ value will have to be solved for each $\bmu$ value.

\subsection{Steady{-state} Navier-Stokes equations}
\label{sec:Nstokes}
Next, we consider the steady{-state} incompressible NS equations \cite{deparis2009reduced,ballarin2015supremizer}. The system becomes nonlinear, but other than that the setup is the same as the Stokes equation. It reads: 
 \begin{align} 
 \begin{split}            
- \frac{1}{\text{Re}}(u_{xx} + u_{yy}) + (u^2)_x +(uv)_y +\partial_x p &= 0, ~ \text{in} ~~\Omega,  \\
 - \frac{1}{\text{Re}}(v_{xx} + v_{yy})+ (uv)_x + (v^2)_y +\partial_y p &=  0, ~ \text{in} ~~\Omega,  \\
  u_x + v_y &= 0, ~ \text{in} ~~ \Omega.  
  \end{split}
  \label{eq:Nstokes}
 \end{align}
Homogeneous boundary conditions are considered on all but the top side of the boundary where $x$-velocity $u$ is parameterized by the boundary control parameter $\nu$, $u = 1 + \nu x$. The MAC scheme is applied in the same fashion on a staggered grid. We adopt the Picard's iteration to solve the resulting nonlinear system. The detailed linearized iterative scheme is as follows
\begin{align}   
\begin{split}          
 - \frac{1}{\text{Re}}\left(\frac{u^{\N, \ell}_{i-1,j} -2u^{\N, \ell}_{i,j} +u^{\N, \ell}_{i+1,j}}{h_x^2} +\frac{u^{\N, \ell}_{i,j-1} -2u^{\N, \ell}_{i,j} + u_{i,j+1}^{\N, \ell}}{h_y^2}\right) + \frac{p^{\N, \ell}_{i+1,j} - p^{\N, \ell}_{i,j}}{h_x} &\\
            +  u_{i,j}^{\N, \ell-1}\frac{u^{\N, \ell}_{i+1,j} -u^{\N, \ell}_{i-1,j}}{2h_x} + \frac{v_{i,j}^{\N,\ell-1} +v_{i,j-1}^{\N,\ell-1}+v_{i+1,j}^{\N,\ell-1} +
 v_{i+1,j-1}^{\N,\ell-1}}{4} \frac{u^{\N, \ell}_{i,j+1} -u^{\N, \ell}_{i,j-1}}{2h_y}  & = 0,\\
-\frac{1}{\text{Re}}\left(\frac{v^{\N, \ell}_{i-1,j} -2v^{\N, \ell}_{i,j} +v^{\N, \ell}_{i+1,j}}{h_x^2} +\frac{v^{\N, \ell}_{i,j-1} -2v^{\N, \ell}_{i,j} +v^{\N, \ell}_{i,j+1}}{h_y^2}\right) + \frac{p^{\N, \ell}_{i,j+1} - p^{\N, \ell}_{i,j}}{h_y} &\\
              + \frac{u_{i,j}^{\N,\ell-1}+u_{i-1,j}^{\N,\ell-1}+u_{i-1,j+1}^{\N,\ell-1}+u_{i,j+1}^{\N,\ell-1}}{4} \frac{v^{\N, \ell}_{i+1,j}-v^{\N, \ell}_{i-1,j}}{2h_x}+ v_{i,j}^{\N,\ell-1} \frac{v^{\N, \ell}_{i,j+1} -v^{\N, \ell}_{i,j-1}}{2h_y}   & = 0,\\
     \frac{u^{\N, \ell}_{i,j}-u^{\N, \ell}_{i-1,j}}{h_x} + \frac{v^{\N, \ell}_{i,j} - v^{\N, \ell}_{i,j-1}}{h_y} & = 0,
 \end{split}
   \label{eq:Nstokes:scheme}
\end{align}
where $\{u, v, p\}_{i,j}^{\N, \ell}$ denote each variable's value at the $(i,j)$-location of the corresponding grid during the $\ell^{\rm th}$ Picard iteration. Concatenating the high-fidelity solutions at the $\ell^{\rm th}$ iteration by $\bu^\ell = [u^{\N, \ell}; v^{\N, \ell}; p^{\N, \ell}]$, we denote the matrix form of this numerical scheme by
\begin{align}
A_{\rm NS}^\N(\bu^{\ell - 1}; \bmu) \bu^{\ell} = f^\N(\bu^{\ell - 1}; \bmu),
\label{eq:Nstokes:matrix}
\end{align}
here $\bmu = (\text{Re}, \nu)$. The stopping criterion of the Picard iteration is that the norm of the full residual $\lVert A_{\rm NS}^\N(\bu^{\ell}; \bmu) \bu^{\ell} - f^\N(\bu^{\ell}; \bmu)\rVert$ is below a prescribed tolerance.

\subsection{Time-dependent Navier-Stokes equations}
\label{sec:lid}

Finally, we consider the time-dependent NS equations \cite{wang2012proper,stabile2018finite} in the two dimensional domain $\Omega =[0,1] \times [0,1]$,
 \begin{align} 
 \begin{split}            
  \frac{\partial u}{\partial t}- \frac{1}{\text{Re}}(u_{xx} + u_{yy}) +(u^2)_x +(uv)_y +\partial_x p &= 0, ~ \text{in} ~~\Omega,  \\
 \frac{\partial v}{\partial t} - \frac{1}{\text{Re}}(v_{xx} + v_{yy})+ (uv)_x + (v^2)_y +\partial_y p &=  0, ~ \text{in} ~~\Omega,  \\
  u_x + v_y &= 0, ~ \text{in} ~~ \Omega,   
  \end{split}
  \label{eq:lid}
 \end{align}
with a given initial condition. As in previous cases, homogeneous boundary conditions are considered on all but the top of the boundary where $x$-velocity $u$ is designated to be $u = 1+\nu x$. The spatial MAC scheme for this case is identical to \eqref{eq:Nstokes:scheme} which means we are solving 
\begin{align}
\left(
\begin{tabular}{c}
$u^{\N,\ell}$\\
$v^{\N,\ell}$\\
$0$
\end{tabular}
\right)_t
+A_{\rm NS}^\N(\bu^{\ell - 1}; \bmu) \bu^{\ell} = f^\N(\bu^{\ell - 1}; \bmu),
\label{eq:NstokesT:matrix}
\end{align}
at each Picard iteration. With backward Euler being our time discretization for simplicity, we have the following fully discretized scheme
\begin{align}
\left(
\begin{tabular}{c}
$u^{\N,\ell}_{(k+1)}$\\
$v^{\N,\ell}_{(k+1)}$\\
$0$
\end{tabular}
\right)
+\tau \, A_{\rm NS}^\N(\bu^{\ell - 1}_{(k+1)}; \bmu) \bu^{\ell}_{(k+1)} = \tau  f^\N(\bu^{\ell - 1}_{(k+1)}; \bmu) + 
\left(
\begin{tabular}{c}
$u^{\N,\ell}_{(k)}$\\
$v^{\N,\ell}_{(k)}$\\
$0$
\end{tabular}
\right).
\label{eq:NstokesT:matrix2}
\end{align}
Here $\tau$ is the time stepsize, subscript $(k)$ denotes the $k^{\rm th}$ timestep value of the corresponding variable while the superscript $\ell$ continues to index the Picard iteration.

\section{The reduced over-collocation method with {adaptive enrichment}}
\label{sec:R2-ROC-Alg}

Following the FOMs presented in the last section, we introduce our proposed Reduced Order Model (ROM) for these three parametric fluid dynamics equations. It is an augmentation of the Reduced Over Collocation (ROC) approach introduced in \cite{ChenJiNarayanXu2020,ChenSigalJiMaday2021} that was developed for steady-state and time-dependent nonlinear problems. The trademark feature {of ROC} is the immunity of the degradation in online efficiency suffered by regular RBM as a result of the EIM-like expansion of the nonlinear and nonaffine terms. 
In what follows, we provide a review of R2-ROC following \cite{ChenSigalJiMaday2021} in Section \ref{sec:r2roc} and then describe in Section \ref{sec:augment} the proposed augmentation, {an adaptive enrichment of the collocation points that is automatically turned for the time-dependent case}.

\subsection{Reduced residual reduced over-collocation (R2-ROC) Method}
\label{sec:r2roc}
{The R2-ROC algorithm is applicable to both steady-state and time-dependent equations. To illustrate it,} we consider the following general time dependent nonlinear and nonaffine PDE:
\begin{align}
\begin{split}
 \frac{\partial u}{\partial t}-\calP(u, t, \bx; \bmu)& = g,~\bx\in \Omega,~t \in [0,T]\\
 \end{split}
 \label{eq:pde}
\end{align}
with necessary boundary and initial conditions. Here $\calP$ is a nonlinear differential operator of $u$ that is parameterized, perhaps in a  nonaffine fashion, by $\bmu$, and $g$ is a forcing term. 
We assume access to a high fidelity solver for \eqref{eq:pde} that is written in a nodal form. To fix the notation, we use $X^\N$ to denote all the grid points, $\calP_\N$  and $g^\N$ represent their discretized counterparts of $\calP$ and $g$ respectively. 
The FOM solution corresponding to parameter $\bmu$ at time $t$ is denoted by $u^\N(t, \bx; \bmu)$ which we assume is close enough to the exact solution $\bu(t, \bx; \bmu)$ for us to adopt as a reference for the ROM. 
Lastly, we will be referring to the (full or reduced version of the) residual when a surrogate solution $v$ is found. This residual is also written in a nodal form.
\begin{align}
r(\bv) = \bv_t -\calP_\N(\bv, t; \bmu) -g^\N.
\end{align}
Now we are ready to briefly review the R2-ROC algorithm. It has two components: an online (reduced) solver of size $n$ that is between $1$ and $N$ with $N$ usually much smaller than $\N$, and an offline training component which repeatedly calls the online solver of increasing size $n$ to build up a surrogate solution space from scratch dimension-by-dimension. 
\subsubsection*{Online solver}
Given the reduced space $W_n$ and a collocation set $X^m \, (\subset X^\N)$ of cardinality $m$ that is comparable to $n$, R2-ROC identifies a surrogate solution for any specific parameter $\bmu$ in the following form
\begin{align}
\widehat{\bu}_n(\bmu, t) = W_{n} \bc_n (\bmu, t).
\end{align}
Here, for simplicity of notation, we also adopt $W_{n}$ for the snapshot matrix whose column space forms the reduced space $W_{n}$.
The unknown coefficients $\bc_n(\bmu, t) \in {\mathbb R}^{n \times 1}$ is obtained by minimizing a subsampled residual
\begin{align}
\bc_n(\bmu, t) = \argmin_{\omega} \left\lVert P_* \left( \left(W_n \omega\right)_t -\calP_\N(W_n\omega, t; \bmu) -g^\N \right)\right\rVert_{\ell^2(\mathbb{R}^m)}.
\label{eq:pde:reduced_t}
\end{align}
The subsampling matrix $P_* \in \mathbb{R}^{m \times \N}$, RB space $W_{n}$, and the reduced collocation set $X^m$  will be generated in the offline process that is described next. 

\subsubsection*{Offline training}

The offline component utilizes  a parameter-time greedy framework \cite{grepl2007efficient,grepl2005posteriori} to iteratively construct  the reduced basis space $W_n$ and subsequently enrich the collocation set $X^m$ which determines the subsampling matrix $P_*$.  
To describe this iterative procedure, we need to fix more notations. 
We discretize the time domain by ${\mathcal T}_f \coloneqq \{t_i: \, i =0, \cdots, \calN_t\}$ with $t_0=0$ being the initial time and $\calN_t =T/ \tau$ where $\tau$ is the temporal step-size and $T$ the final time. The set of reduced time nodes, denoted by ${\mathcal T}_r$, is a subset of ${\mathcal T}_f$ and will be gradually enriched as the new bases are generated. 
In addition, we need the GEIM \cite{MadayMulaTurinici_GEIMSIAM} interpolating functional as described in \cite{ChenSigalJiMaday2021}. It is defined for any admissible function $v(\bx)$ any $\bx \in \Omega$, $t \in [0, T]$, and $\bmu \in \calD$
\begin{equation}
\sigma_{\bx, t}^{\bmu}(v) = v_t -\calP(v(t,\bx;\bmu), t;\bmu).
\label{eq:functional}
\end{equation}
Here, we will also adopt $\sigma_{\bx, t}^{\bmu}(\cdot)$ for its discretized form.
The algorithm judiciously identifies different parameter-time pairs 
\[
\left\{(\bmu^1, t^1_{\bmu^1}), \dots, (\bmu^1, t^{k_{\bmu^1}}_{\bmu^1}); \dots \dots (\bmu^n, t^1_{\bmu^n}), \dots, (\bmu^n, t^{k_{\bmu^n}}_{\bmu^n})\right\}
\]
one-by-one (but not necessarily in the order above), and construct the reduced basis space via the corresponding snapshots. 
With these notations set, we start the greedy procedure with a randomly chosen $\bmu^1$ and obtain the 
snapshots $\{\bu(t_i, X^\N; \bmu^1)\}_{i=0}^{\calN_t}$ by the high fidelity algorithm. ${\mathcal T}_r$ is initiated by the time instant when the corresponding snapshot has the largest variation. That is,
\[
  {\mathcal T}_r = \{t_{\bmu^1}^{1}\} \mbox{ where } t_{\bmu^1}^{1} = \argmax_{t \in {\mathcal T}_f} \left( \max_{x \in X^\N} \bu(t, X^\N; \bmu^1)- \min_{x \in X^\N} \bu(t, X^\N; \bmu^1) \right).
\]
The RB space $W_1$ is then set as $W_1 =\{u_1\}= \{\bu(t_{\bmu^1}^{1}, X^\N; \bmu^1)\}$, and the first collocation point chosen as the GEIM point of the first basis $\bx_\ast^1 = \argmax_{x \in X^\N} |\sigma_{\bx, t_{\bmu^1}^{1}}^{\bmu^1}(u_1)|$ \cite{ChenSigalJiMaday2021}. When the first pair $(\bmu^1, t^1_{\bmu^1})$ is determined, 
we use the online procedure described above to obtain an RB approximation $\widehat{\bu}_{n}(\boldsymbol{\mu},t)$ for each parameter $\bmu$ in $\Xi_{\rm train}$ (a discretization of the parameter domain $\calD$) 
and compute its error estimator $\Delta_n^{RR_t}(\bmu)$. 
\begin{align}
\Delta_n^{RR_t}(\bmu) \coloneqq \sum_{t \in {\mathcal T}_f}  \varepsilon^{RR}(t; \bmu)  \,\, \mbox{ with } \,\, \varepsilon^{RR}(t; \bmu) \coloneqq \lVert P_\ast \br_n(t; \bmu)\rVert_\infty.
\label{eq:delta_t}
\end{align}
Here, 
\begin{equation}
\br_n(t; \bmu) = W_n \bc_t(\bmu) -\calP_\N(W_n\bc(\bmu), t; \bmu) -g^\N,
\label{eq:fullresidual:rb}
\end{equation}  
is the {\em full} residual for the current RB approximation $\widehat{\bu}_n(\bmu, t)$ at time $t$ of parameter $\bmu$. $P_\ast \br_n(t; \bmu)\in \mathbb{R}^{m \times 1}$ then represents its {\em reduced} (subsampled) version.\footnote{The conventional error estimate calculates the negative-order norm of the residual and scales it by the (parametric) stability factor. It is challenging to compute for the nonlinear and nonaffine case with EIM expansion, the involvement of the successive constraint method \cite{HuynhSCM,HKCHP} used to efficiently estimate the parametric stability factor, and the delicacy of evaluating the residual norm even for the linear problem \cite{Casenave2014_M2AN,JiangChenNarayan2019}. This simple error estimator based on the reduced residual \cite{ChenSigalJiMaday2021} has shown to be promising for nonlinear and nonaffine problems without the need of EIM expansion.} Then the next $(\bmu, t)$ pair will be determined in the following greedy fashion \cite{ChenSigalJiMaday2021}:
\begin{enumerate}
\item {\bf Greedy in $\bmu$:} 
The greedy choice for the $\bmu$-component of the $(\bmu, t)$ pair is through maximizing $\Delta_n^{RR_t}(\bmu)$ over the training set $\Xi _{\rm train}$:
\begin{align}
  \bmu^{n+1} = \argmax_{\bmu \in \Xi _{\rm train}} \Delta_n^{RR_t}(\bmu),~~
  \Delta_n = \Delta_n^{RR_t}(\bmu^{n+1}).
  \label{eq:errorestimator:rb}
\end{align}

\item {\bf Greedy in $t$:} 
Once we have a new $\bmu^{n+1}$, we calculate the ROC approximations $\widehat{\bu}_n(\bmu^{n+1}, t) = W_{n} \bc_n (\bmu^{n+1}, t)$ for all time levels $t \in {\mathcal T}_f$, and then the greedy $t$-choice is given as
\begin{align}
t_{\bmu^{n+1}}^{k_{\bmu^{n+1}}} \coloneqq \argmax_{t \in {\mathcal T}_f}\left\{\varepsilon^{RR}( t;\bmu) \coloneqq {\lVert P_* \br_n (t;\bmu^{n+1})\lVert}_\infty\right\}, \mbox{ and } {\mathcal T}_r = {\mathcal T}_r \bigcup \{t_{\bmu^{n+1}}^{k_{\bmu^{n+1}}}\}.
\end{align}
\end{enumerate}
{\bf $X^m$ expansion:} With the new selected $(\bmu^{n+1}, t_{\bmu^{n+1}}^{k_{\bmu^{n+1}}})$, we solve for the truth approximations $\bu(t, X^\N; \bmu^{n+1})$ for $t \le t_{\bmu^{n+1}}^{k_{\bmu^{n+1}}}$. We then obtain the first additional collocation point  $\{\bx^{n+1}_*\}$ from the GEIM process of $\bu(t_{\bmu^{n+1}}^{k_{\bmu^{n+1}}}, X^\N; \bmu^{n+1})$, see Algorithm \ref{alg:GEIM}, and the second additional point $\{\bx^{n}_{**}\}$ by the EIM process of the full residual $\br_n(t_{\bmu^{n+1}}^{k_{\bmu^{n+1}}}; \bmu^{n+1})$, see Algorithm \ref{alg:EIM}. 
The analysis in \cite{ChenSigalJiMaday2021} demonstrates the importance of retaining these two sets of points for producing accurate approximations for both the reduced solution and the residual corresponding to each parameter.

We end this subsection by noting that the GEIM interpolating functional for the steady-state problems is defined as: for any admissible function $v(\bx)$, any $\bx \in \Omega$, and $\bmu \in \calD$: $\sigma_{\bx}^{\bmu}(v) = \calP(v(\bx;\bmu);\bmu)$.
The algorithm now only needs to identify different parameters $
\left\{\bmu^1, \dots, \dots \bmu^n \right\}$
and the error estimator is calculated by $\Delta_n^{RR}(\bmu) \coloneqq \lVert P_\ast \br_n( \bmu)\rVert_\infty$  with $\br_n(\bmu) = \calP_\N(W_n\bc(\bmu); \bmu) -g^\N$
being the {\em full} residual for the current RB approximation $\widehat{\bu}_n(\bmu)$ of parameter $\bmu$. The $(n+1)$-th parameter is still chosen by maximizing $\Delta_n^{RR}(\bmu)$ over the training set.
%\begin{align}
%  \bmu^{n+1} = \argmax_{\bmu \in \Xi _{\rm train}} \Delta_n^{RR}(\bmu).
%  \label{eq:errorestimator:rb:steady2}
%\end{align} 

\begin{algorithm}[h]
\begin{algorithmic}[1]
\vspace{0.5ex}
\State \textbf{Input:} current basis space $W_n$, interpolation operators $\{\sigma_1, \cdots, \sigma_n\}$, collocation set $X^m$, current subsampling operator matrix $P_*$, parameter-time pair $(\bmu^{n+1}, t_{\bmu^{n+1}}^{k_{\bmu^{n+1}}})$ and new basis $\xi_{n+1}$.
\State Compute a generalized interpolatory residual for vector $\xi_{n+1}$:

Find $\xi_{n+1} = {\xi_{n+1}} - W_{n}\alpha$ with $\{\alpha_i\}$ determined by {requiring that} $\sigma_i(\xi_{n+1})=0$ for $i= 1, \ldots, n$.

Find $\bx^{n+1}_*=\argmax_{\bx \in X^\N/X^m} \left|\sigma_{\bx, t_{\bmu^{n+1}}^{k_{\bmu^{n+1}}}}^{\bmu^{n+1}}(\xi_{n+1})\right|$.

Define $\sigma_{n+1}(\cdot) := \sigma_{\bx^{n+1}_*, t_{\bmu^{n+1}}^{k_{\bmu^{n+1}}}}^{\bmu^{n+1}}(\cdot)$.
$\xi_{n+1}=\xi_{n+1}/\sigma_{n+1}(\xi_{n+1})$ and let $i_{1}$ be the index of $\bx^{n+1}_*$ in $X^\N$.
\State \textbf{Output:} $\{\sigma_1, \cdots, \sigma_{n+1}\}$, $i_{1}$, $\xi_{n+1}$, $\{\bx^{n+1}_*\}$.
\end{algorithmic}
\caption{\sf GEIM algorithm {for the snapshots}}
\label{alg:GEIM}
\end{algorithm}

\begin{algorithm}[h]
\begin{algorithmic}[1]
\vspace{0.5ex}
\State \textbf{Input:} current basis space $W_n$, current subsampling operator matrix $P_*$, interpolatory residual vectors $\{r_1, \cdots, r_{n-1}\}$, collocation set $X^m$, residual-based collocation set $X_r^{n-1}$, collocation $\bx_{*}^{n+1}$ and parameter-time pair $(\bmu^{n+1}, t_{\bmu^{n+1}}^{k_{\bmu^{n+1}}})$.
\State Calculate the full residual $r_{n}(t_{\bmu^{n+1}}^{k_{\bmu^{n+1}}}, \bmu^{n+1})$ by eq.~\eqref{eq:fullresidual:rb}. 

\State Compute an interpolatory residual $r_{n}$: 

Find $r_{n} = r_{n} - \sum_{j=1}^{n-1}\alpha_j r_j$ with $\{\alpha_j\}$ determined by {requiring that} $r_{n}(X^{n-1}_r)=0$. 

Find $\bx_{**}^n=\argmax_{\bx \in X^\N/\left\{X^m, \bx_{*}^{n+1}\right\}} |r_{n}|$.

Calculate $r_{n}=r_{n}/ r_{n}(\bx_{**}^n)$, and let $i_2$ be the index of $\bx_{**}^n$ in $X^\N$.

%\\
\State \textbf{Output:} $\{r_1, \cdots, r_{n-1}, {r_n}\}$, $i_{2}$, $\{\bx^{n}_{**}\}$. 
\end{algorithmic}
\caption{\sf EIM algorithm for the residuals}
\label{alg:EIM}
\end{algorithm}

\subsection{Adaptive {enrichment of R2-ROC}}
\label{sec:augment}
The {original} R2-ROC algorithm, if applied {directly} to the more complicated unsteady fluid flow problems targeted in this paper, leads to surrogate solutions of oscillating accuracy. This oscillation is usually more pronounced when the RB dimension is low. 
{The instability / lack of robustness is usually ameliorated by enriching the test space in Finite Element Method. This observation and the fact that the classical RBM and other linear MOR methods adopt global basis functions motivate us to design} an adaptive strategy to {enrich the set of} collocation points during the Offline process. This idea is realized by {calculating a {\em robustness indicator} and then systematically} selecting more collocation points {in a try-and-error fashion} when a new basis is determined {until the robustness indicator falls below a prescribed tolerance}. 

{To present the detail of our adaptive enrichment, we first recall that the} $n$th parameter-time pair is determined based on $\Delta_{n-1}$. {Once it is obtained,} we substitute the $(n-1)$ dimensional ROC approximation, {at the chosen parameter-time pair,} into the {full order model for} its residual $\br_{n-1}$,
$$\br_{n-1} = {(\widehat{\bu}_{n-1})}_t(t_{\bmu^{n}}^{k_{\bmu^{n}}};\bmu^{n}) -\calP_\N(\widehat{\bu}_{n-1}(t_{\bmu^{n}}^{k_{\bmu^{n}}};\bmu^{n}), t_{\bmu^{n}}^{k_{\bmu^{n}}}; \bmu^{n}) - g^\N,$$
{After an EIM orthonormalization with previous such residuals,} 
the (original) ROC algorithm only selects one point, {the maximizer of $|r_{n-1}|$, in addition to the} one point from the GEIM process for the solutions. {It then builds the $n$-dimensional ROC space and proceed to select the $(n+1)$th parameter-time pair, based on $\Delta_n$. }

{The essence of our augmentation is to define a {\em robustness indicator}
\[
\rho_\Delta^n =\frac{\Delta_n}{\Delta_{n-1}}, 
\]
and control the accuracy oscillation in the surrogate solutions by making sure that $\rho_\Delta^n$ is essentially non-oscillatory. We note that this is a natural choice whenever the effectivity index of $\Delta_n$ remains close to being a constant. This is achieved by revisiting the EIM-orthonormalized residual $\br_{n-1}$ and extract additional collocation points which, in turn, refines the $n$-dimensional ROC solver. We then update error indicator $\Delta_n$ and the robustness indictor $\rho_\Delta^n$. This process is repeated until $\rho_\Delta^n$ is smaller than a prescribed tolerance $\gamma$ or the total number of additional collocation points surpasses a prescribed cap $n_{adap}^{max}$. 
These {\em adaptively enriched} points are sampled according to an inverse CDF (Cumulative Distribution Function) strategy. Specifically, we generate the CDF $\calF_{|\br_{n-1}|}$ of the (discretized) domain of $\br_{n-1}(\bx)$ according to the function's absolute value and then perform a uniform inverse sampling of $n_{adap}$ collocation points from the top $p_{adap}$ percentile. Denoting these additional} collocation points by $\bx_{+}^{n-1}$, {we have
\begin{equation}
\bx_{+}^{n-1} = \calF_{|\br_{n-1}|}^{-1}\left (1-p_{adap}: \frac{p_{adap}}{n_{adap} - 1}: 1\right),
\label{eq:add_colo_formula}
\end{equation}
where $a:h:b$ is the Matlab notation of sampling the interval $[a,b]$ with a stepsize $h$.
}

The second set of {residual-based} collocation points is then updated as $X^{n-1}_r =X^{n-1}_r \cup \{\bx_{+}^{n-1}\}$,  {and the whole} collocation set updated as $X^m =X^m \cup \{\bx_{+}^{n-1}\}$. Correspondingly, a new $n$th parameter will be determined based on this updated error estimator $\Delta_{n}$ {as a consequence of the reduced collocation set being enriched}. The adaptive algorithm is displayed in Algorithm \ref{alg:rectified} and the complete R2-ROC {with adaptive enrichment} for time dependent problems is shown in Algorithm \ref{alg:c:plus:offline:time}.
{We note that, when the robustness indicator $\rho_\Delta^n$ stays below $\gamma$, the adaptively enriched R2-ROC} will degrade to the original R2-ROC algorithm. {As our numerical results will show  in the next section, the original R2-ROC suffices for the steady-state Stokes and Navier-Stokes equations. That is, $2N-1$ points for $N$-dimensional RB spaces produce stable surrogate solutions for those equations. We also note that $\gamma$ represents the allowable oscillation and should be set by the practitioners for their specific problems. We take the following} piecewise function in this paper
\[
\gamma = \left\{
\begin{array}{rcl}
5& & {n<6},\\
2& & {n>=6}.
\end{array}\right.
\]

\begin{algorithm}[h]
\begin{algorithmic}[1]
\vspace{0.5ex}
\State \textbf{Input:} control parameters $\gamma$, {$n_{adap}^{max}$, $n_{adap}$ and $p_{adap}$}, error estimators $\Delta_{n-1}$ and $\Delta_{n}$, full residual $r_{n-1}$, residual-based collocation set $X_r^{n-1}$, solution-based collocation set $X^{n}_s$, current collocation set $X^m$, interpolatory operator matrix $P_*$.
\State \mbox{\textbf{While}} $\Delta_{n} > \gamma \Delta_{n-1}$ \& $n_{adap}<n_{adap}^{max}$ 

\State  \quad $n_{adap} = n_{adap}+10$.
\State  \quad Select $n_{adap}$ elements from the {top $p_{adap}$ percentile of the domain} of $|r_{n-1}|$, $\bx_{+}^{n-1}$, {according to \eqref{eq:add_colo_formula}}.
\State  \quad Update collocation set $X^{n-1}_r = X^{n-1}_r \cup \bx_{+}^{n-1}$, $X^m = X^{n}_s \cup X_r^{n-1}$, $m= \#\{X^{n}_s,X^{n-1}_r\}$, 

and the operator matrix $P_* = [P_*; e_{i_+}^T]$, $i_+$ represents the index of $\bx_{+}^{n-1}$ in $X^\N$. 

\State \quad Solve the reduced problem for $\bc_{n} (\bmu, t_k)$ by eq.~\eqref{eq:pde:reduced_t}. Obtain an updated error estimator 

$\Delta_{n}$ and a new $\bmu^{n+1}$ by eq.~\eqref{eq:errorestimator:rb}.

\State \mbox{\textbf{End While}}

\State \textbf{Output:} $X^{n-1}_r$, $X^m$, $P_*$, $\Delta_{n}$, $\bmu^{n+1}$. 
\end{algorithmic}
\caption{\sf Adaptive {enrichment} algorithm}
\label{alg:rectified}
\end{algorithm}

\begin{algorithm}[h]
\begin{algorithmic}[1]
\vspace{0.5ex}
\State \textbf{Input:} control parameters $\gamma$, {$n_{adap}^{max}$, and $p_{adap}$}.

\State Choose $\bmu^1$, and set $k_{\bmu^1} = 1$. 
$t_{\bmu^1}^{k_{\bmu^1}} = \argmax_{t \in {\mathcal T}_f} \left( \max_x \bu(t, X^\N; \bmu^1)- \min_x \bu(t, X^\N; \bmu^1) \right)$ be the first temporal node. Define $\xi_1\coloneqq \bu(t_{\bmu^1}^{k_{\bmu^1}}, X^\N;\bmu^1)$. 
\State Find $\bx_*^1=\argmax_{\bx \in X^\N} \left|\sigma_{\bx, t_{\bmu^1}^{k_{\bmu^1}}}^{\bmu^1}(\xi_1)\right|$. Define $\sigma_1(\cdot) := \sigma_{\bx_*^1, t_{\bmu^1}^{k_{\bmu^1}}}^{\bmu^1}(\cdot)$ and calculate $\xi_1 = \xi_1/\sigma_1(\xi_1)$. Let $P_* = [e_{i_1}]^T$, where $i_1$ is the index of $\bx_*^1$ in $X^\N$.
\State Initialize $m = n = 1, \, X^m = X^n_s =\{\bx_*^1\}$, $W_1 = \left\{\xi_1 \right\}$, $n_{adap} =0$ and $X_r^0 = \emptyset$. 
\State \mbox{\textbf{For}} $n = 2,\ldots, N$
\State \quad Solve the reduced problem for $\bc_{n-1} (\bmu, t_k)$ by eq.~\eqref{eq:pde:reduced_t}. 
\State \quad Find $\bmu^n = \argmax_{\bmu \in \Xi_{\rm train}}  \Delta_{n-1}^{RR_t}(\bmu)$, $t_{\bmu^n}^{k_{\bmu^n}}= \argmax_{t \in {\mathcal T}_f}\varepsilon^{RR}( t;\bmu^n)$ and 
$\Delta_{n-1} = \Delta_{n-1}^{RR_t}(\bmu^n)$.
\State \quad If $\rho_\Delta^{n-1} > \gamma$, $n_{adap} =0$, perform Algorithm \ref{alg:rectified} to obtain a rectified error estimator $\Delta_{n-1}$ 

and determine a new $\bmu^n$.
\State \quad Solve $\xi_n=\bu (t_{\bmu^n}^{k_{\bmu^n}}, X^\N;\bmu^n)$. 
\State \quad Perform Algorithms \ref{alg:GEIM}-\ref{alg:EIM} and \mbox{update} $W_{n} = \{W_{n-1},\xi_n\}$, $X^n_s= X^{n-1}_s \cup \bx_*^n$, $X^{n-1}_r= X^{n-2}_r \cup \bx_{**}^{n-1}$, 

$m= \#\{X^n_s, X^{n-1}_r\}$, $X^m= X^n_s \cup X_r^{n-1}$, $P_* = [P_*; \,\, (e_{i_1})^T; \, \, (e_{i_2})^T]$. %\\%[0.5ex]

\State  \mbox{\textbf{End For}}

\State \textbf{Output:} $W_{n}$, $X^n_s$, $X^{n-1}_r$, $X^m$, $P_*$.
\end{algorithmic}
\caption{\sf R2-ROC algorithm {with adaptive enrichment} for time dependent problems}
\label{alg:c:plus:offline:time}
\end{algorithm}

\section{Numerical results}
\label{num:final}
In this section, we test the proposed R2-ROC method {with adaptive enrichment} on the three fluid problems described in Section~\ref{sec:fluids}. In order to evaluate the accuracy of our scheme, we measure the relative error $E(n)$ when $n$-dimensional RB space is used for parameters from a given testing set $\Xi_{\rm test}$. It is defined as follows for the steady state problems:
\begin{equation}
E(n) =  \max_{\bmu \in \Xi_{\rm test}}\left\{\frac{\| \bu(\bmu) - \widehat{\bu}_n(\bmu)\|_\infty}{\|\bu(\bmu)\|_{L^\infty(\Omega)}}\right\} \mbox{ where } ||\bu||_{L^\infty(\Omega)} =  \|\bu(\bmu) \|_\infty.
\label{eq:error:steadyns}
\end{equation}
For the time dependent case, we take the $L^\infty$-norm of the quantity above over a discrete time domain ${\mathcal T}_e$. That is,
\begin{equation}
E(n) =  \max_{t \in {\mathcal T}_e} \max_{\bmu \in \Xi_{\rm test}}\left\{\frac{\| \bu(\bmu,t) - \widehat{\bu}_n(\bmu,t)\|_\infty}{\|\bu(\bmu)\|_{L^\infty(\Omega)}}\right\} \mbox{ where } ||\bu||_{L^\infty(\Omega)} = \max_{t \in {\mathcal T}_e} \|\bu(\bmu,t) \|_\infty.
\label{eq:error:timens}
\end{equation}
Here, ${\mathcal T}_e$ can be the full time grid ${\mathcal T}_f$ or a sub-sampled version. In addition, we measure the error of the streamline approximation by:
\begin{align}
E_q(n) = |q^\N - \widehat{q}_n|,
\label{eq:error:streamline}
\end{align} 
where $q^\N$ is the truth streamline function calculated by the full FDM and $\widehat{q}_n$ is its ROC approximation using $n$ basis. Note that the streamline $q$ is related with the velocities $u$ and $v$ via $-\Delta q = \frac{ \partial u}{\partial y} - \frac{\partial v}{\partial x}$.

\subsection{Steady-state Stokes and Navier-Stokes equations}
We first consider the (linear) Stokes and (nonlinear) incompressible Navier-Stokes equations presented in Sections \ref{sec:stokes} and \ref{sec:Nstokes}. The parameter dependence is linear in Reynolds number $\text{Re}$ and nonlinear in the boundary control $\nu$. The parameter domain for $\bmu =(\text{Re}, \nu)$ is $\Omega_p = [10,2000] \times [0,4]$ and $[10,1000] \times [0,2]$  respectively. We first validate our computational model by performing an accuracy test on a sequence of uniform grids for the equation with reference solution as a manufactured solution $u(x,y) = \sin(\pi x), v(x,y) = -\pi \cos(\pi x) y$ for Stokes and a highly refined FDM solution for Navier-Stokes ($n_x =256, n_y =256$). These reference solutions are marked with a superscript $^\text{acc}$ and the $L^2$ norm of the numerical errors are listed in Table~\ref{table:accuracy}, where the second order convergence is observed for both equations.

\begin{table}[thbp]
	\centering
			\begin{tabular}{c|cccc}
		\hline
Equation &$n_x, n_y$ &  $||u^{\text{acc}}-u^{\text{fdm}}||_2$  & $||v^{\text{acc}} -v^{\text{fdm}}||_2$ & $||p^{\text{acc}} -p^{\text{fdm}}||_2$\\ \hline
\multirow{5}{*}{Stokes}
&8 	& 0.5771 & 0.5771 & 0.0258\\ 
&16  & 0.0709 &  0.0709 &  0.0501 \\ 
&32   & 0.0175	& 0.0175 & 0.0125 \\ 
&64 & 0.0044	& 0.0044 & 0.0031 \\ 
&128 & 0.0011	& 0.0011 & 0.000784 \\ \hline
\multirow{5}{*}{Navier-Stokes}
&8 	& 0.3101 & 0.3389 & 1.0508\\ 
&16  & 0.2166 &  0.2108 &  0.4416 \\ 
&32   & 0.0642	& 0.0650 & 0.1749 \\ 
&64 & 0.0159	& 0.0159 & 0.0517 \\  
&128 & 0.0032	& 0.0032 & 0.0112 \\ \hline
	\end{tabular}
	\caption{Accuracy of the MAC scheme for the Stokes ($\text{Re} =10$) and incompressible Navier-Stokes ($\nu =1, \text{Re} =20$) equations, on staggered grids with $n_x \times n_y$ cells.} 
	\label{table:accuracy}
\end{table}
\begin{figure}[thbp]
\centering
\includegraphics[width=0.35\textwidth]{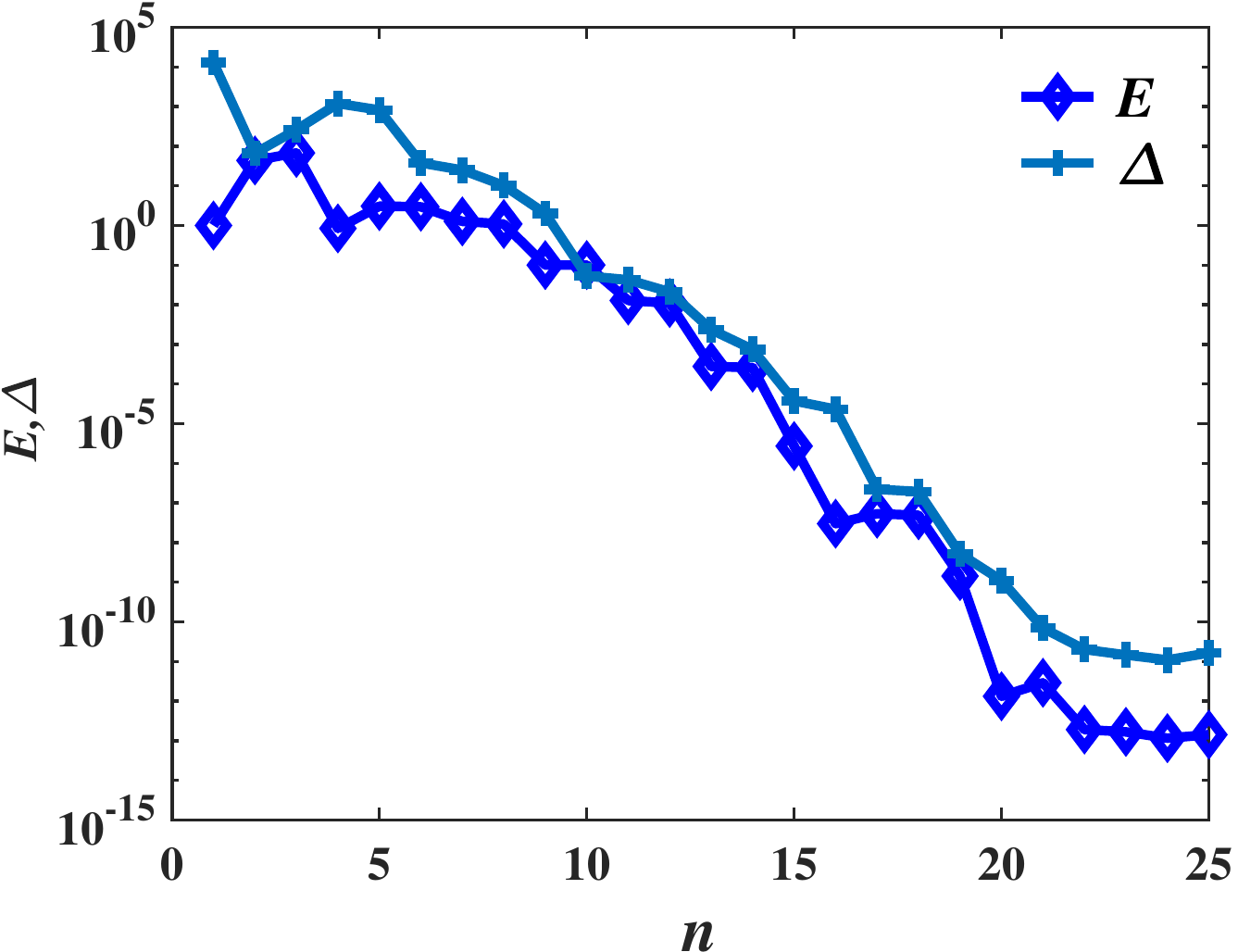}
\includegraphics[width=0.35\textwidth]{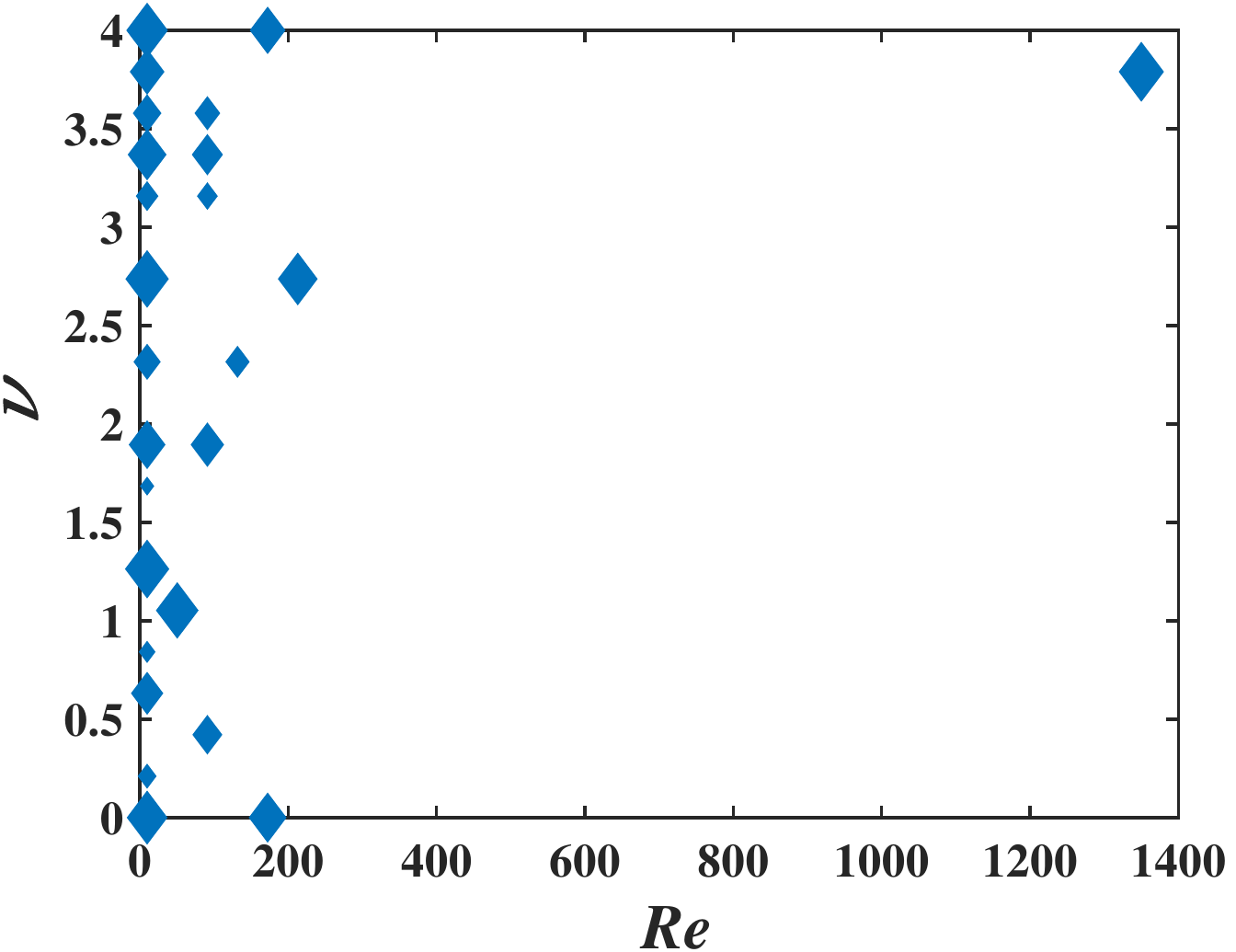}
\includegraphics[width=0.35\textwidth]{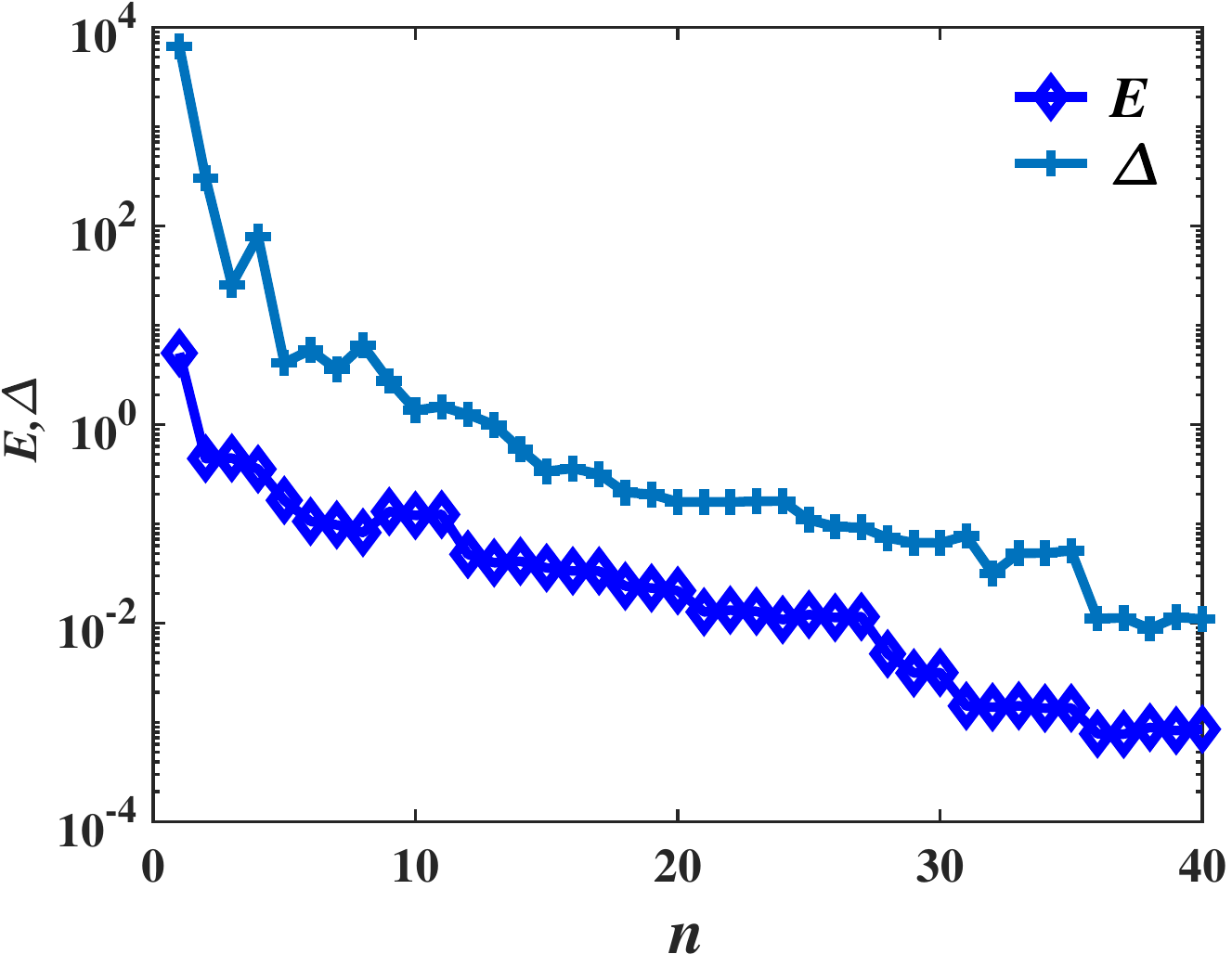}
\includegraphics[width=0.35\textwidth]{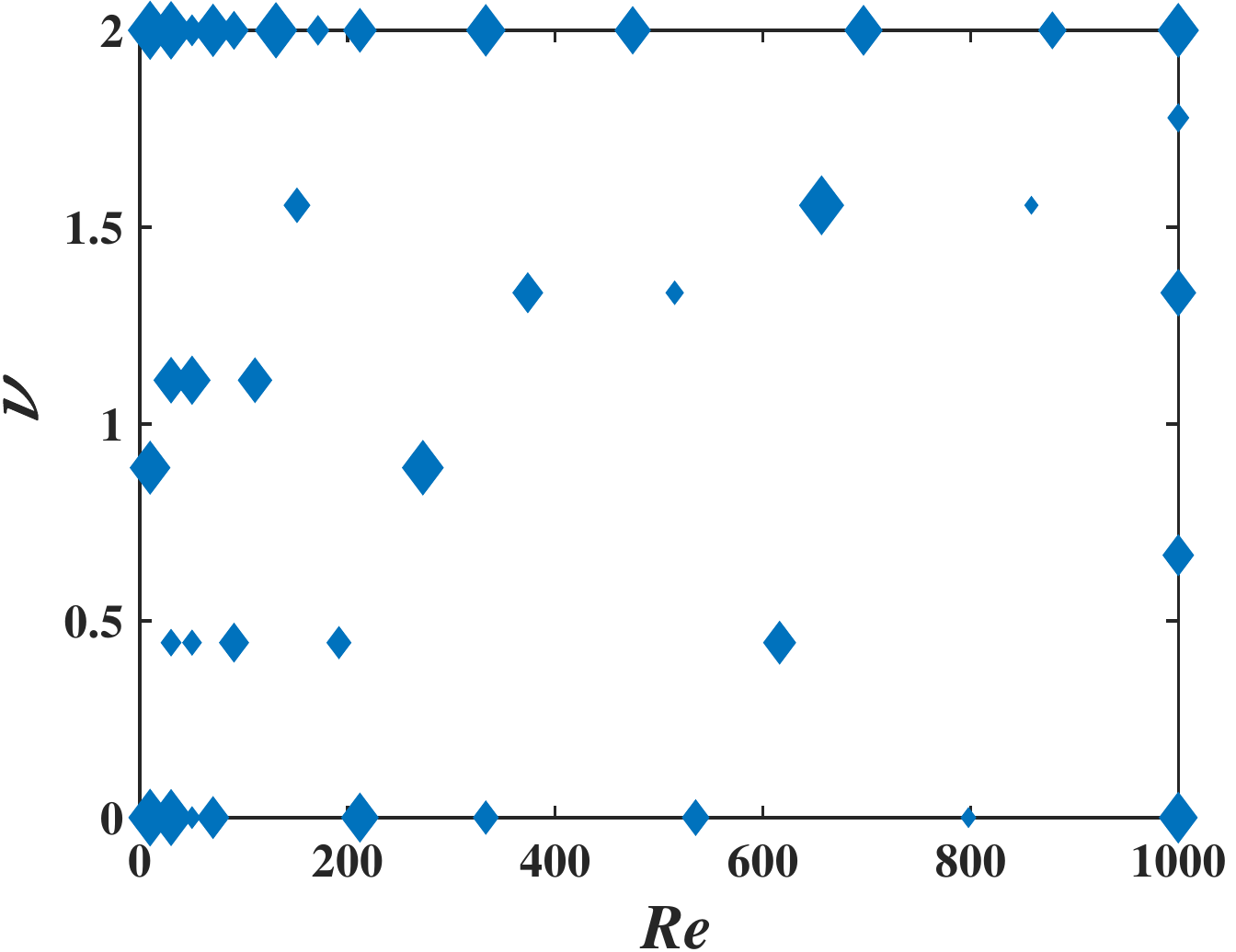}
\caption{\small History of convergence for the relative error and error estimator (left) and the selected parameter values by the ROC algorithm (right). Parameter pair selected earlier has larger symbols. The first parameter is chosen randomly. 
Top is for Stokes and bottom Navier-Stokes.
}
\label{fig:stokes:error} 
\end{figure}

For the ROC approximation, we fix the mesh with $n_x =n_y =200$ ($n_x=120, n_y=120$ for Navier-Stokes) and perform the offline procedure on a training set obtained by discretizing the parameter domain $\Omega_p$ by a $50\times 20 \, (10)$ uniform grid for Stokes (Navier-Stokes). We then test the ROC model on $20 \times 10 \, (5)$ uniform grids from $[10+2.2, 2000 \, (1000) -2.2] \times [0+0.22, 4 \, (2) -0.22]$ for Stokes (Navier-Stokes) that, in particular, are non-overlapping with the training sets. The parameter values chosen by the ROC algorithm during the training process are displayed in Figure \ref{fig:stokes:error} Right, where parameters selected earlier are presented by larger diamonds. It is seen that the selected parameters concentrate in a low-Reynolds-number regime but distribute relatively evenly in the range of $\nu$ for the (linear) Stokes problem. The distribution is more delicate for the Navier-Stokes equation.
The error and error estimator curves are plotted in Figure \ref{fig:stokes:error} Left, both showcasing exponential convergence with respect to basis number $n$. It underscores the effectiveness of our error estimator and the resulting greedy algorithm. 

\begin{figure}[thbp]
\centering
\includegraphics[width=0.3\textwidth]{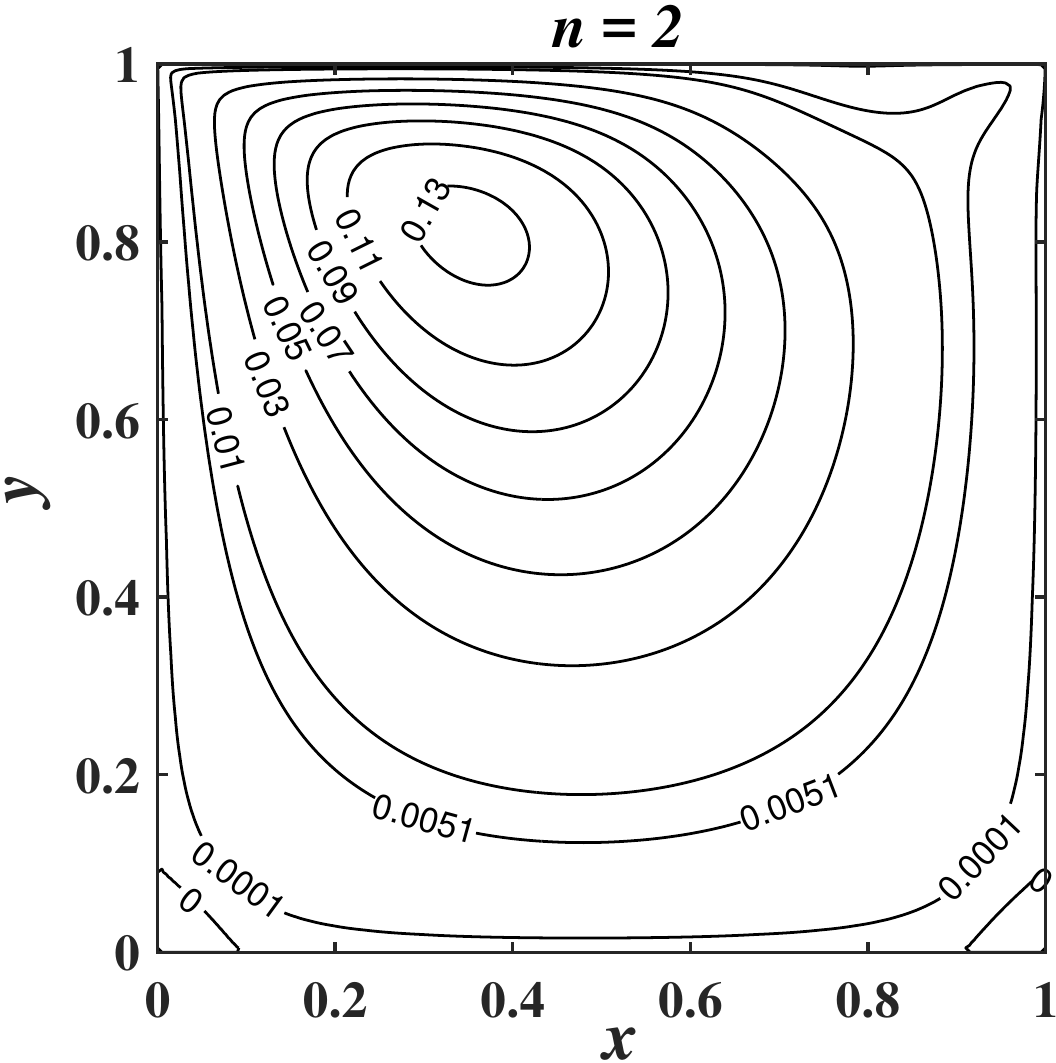} %\hspace{0.4cm}
\includegraphics[width=0.3\textwidth]{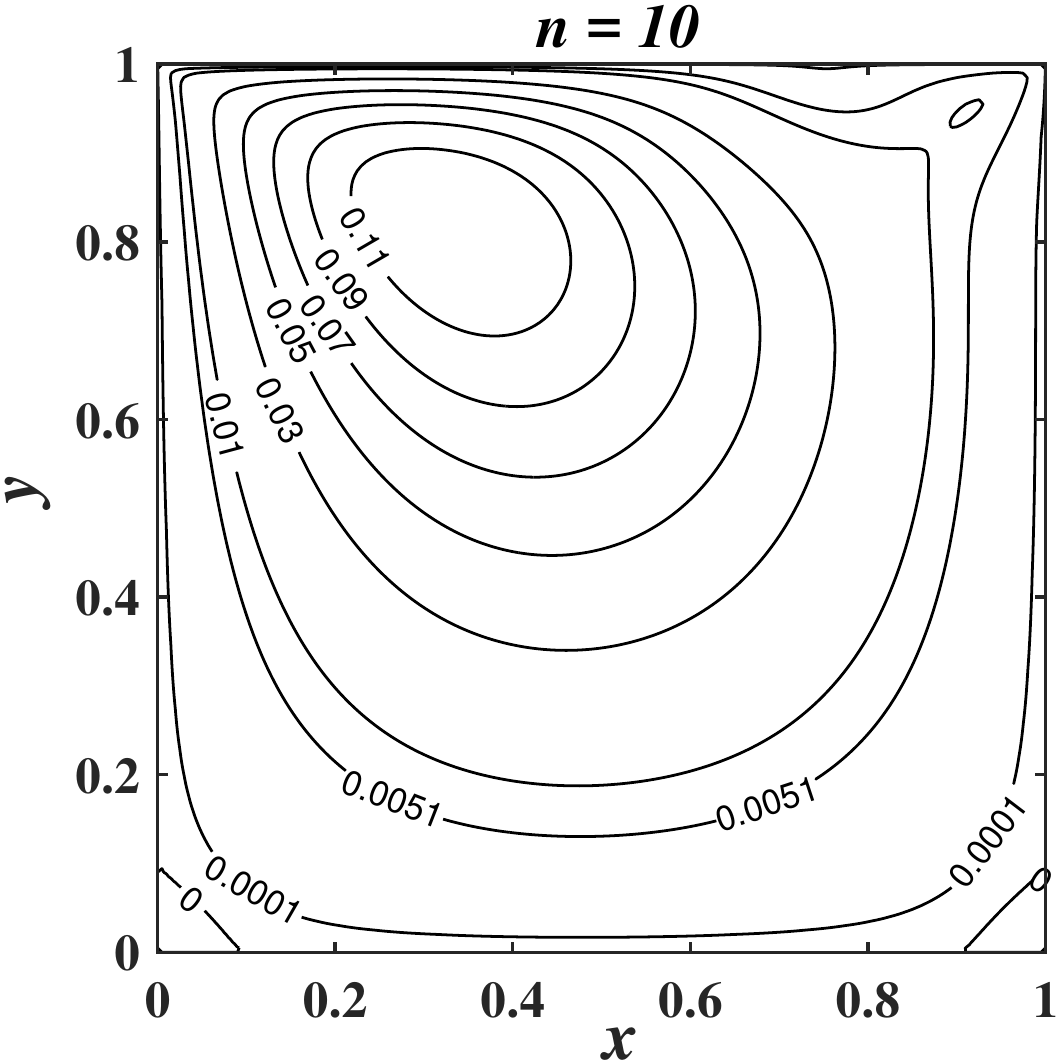} %\hspace{0.4cm}
\includegraphics[width=0.3\textwidth]{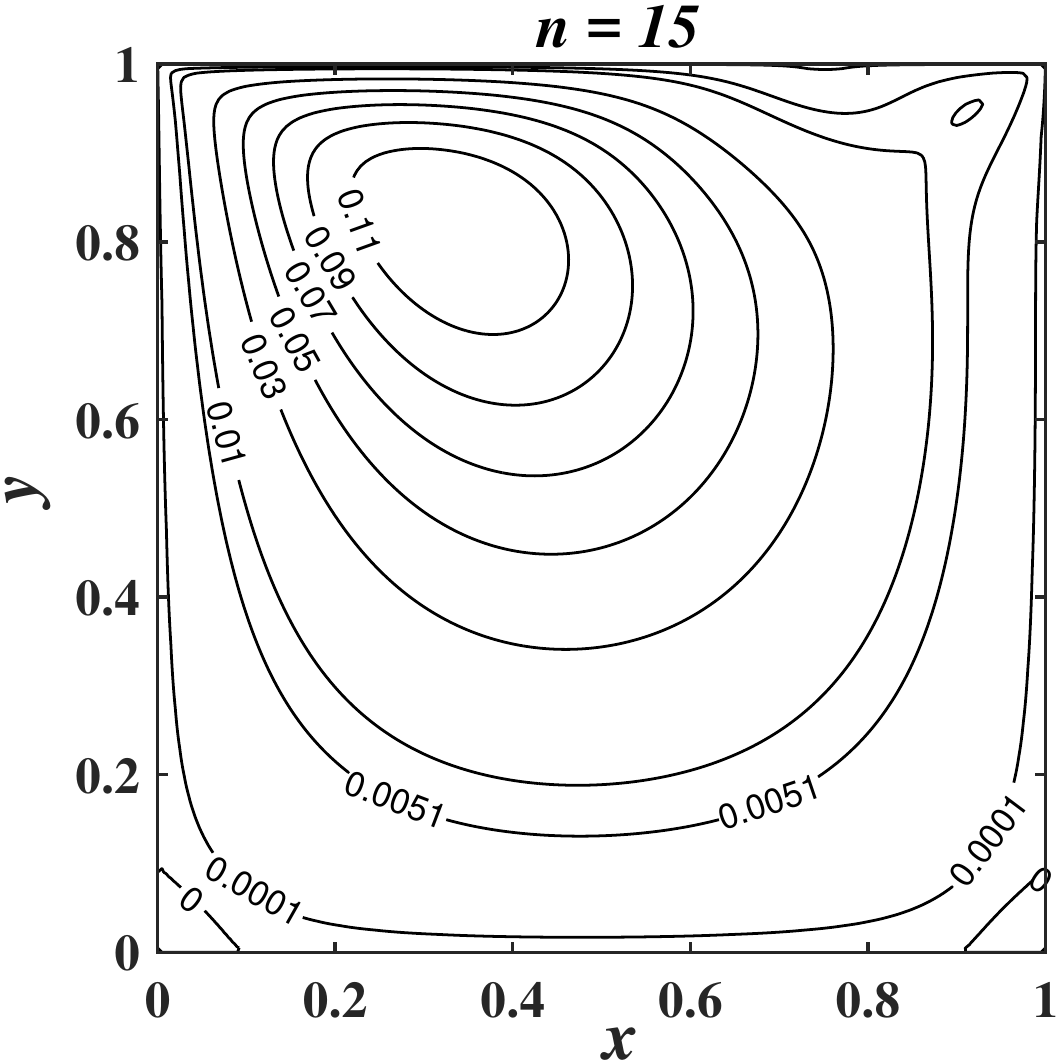} 
\\
\includegraphics[width=0.3\textwidth]{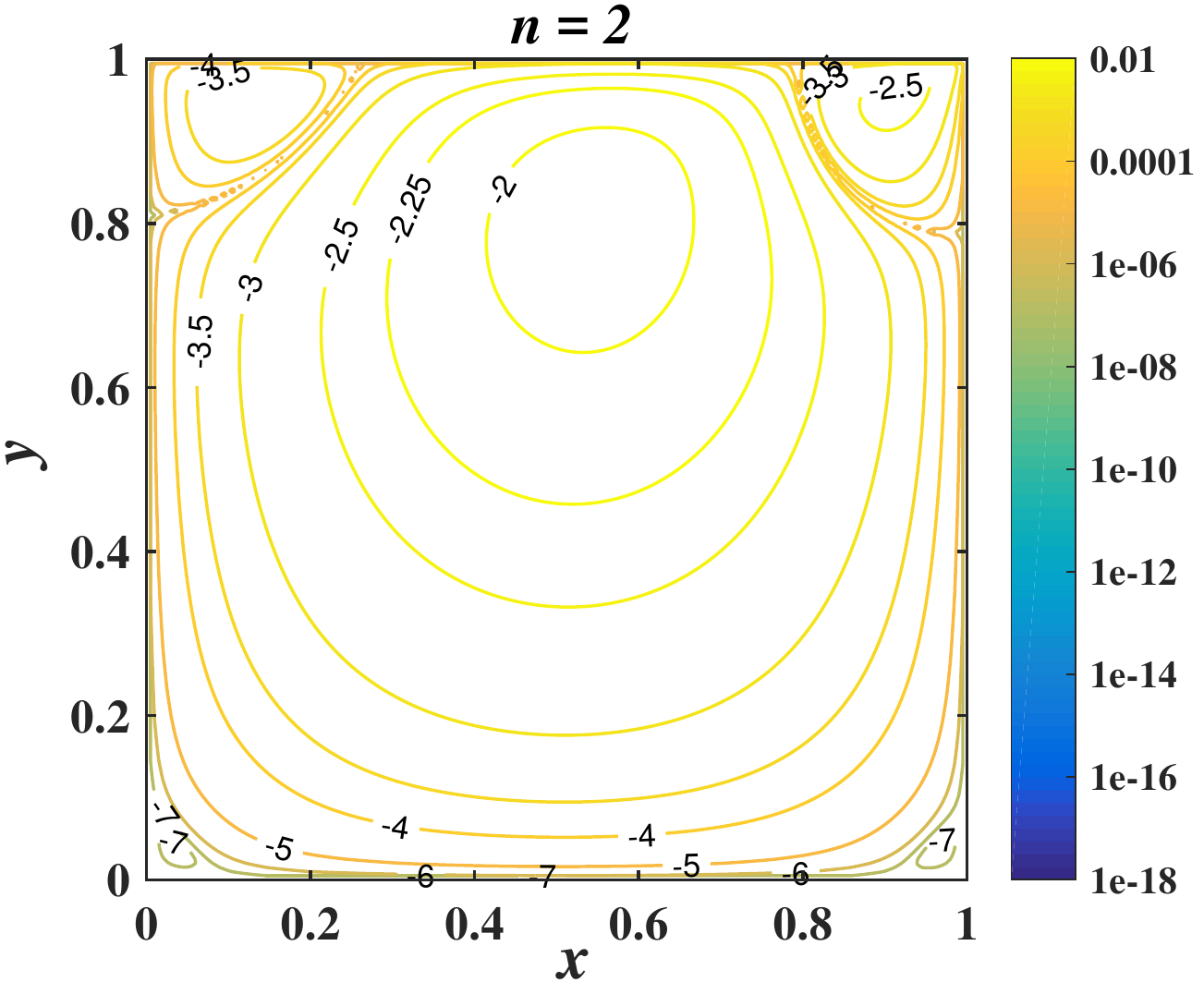}
\includegraphics[width=0.3\textwidth]{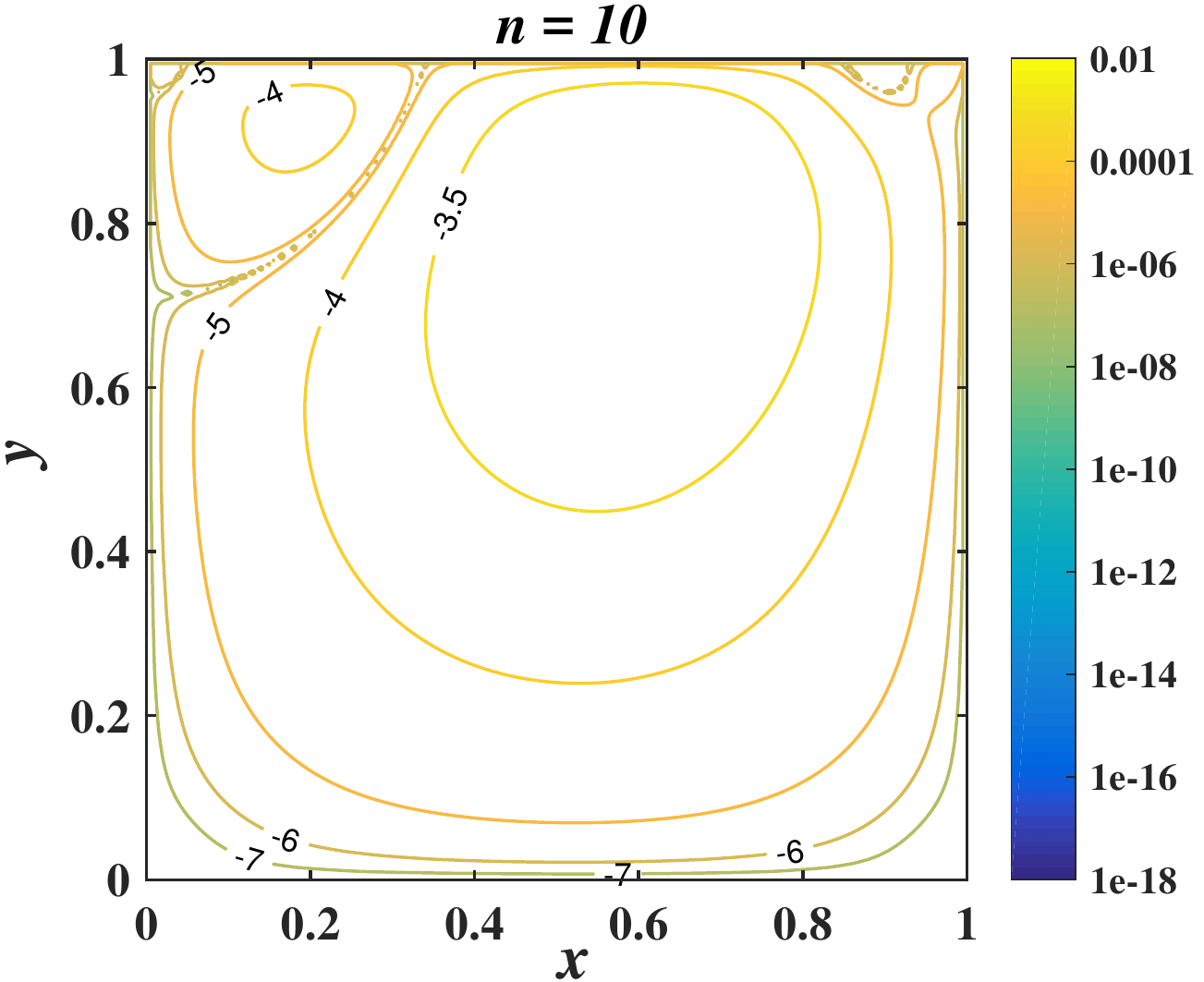}
\includegraphics[width=0.3\textwidth]{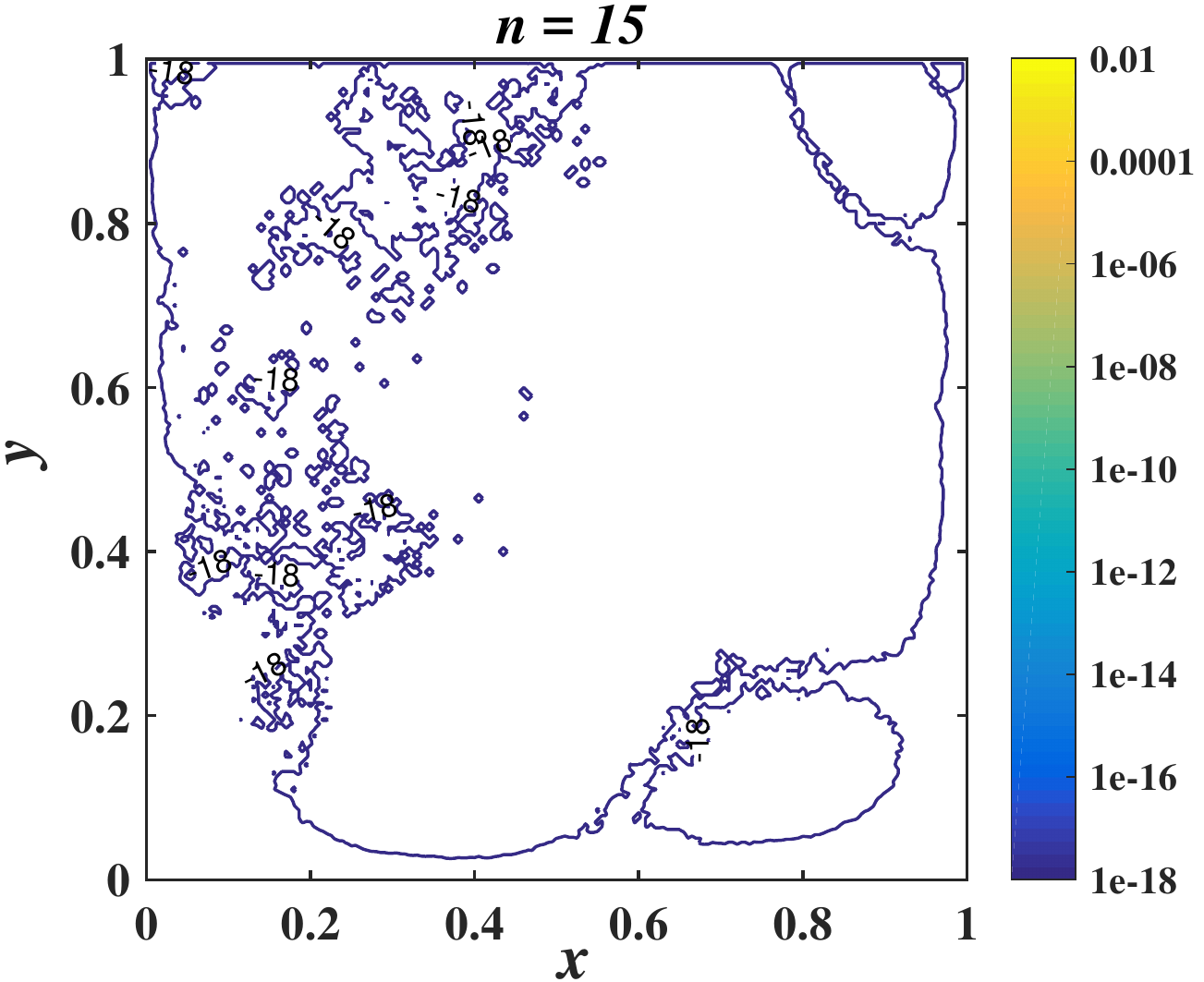}\\
\caption{\small Streamlines of the ROC approximations with {a given} RB dimension $n$ ({first row}), and the {corresponding} absolute streamline errors $E_q(n)$ in logarithmic scales ({second row}) for Stokes. }
\label{fig:stokes:streamline} 
\end{figure}
Specifically, the error curves of Stokes decrease much faster and attain the accuracy of our high fidelity solution with a small number of basis. This is consistent with other {results in the literature}, such as a parameter dependent Stokes problem {studied in} \cite{MadayMulaPateraYano2015}. 
{The spatial domain therein} is $\Omega =[0,1] \times [0,1]$ and the two parameters vary in $[1,8]\times [1,8]$. 
Relative {test} errors for 
both pressure and velocity fields {decay exponentially while stagnating} at about $M=25$ interpolatory terms with the accuracy of the underlying high fidelity solution {being} $10^{-5}\sim 10^{-6}$. 
Authors in \cite{diez2017generalized} use proper generalized decomposition (PGD) method with $40$ PGD terms to approximate solutions of the Stokes flow defined on a parametric domain. The final relative error 
stagnates at $10^{-6}$ for the velocity and $10^{-4}$ for the pressure.
For the nonlinear Navier-Stokes problem with a two-dimensional parameter, the error curve decrease{, as expected,} slower compared to the Stokes problem above. This is also {in line with} other existing literature {for steady-state Navier-Stokes equations with similar} settings. {For example, authors of \cite{ballarin2015supremizer}} approximate the steady-state incompressible NS equations with low Reynolds number by the POD-Galerkin method. 
{Their surrogate solutions for} some representative parameters $\text{Re} =30, 90$ attain $10^{-3} \sim 10^{-4}$ accuracy with about $30\sim 40$ basis {before stagnating}.

\begin{figure}[thbp]
\centering
%\hspace{-0.5cm} 
\includegraphics[width=0.3\textwidth]{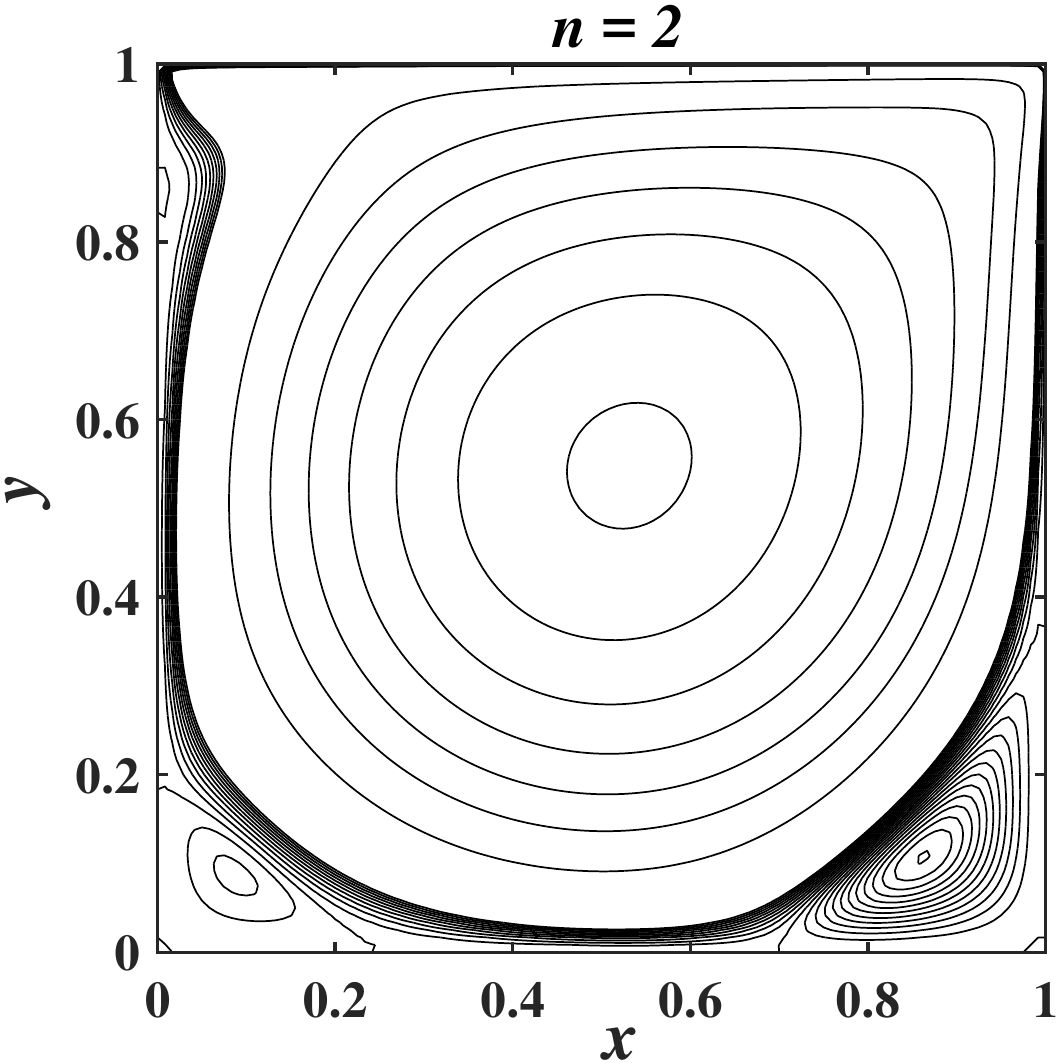} \hspace{0.4cm}
\includegraphics[width=0.3\textwidth]{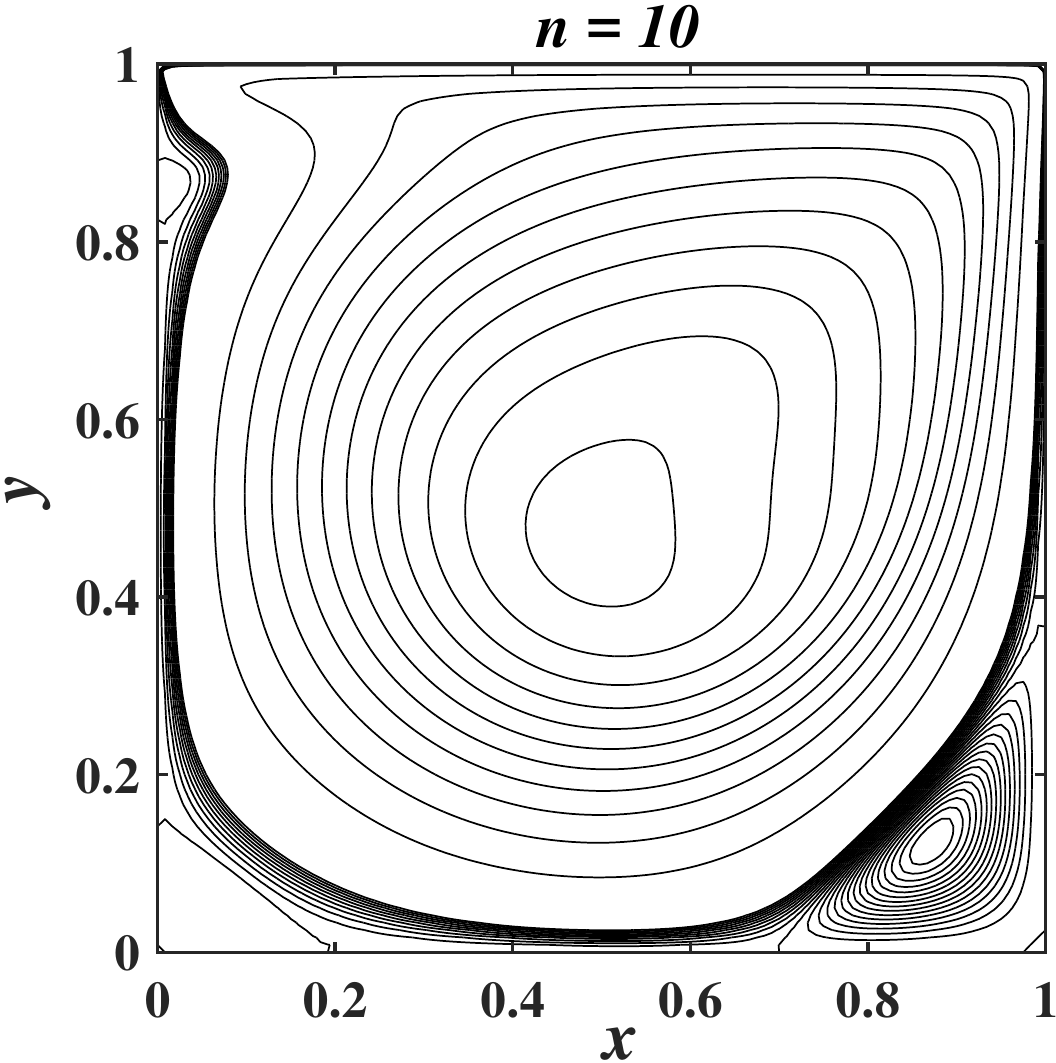} \hspace{0.4cm}
\includegraphics[width=0.3\textwidth]{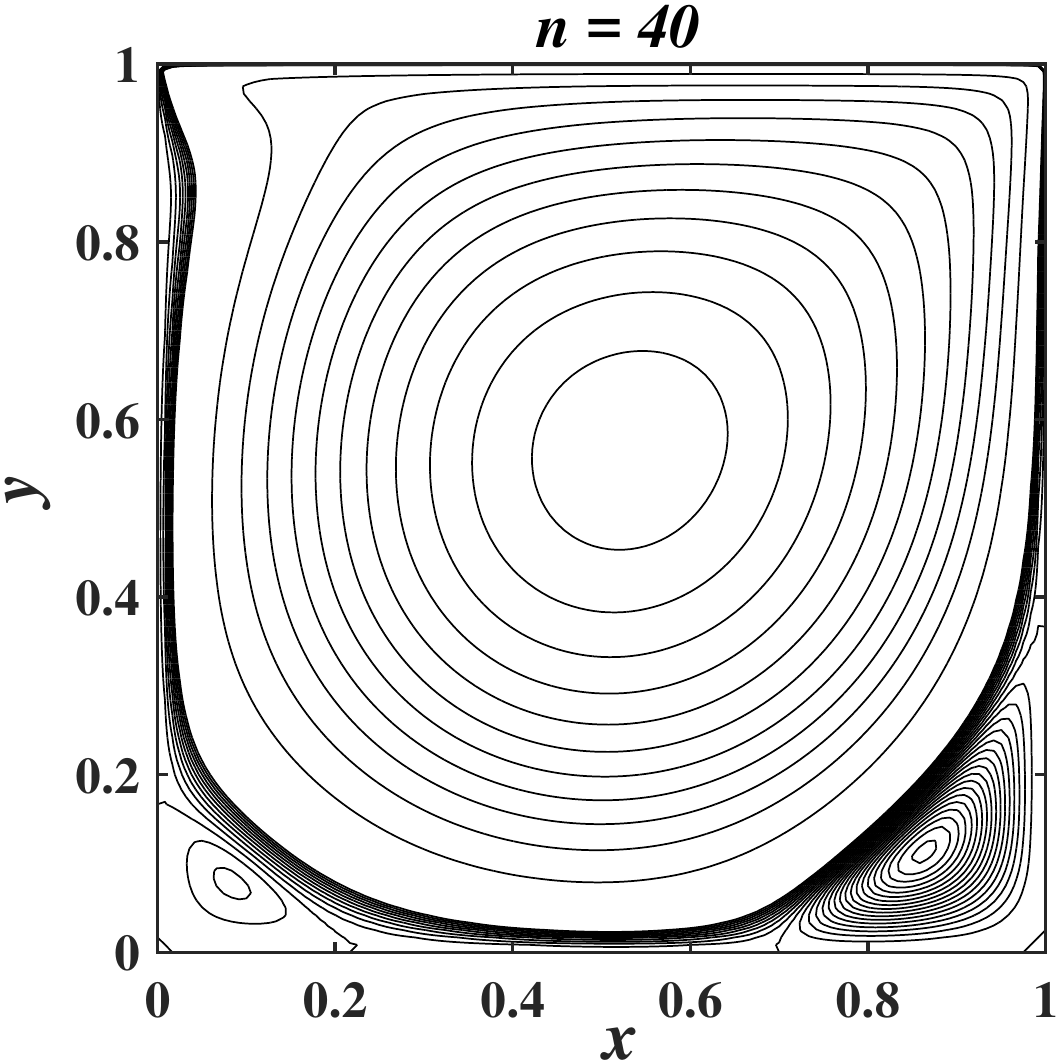}\\
\includegraphics[width=0.32\textwidth]{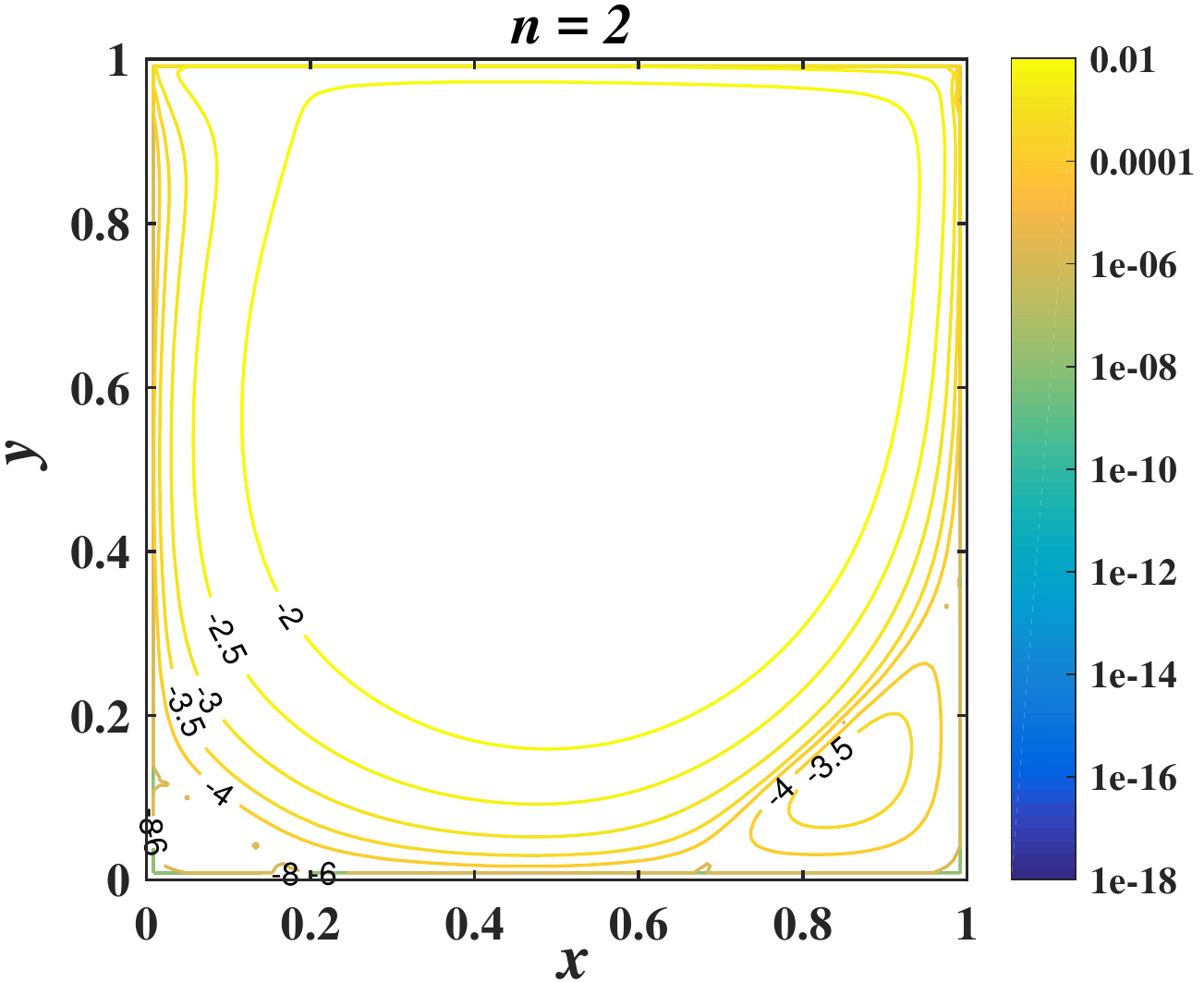}
\includegraphics[width=0.32\textwidth]{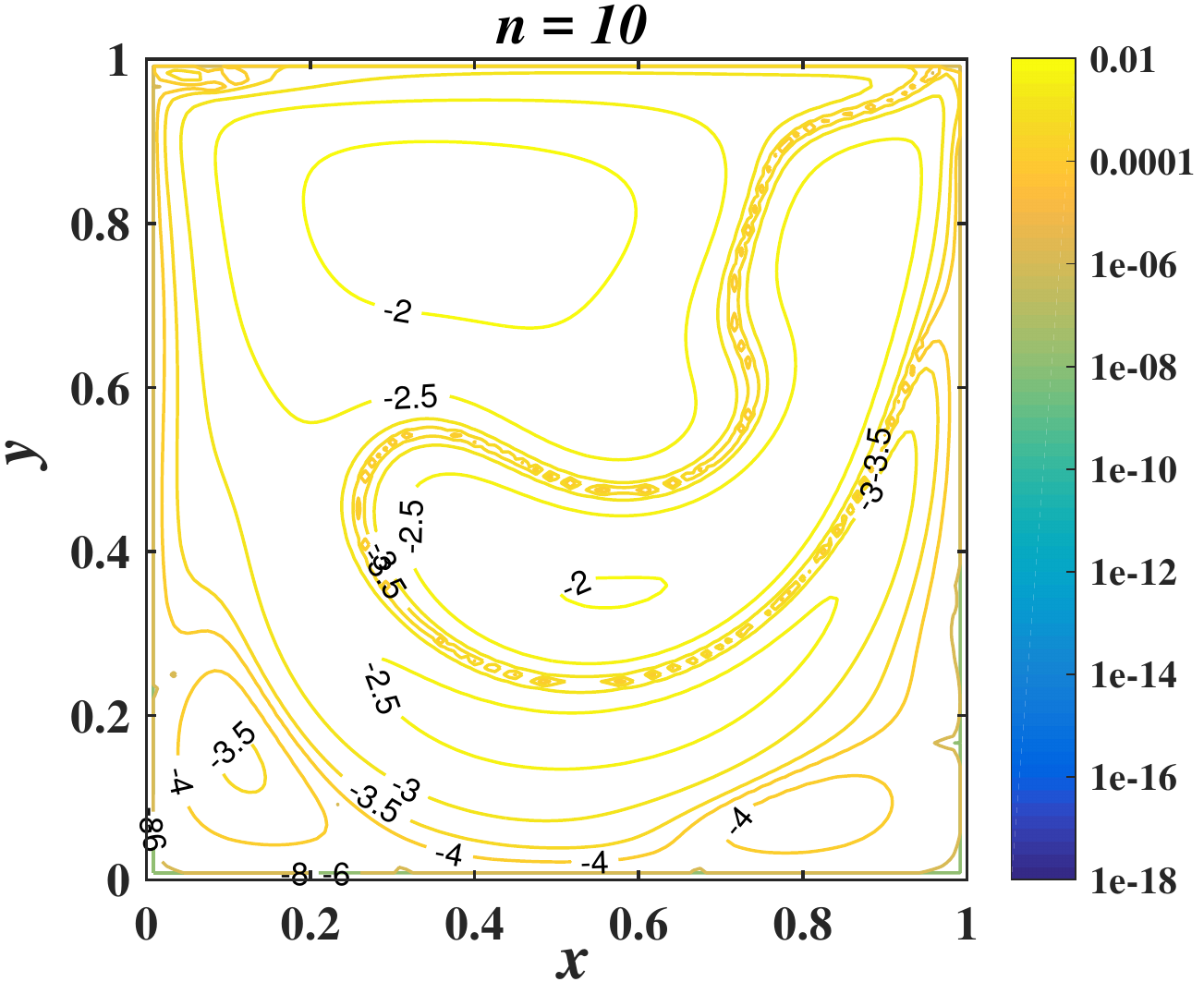}
\includegraphics[width=0.32\textwidth]{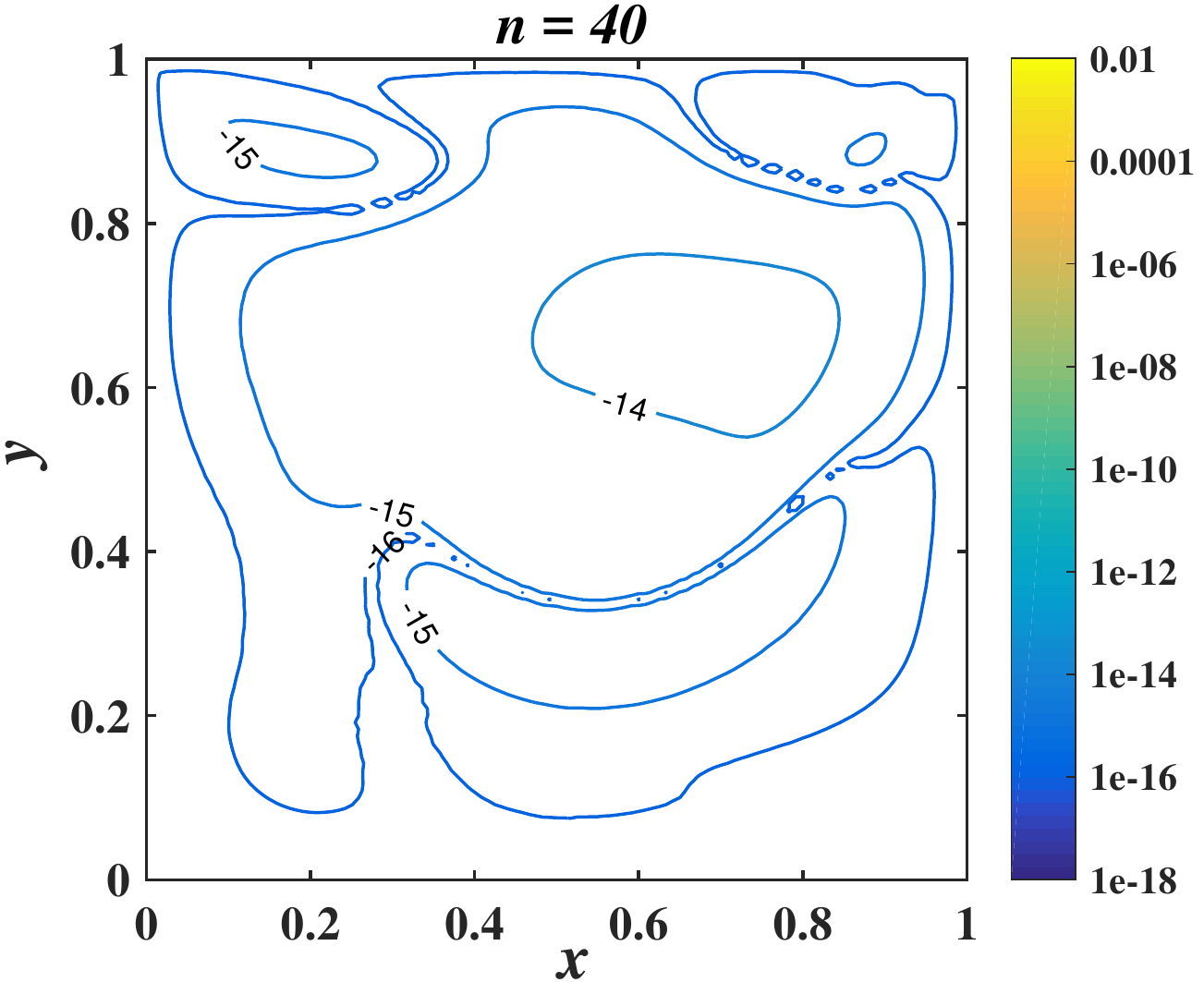}\\
\hspace{-0.5cm} 
\includegraphics[width=0.3\textwidth]{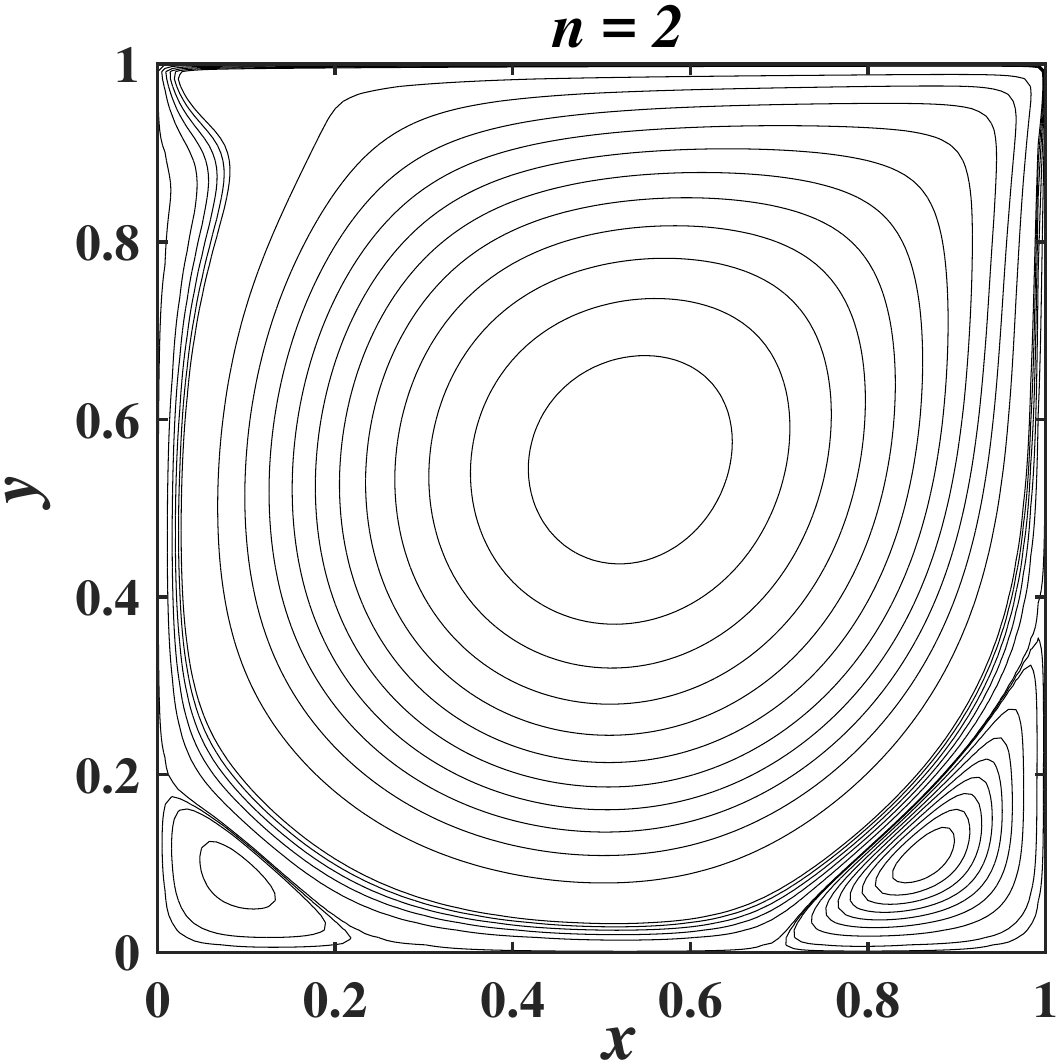} \hspace{0.4cm} 
\includegraphics[width=0.3\textwidth]{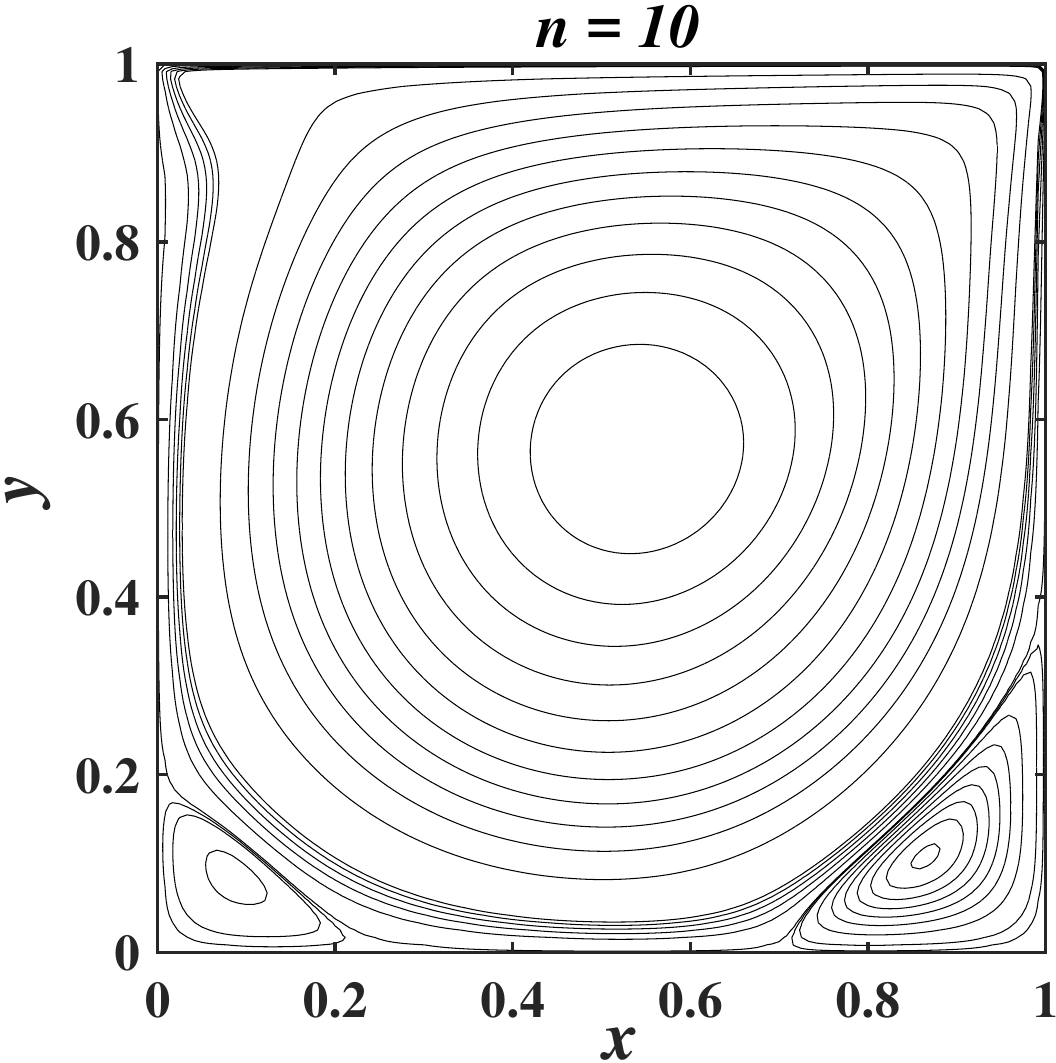} \hspace{0.4cm} 
\includegraphics[width=0.3\textwidth]{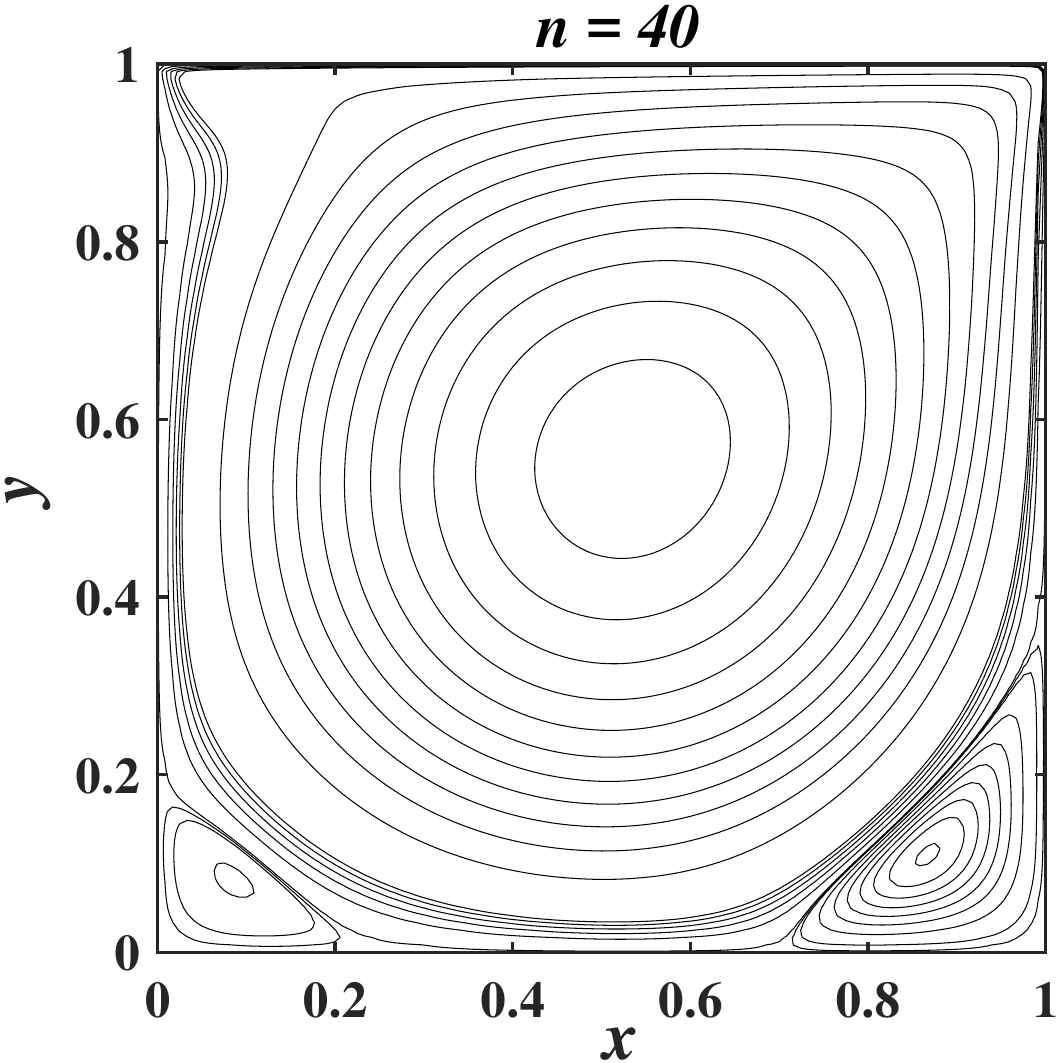}\\
\includegraphics[width=0.32\textwidth]{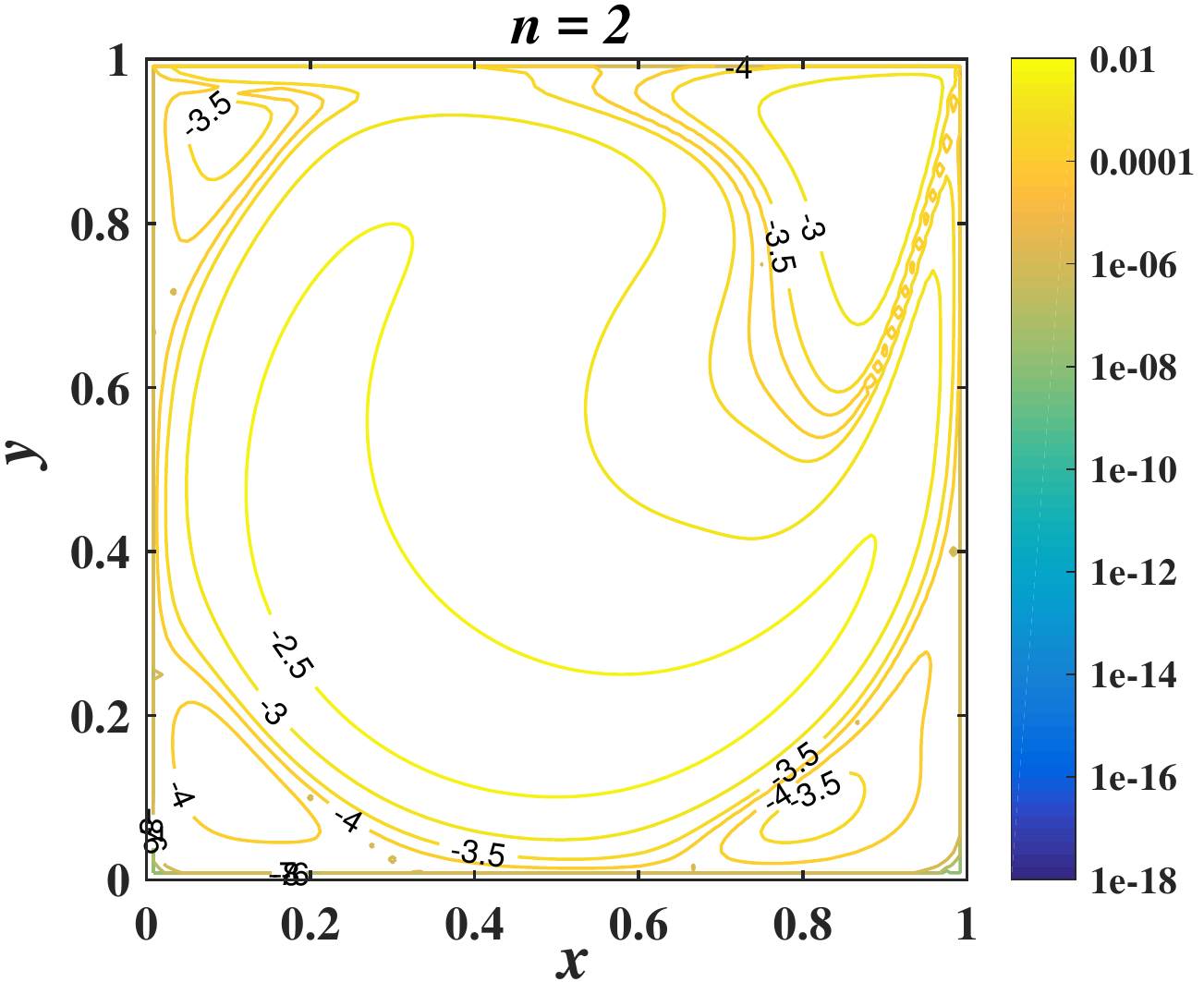}
\includegraphics[width=0.32\textwidth]{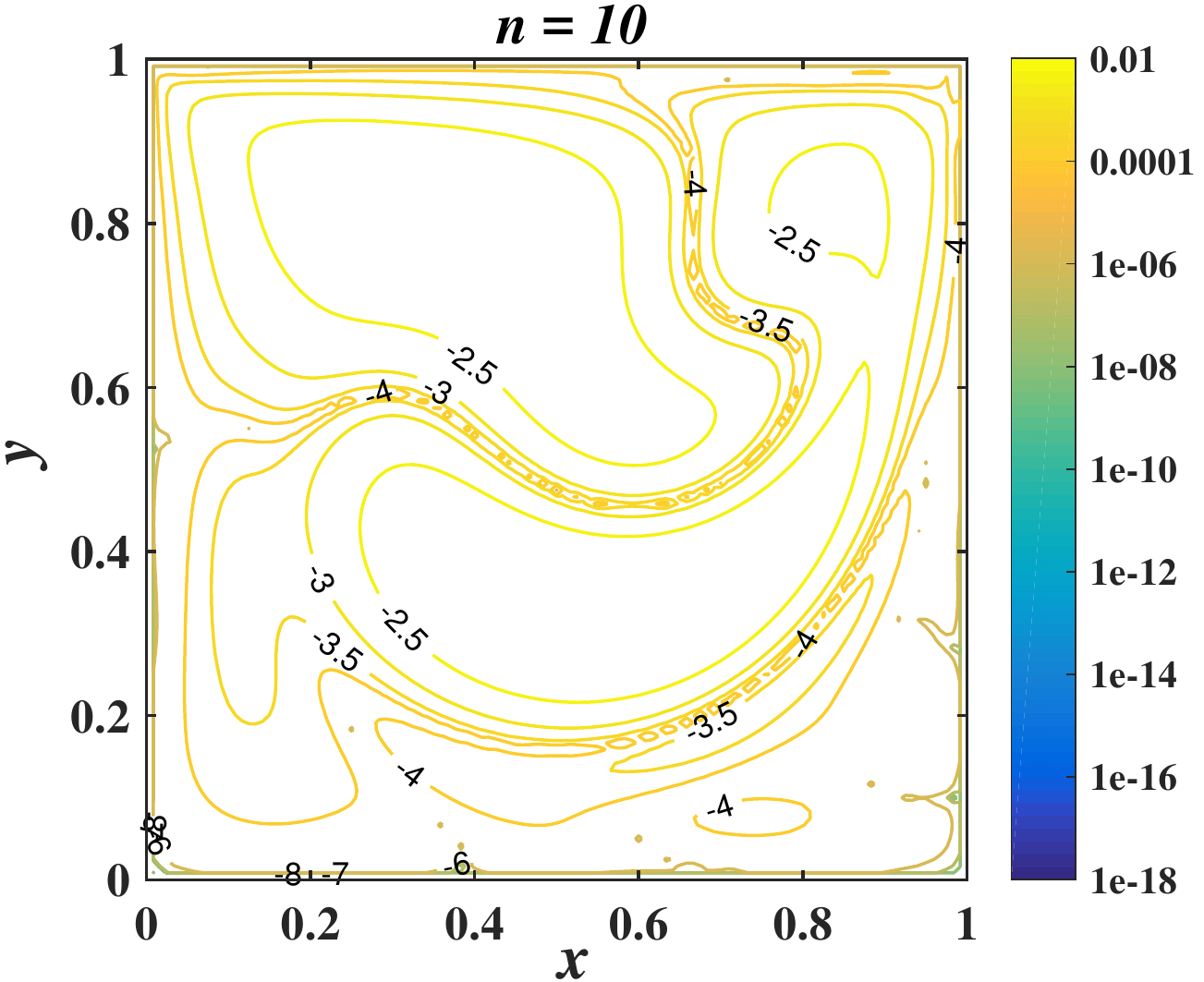}
\includegraphics[width=0.32\textwidth]{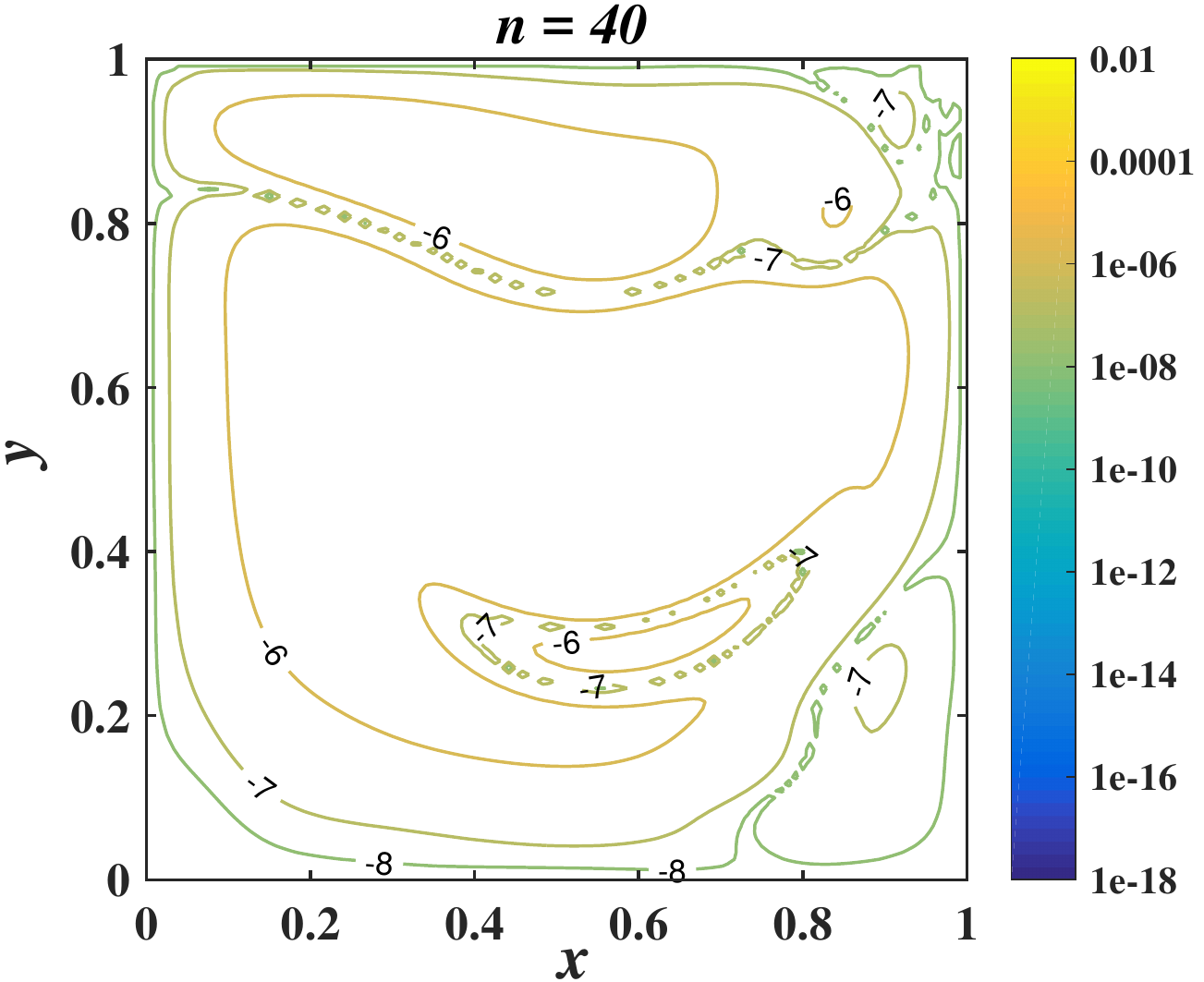}
\caption{Streamlines of the ROC approximations with {a given} RB dimension $n$ (Rows 1 and 3), and the absolute streamline errors $E_q(n)$ in logarithmic scale (Rows 2 and 4) for Navier-Stokes {with parameter values $(Re, \nu) = (1000, 0)$ on the top and $(500, 2)$ at the bottom}. }
\label{fig:Nstokes:streamline} 
\end{figure}
\begin{figure}[thbp]
\centering
\includegraphics[width=0.22\textwidth]{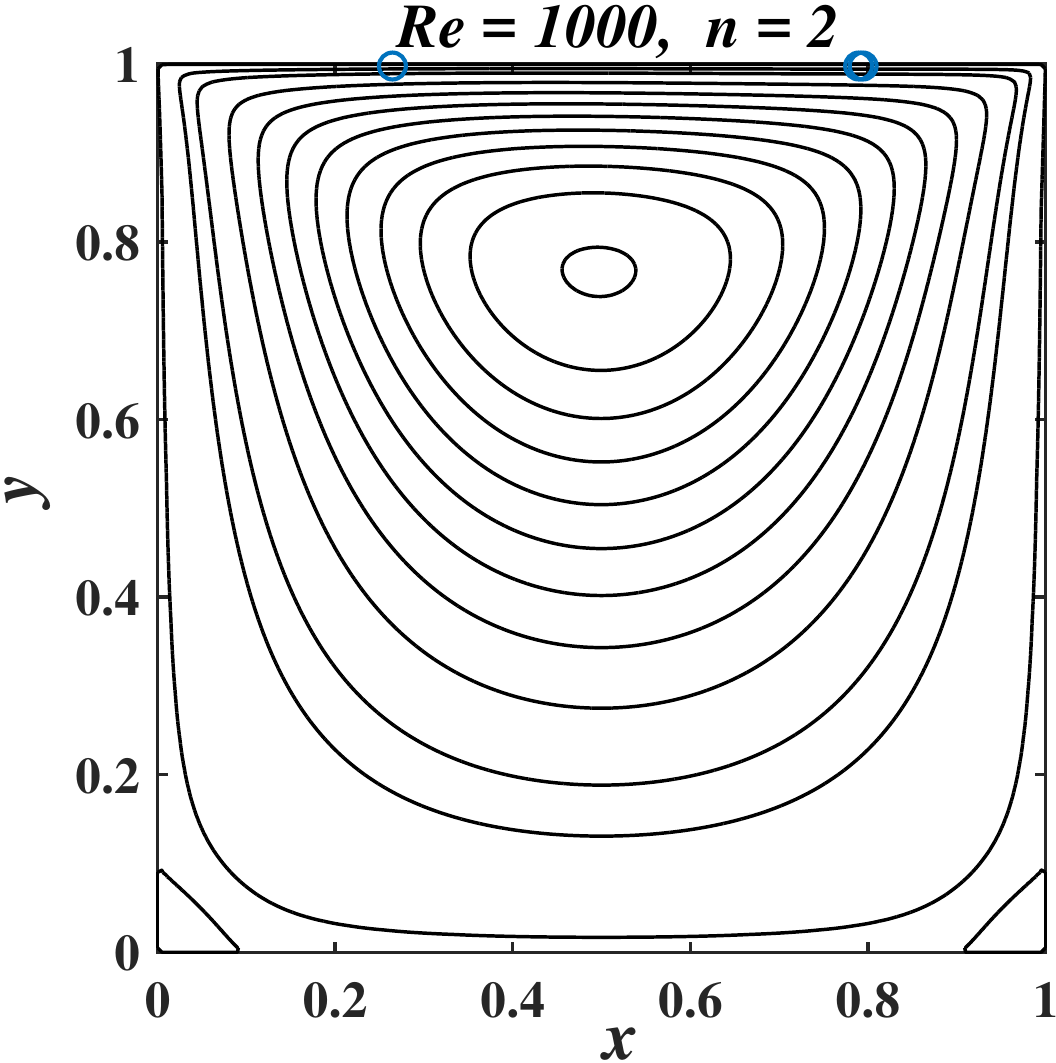}
\includegraphics[width=0.22\textwidth]{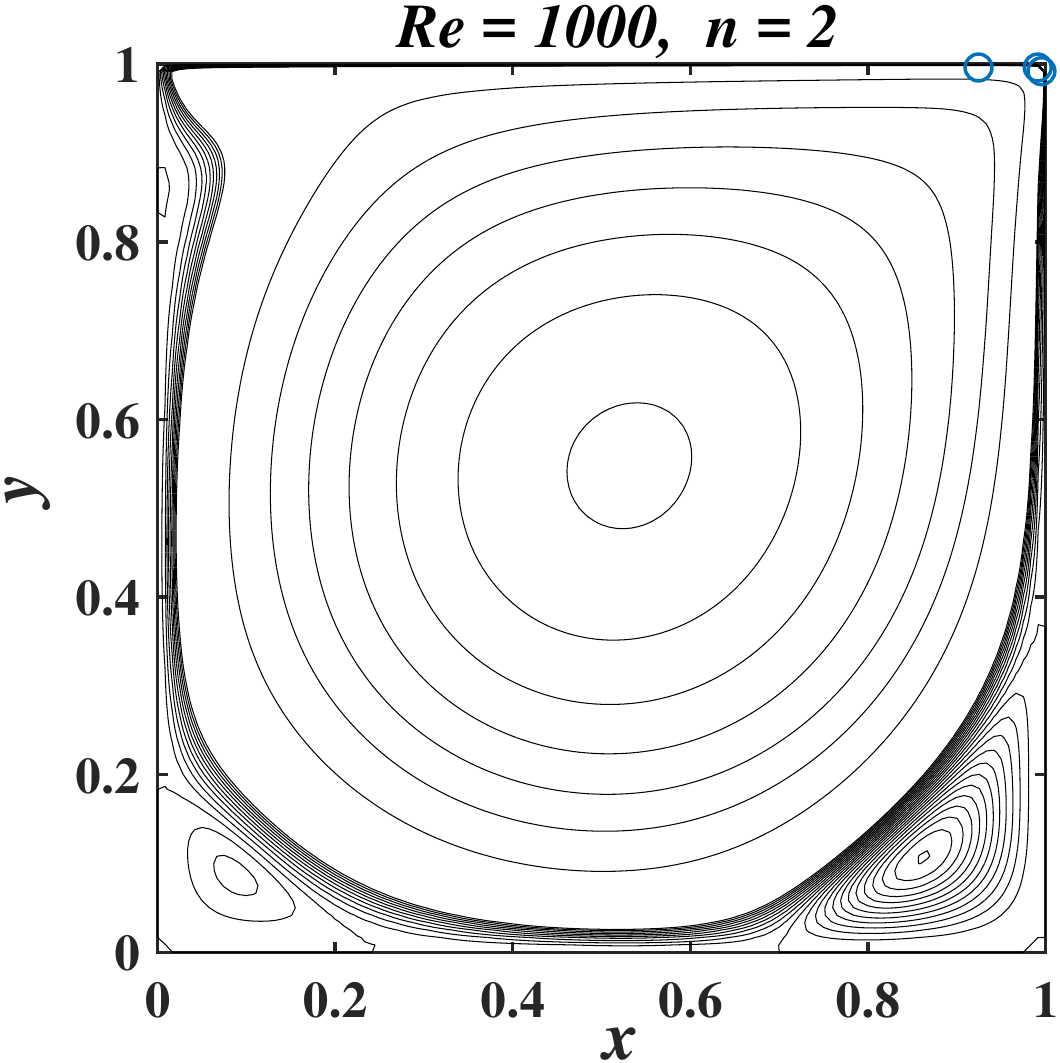}
\includegraphics[width=0.25\textwidth]{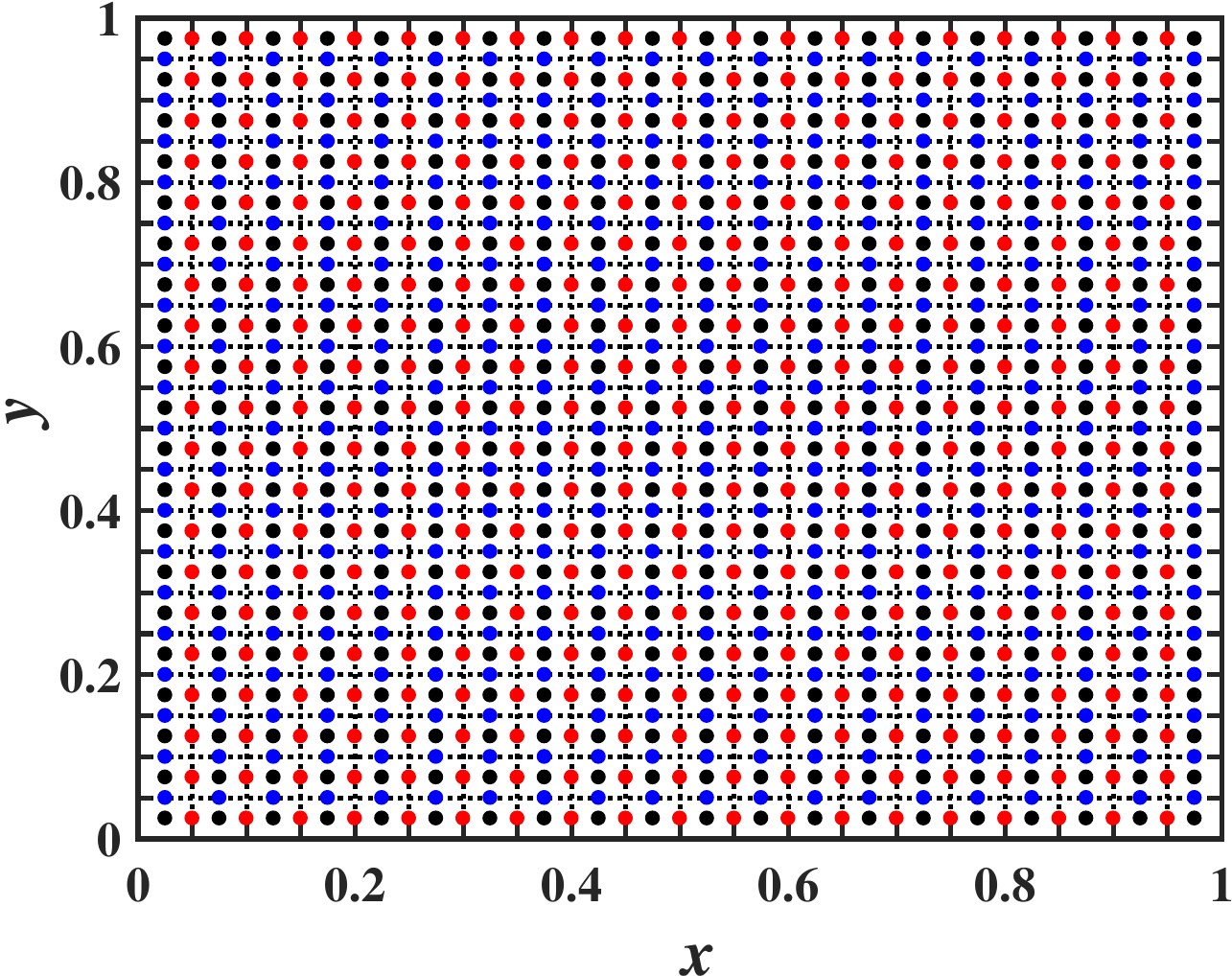}
\includegraphics[width=0.25\textwidth]{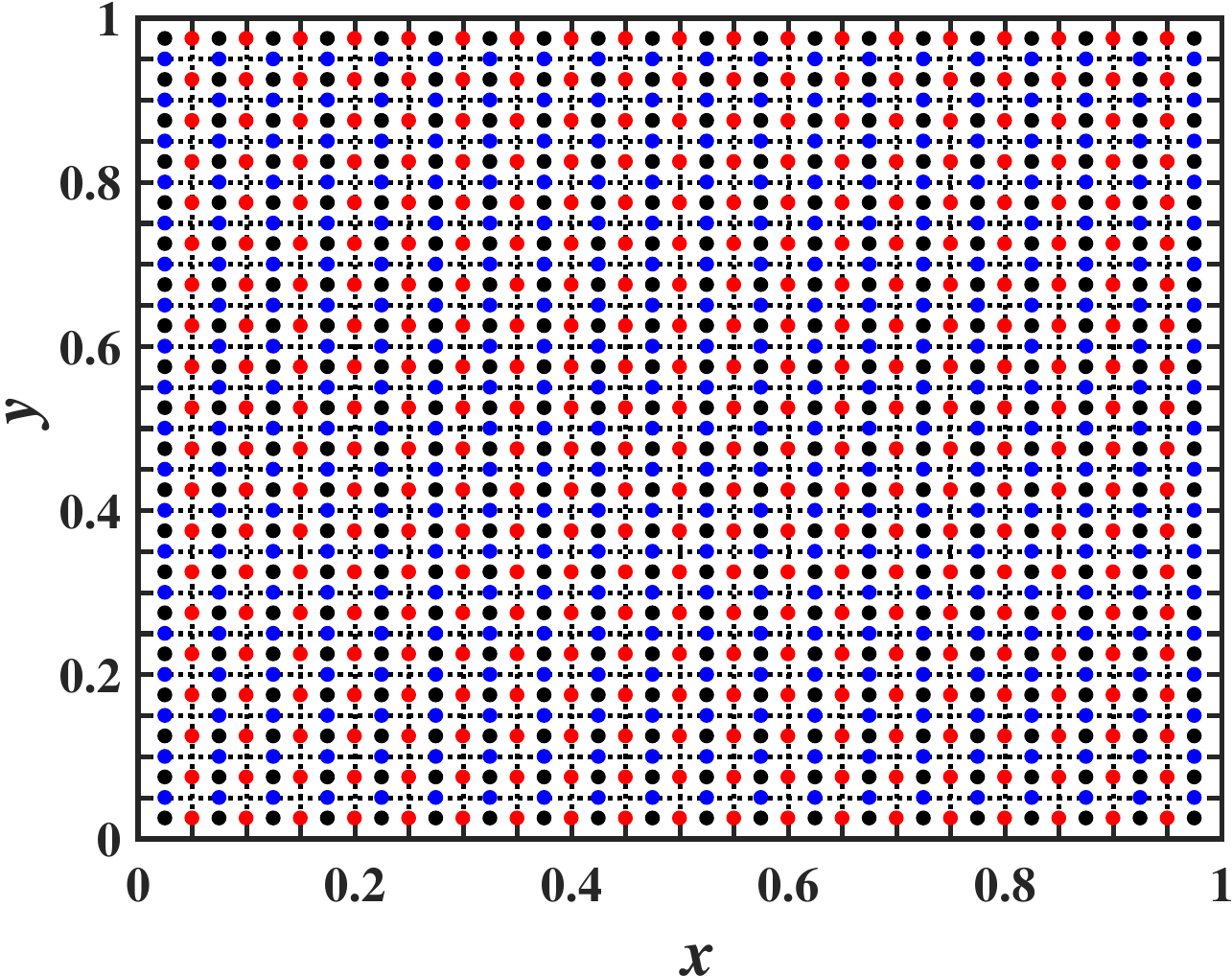}\\
\includegraphics[width=0.22\textwidth]{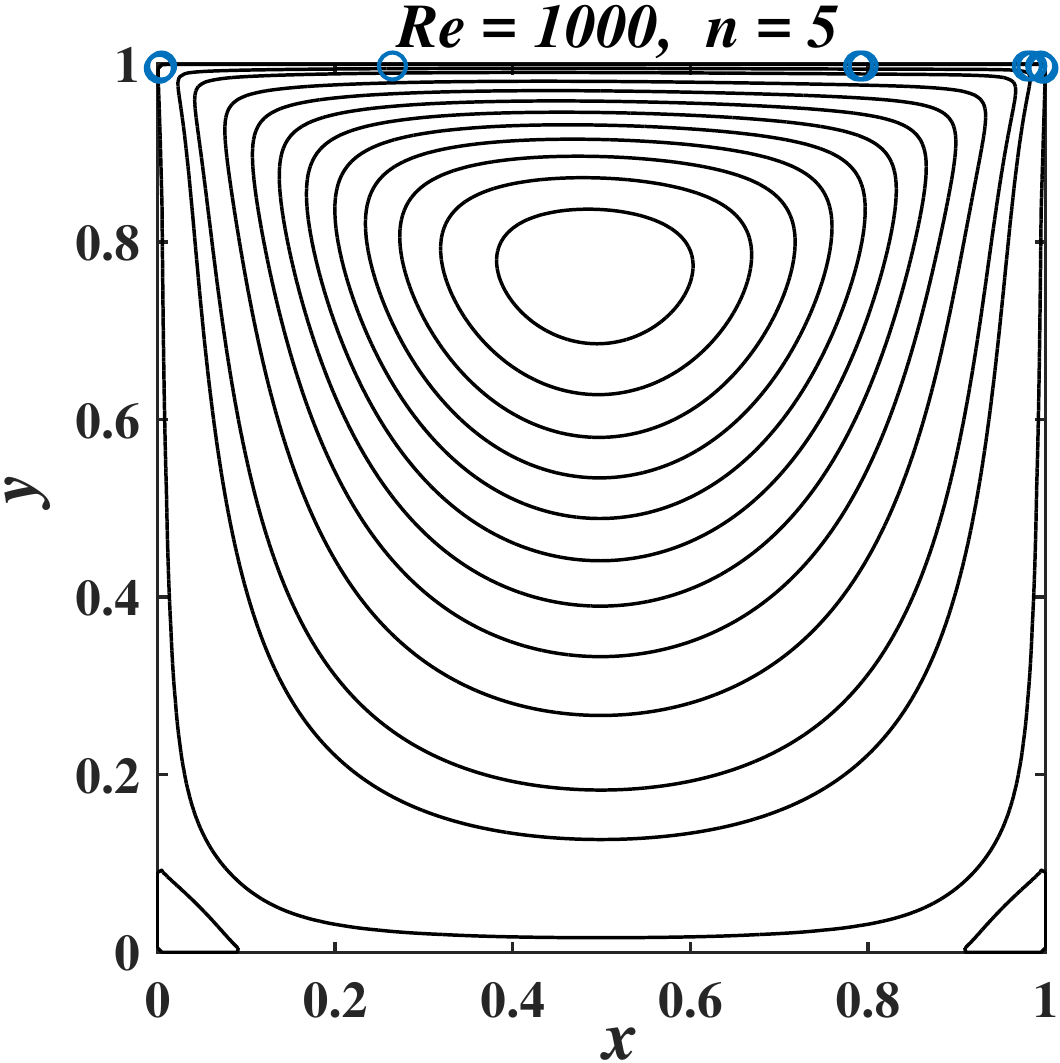}
\includegraphics[width=0.22\textwidth]{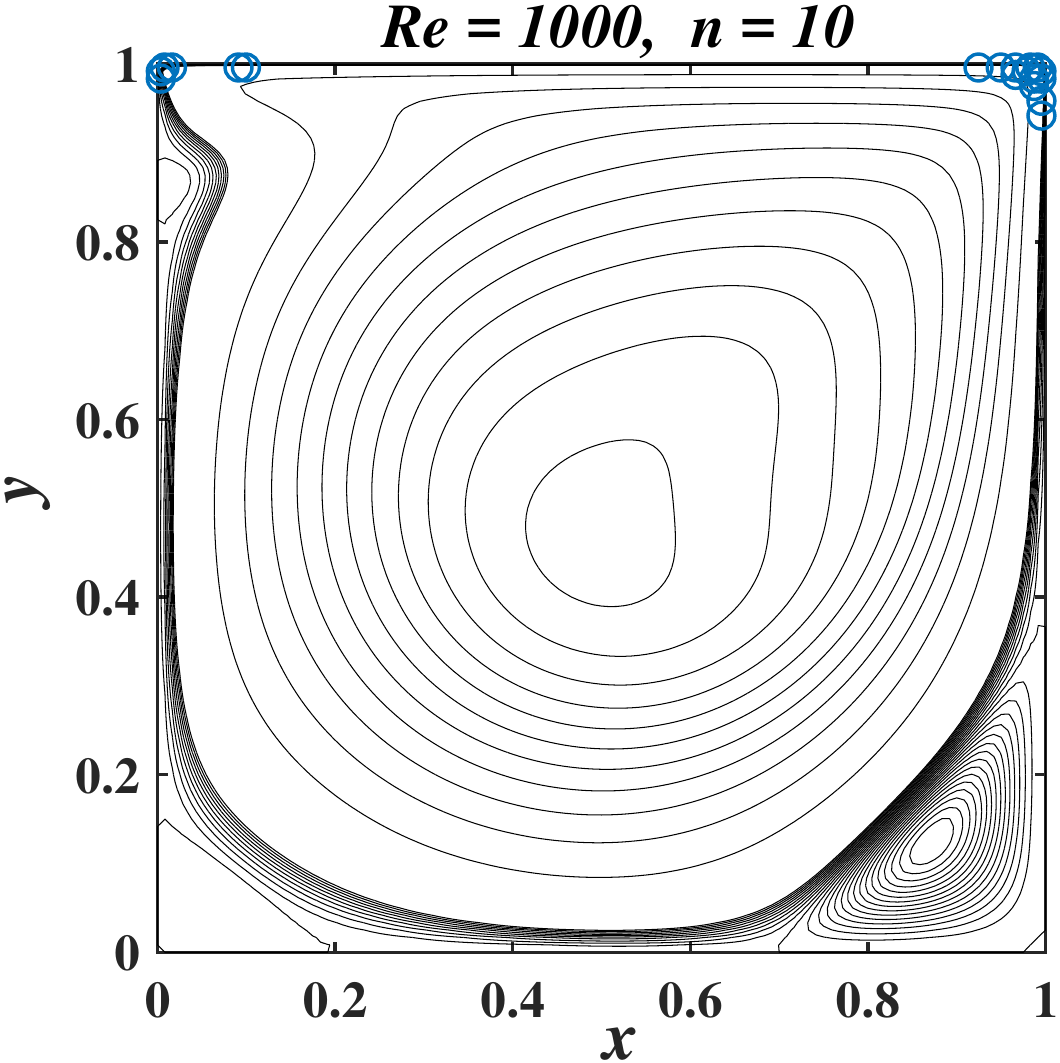}
\includegraphics[width=0.26\textwidth]{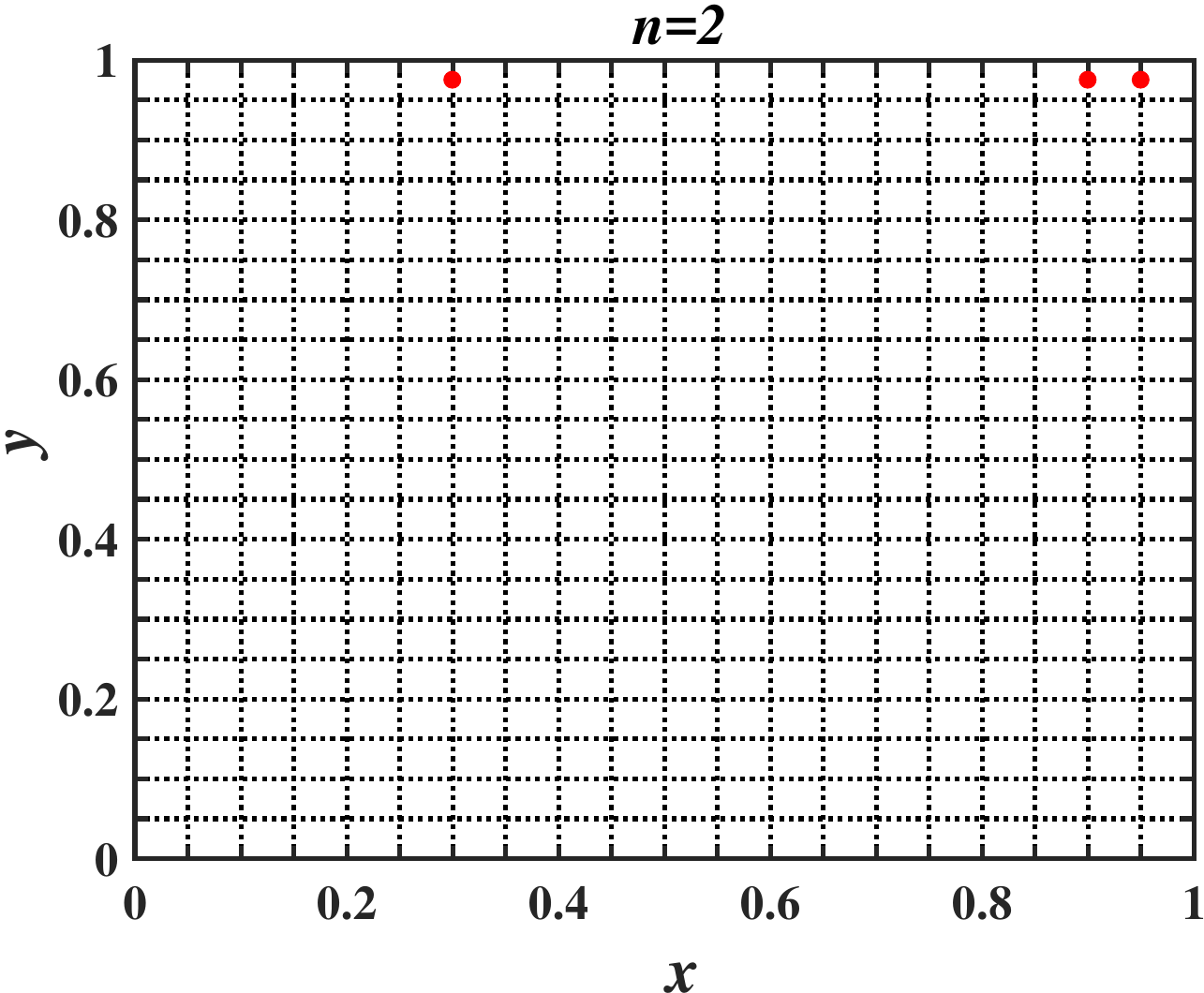}
\includegraphics[width=0.26\textwidth]{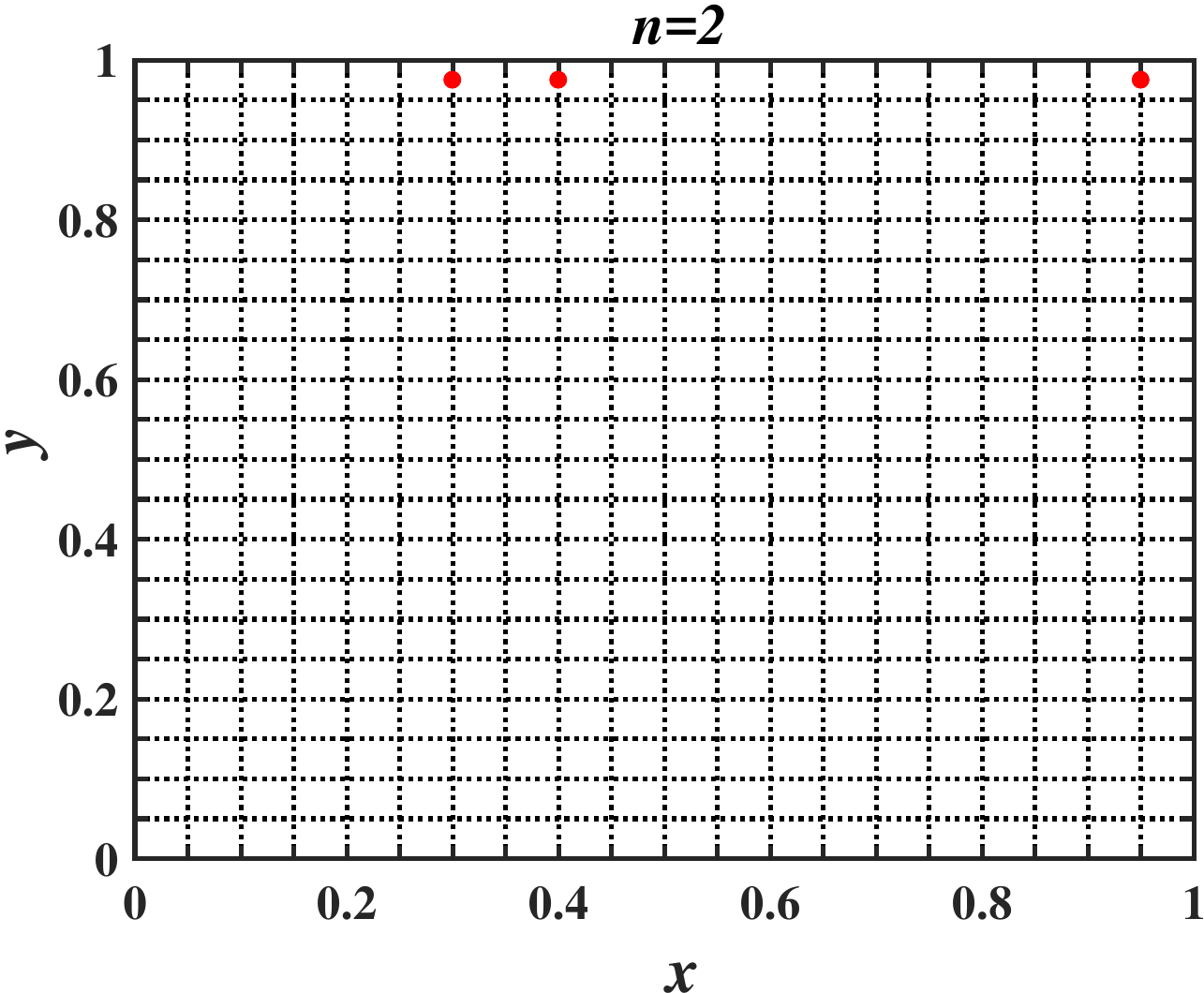}\\
\includegraphics[width=0.22\textwidth]{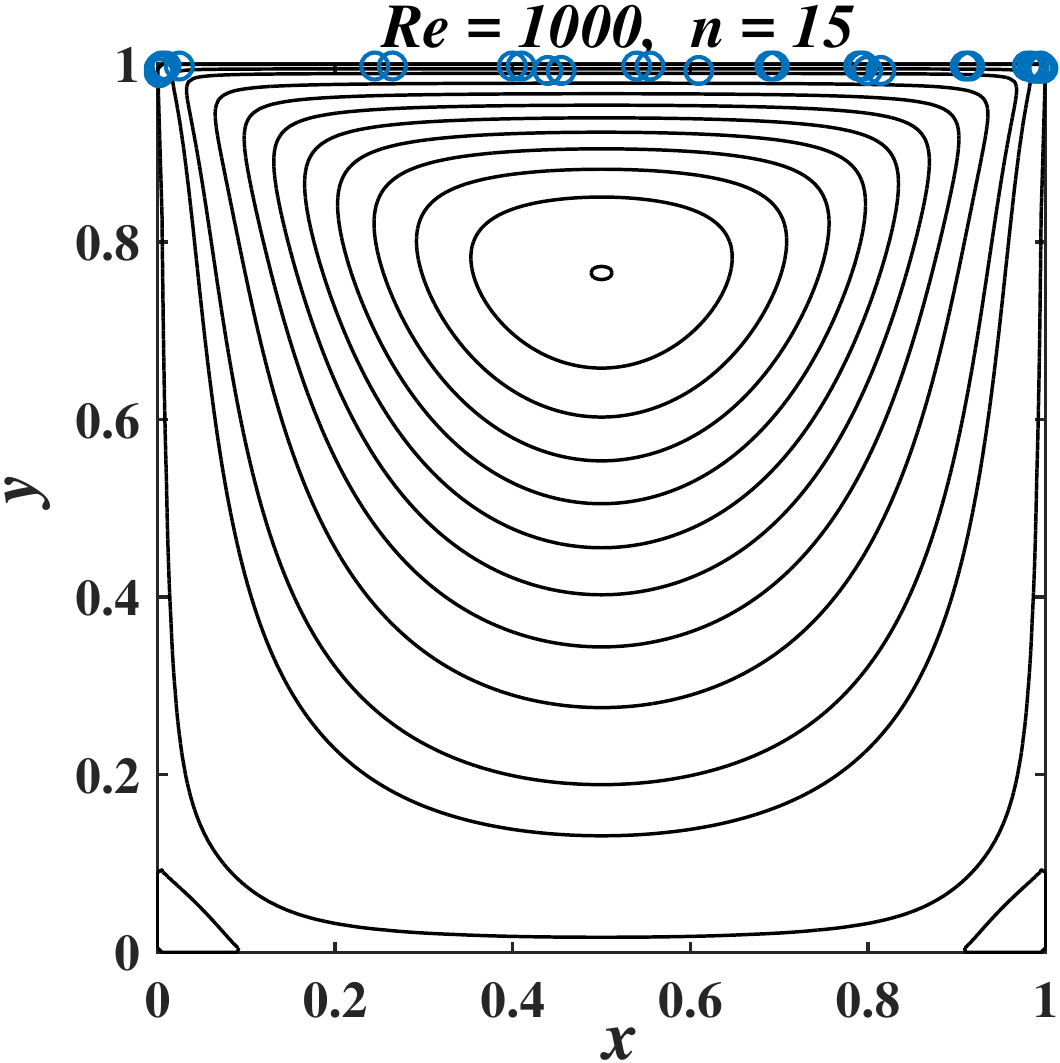}
\includegraphics[width=0.22\textwidth]{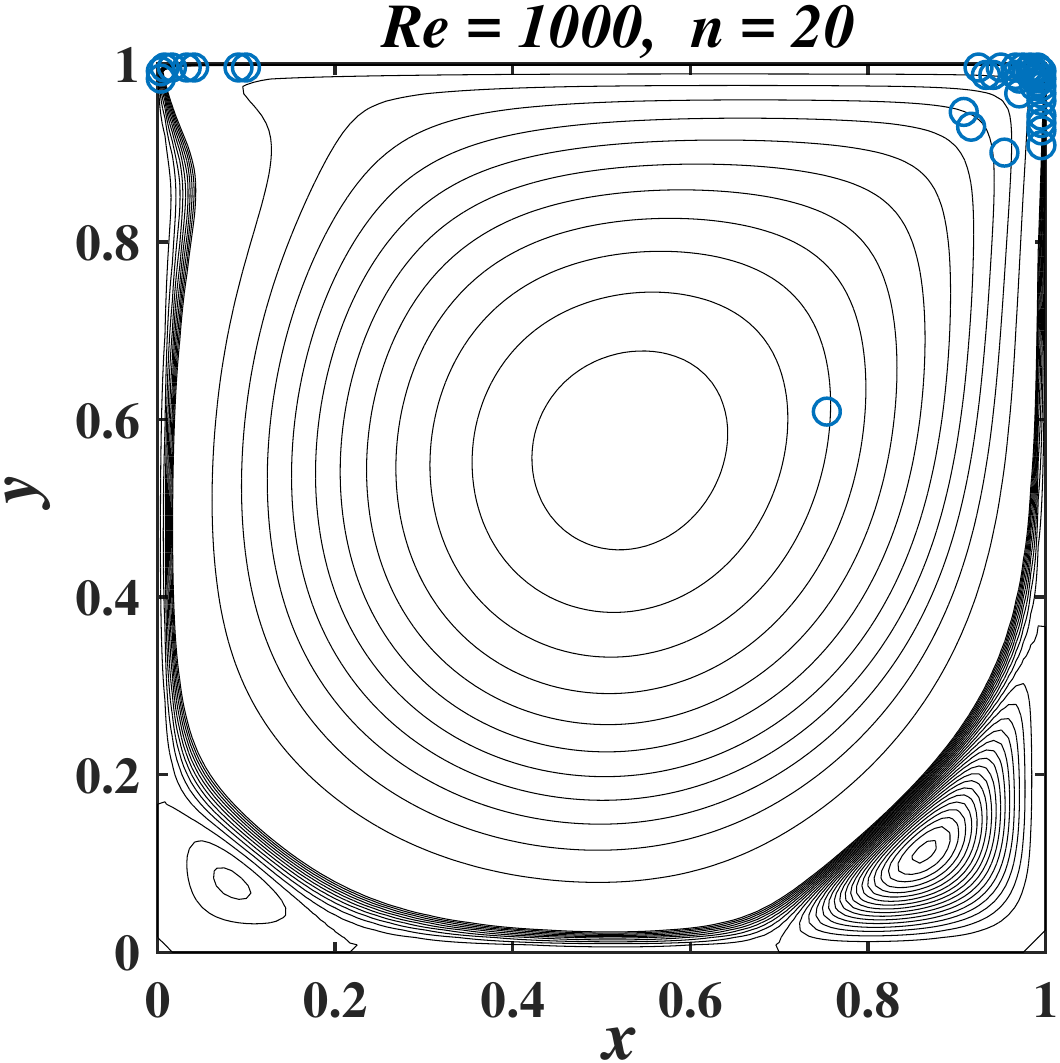}
\includegraphics[width=0.26\textwidth]{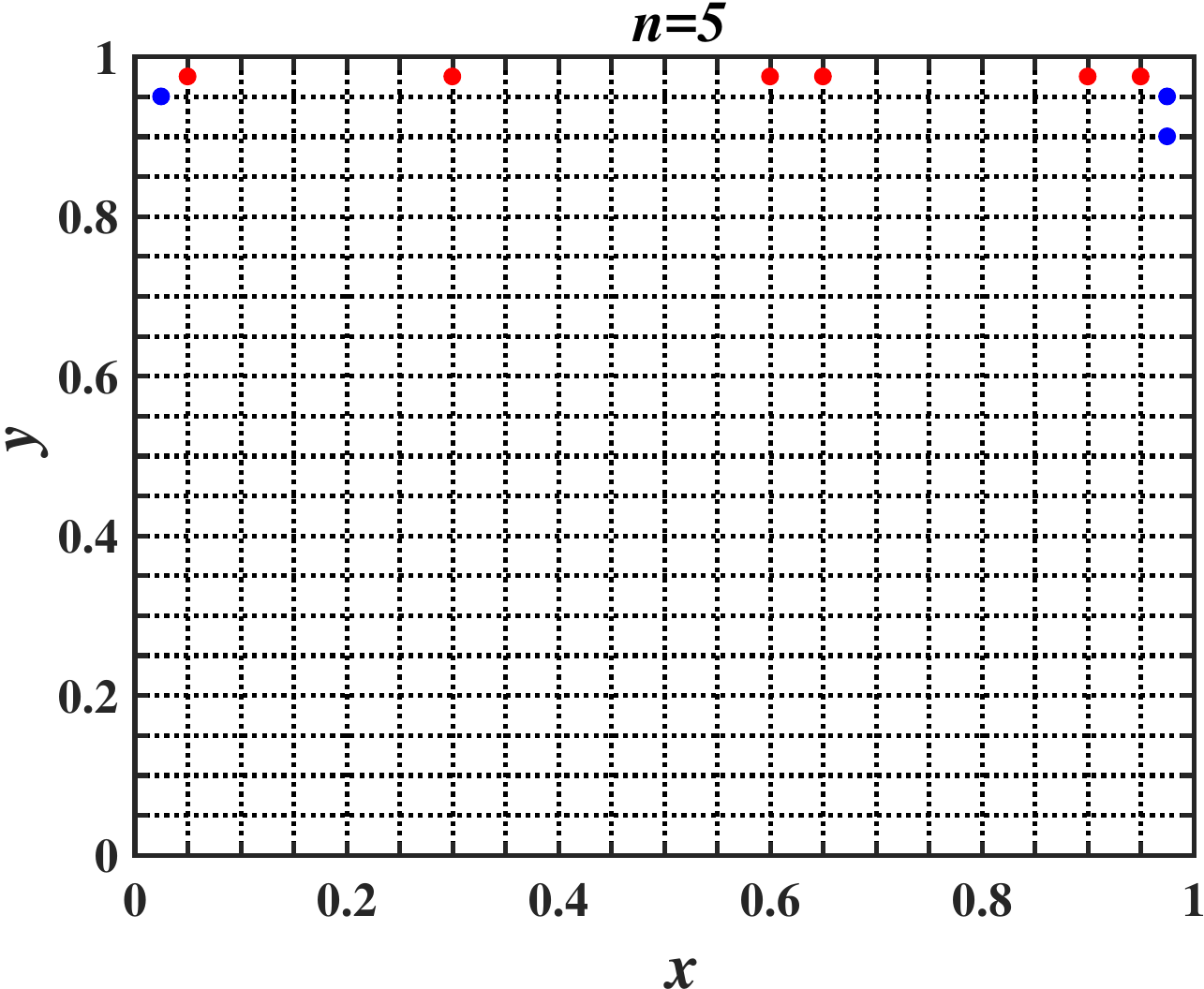}
\includegraphics[width=0.26\textwidth]{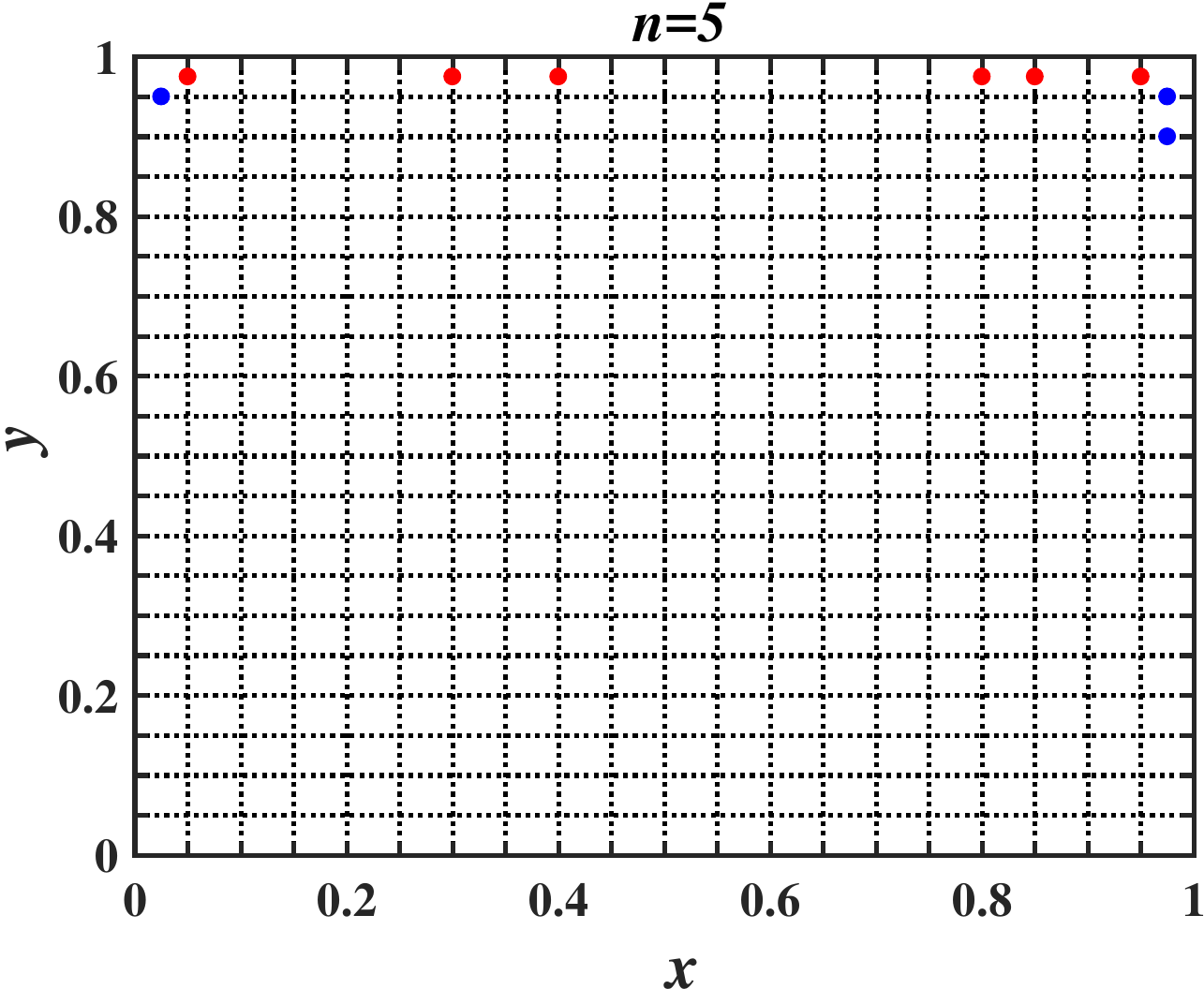}\\
\includegraphics[width=0.22\textwidth]{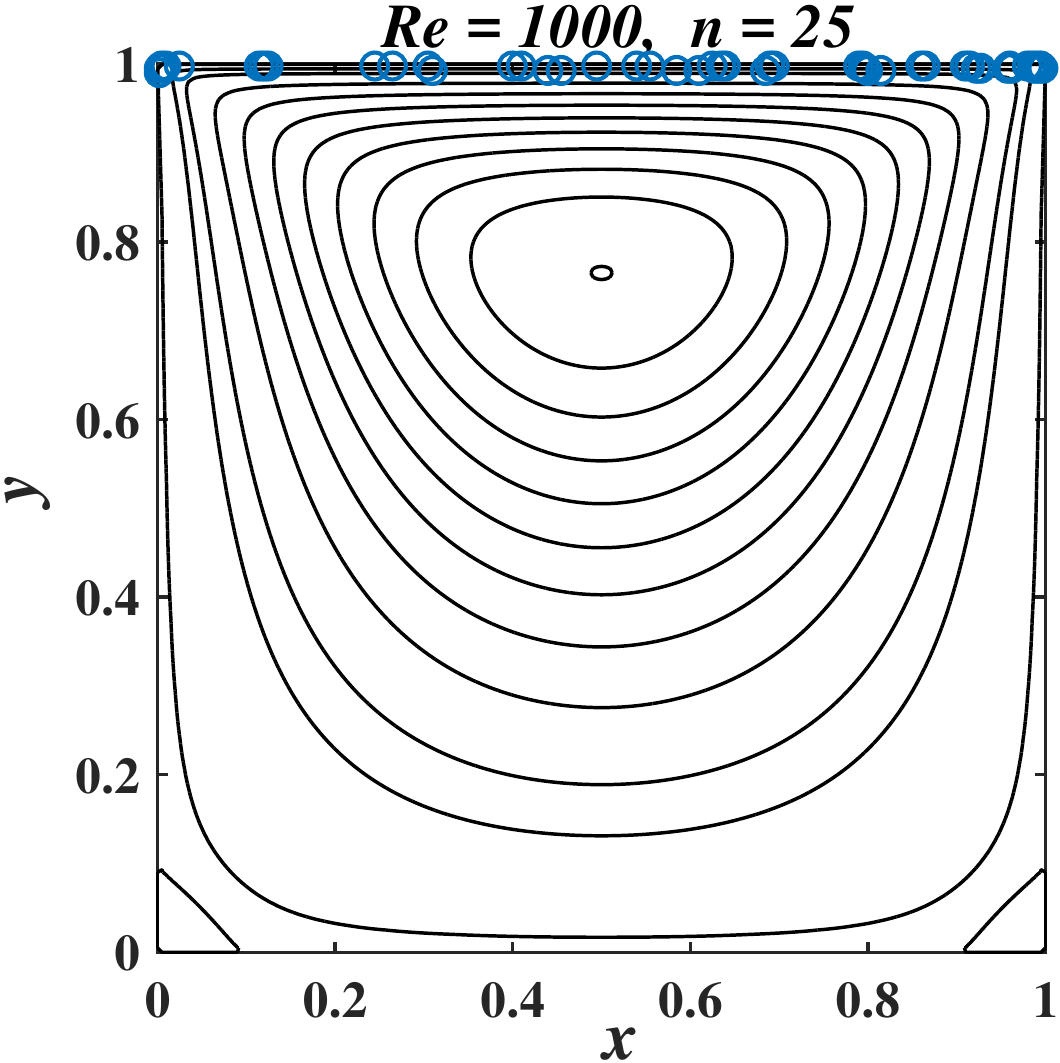}
\includegraphics[width=0.22\textwidth]{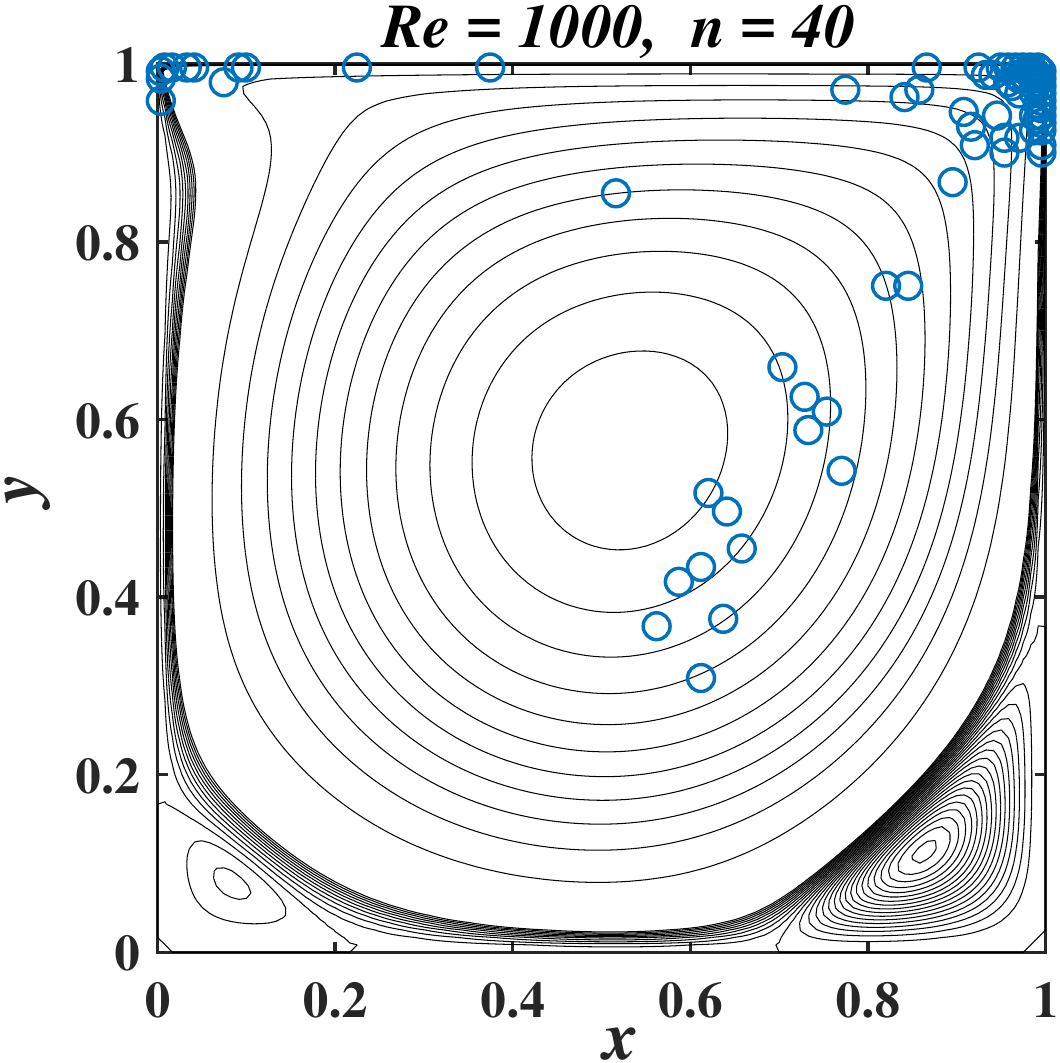}
\includegraphics[width=0.26\textwidth]{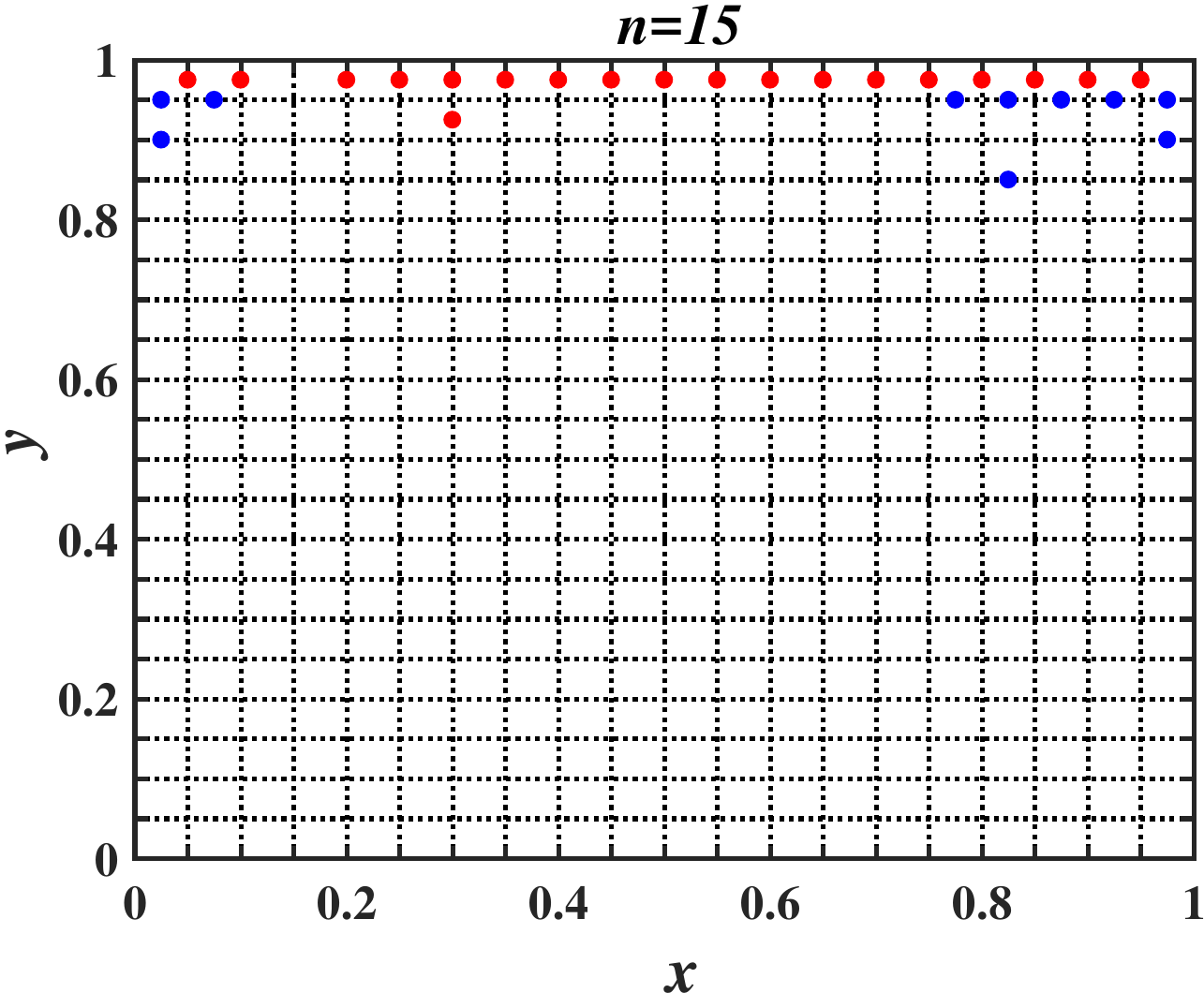}
\includegraphics[width=0.26\textwidth]{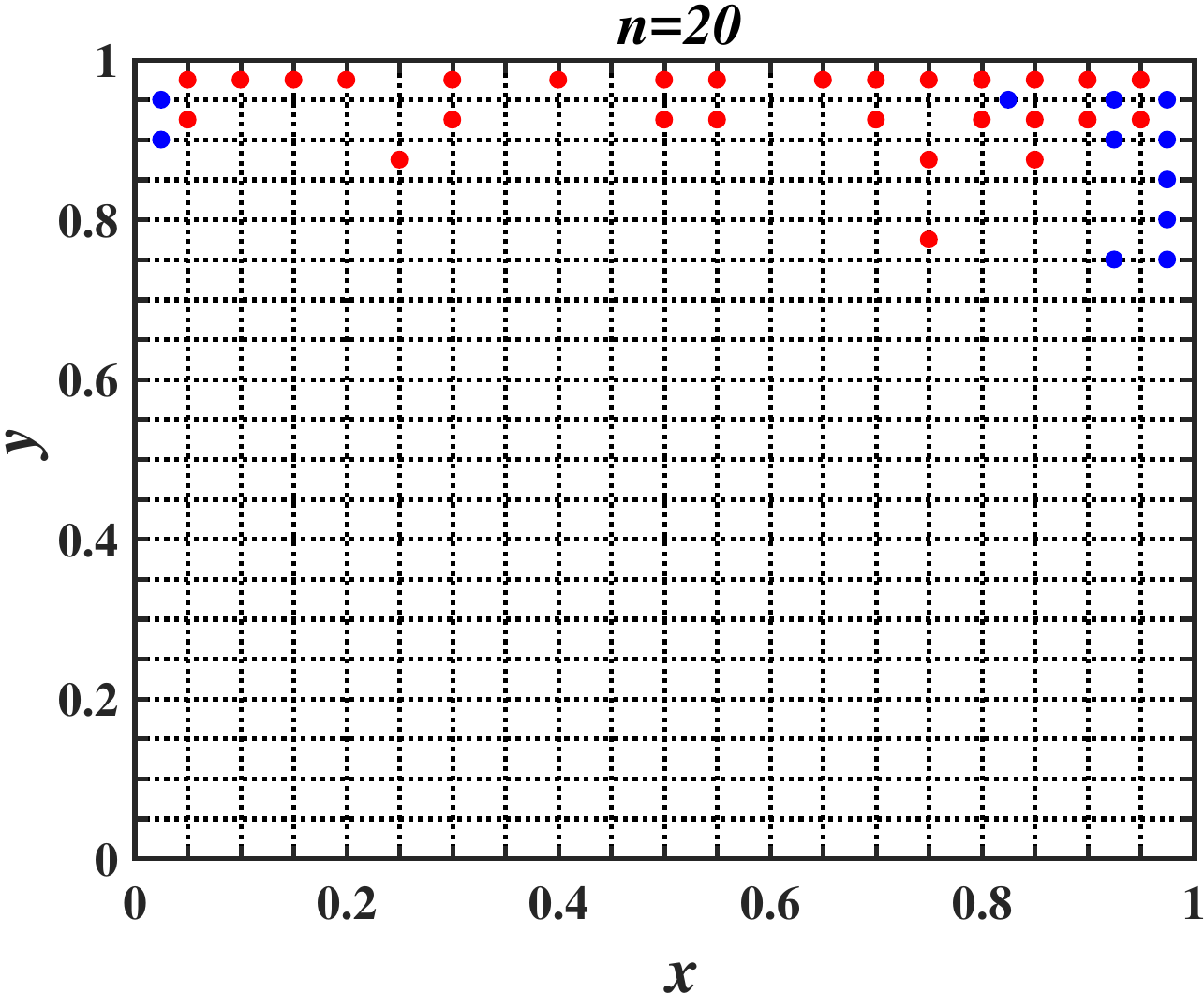}
\caption{Stokes (Column 1) and Navier-Stokes (Column 2) streamlines of ROC approximations at $\text{Re} = 1000, \nu = 0$ with $n_x =n_y =200$ (NS: $n_x =n_y =120$) and various RB dimension $n$. Column 3 and 4: selected collocation points for Stokes and Navier-Stokes cases with training set $[10,50] \times [0,1]$, and $n_x =n_y =20$. The corresponding $2n-1$ collocation points are showed at each subplot with different $n$ and the whole MAC points are also showed in the first subplots in Column 3 and 4, where red dots correspond to $u$-components and blue ones correspond to $v$-components.}
\label{fig:navierstokes:streamline3} 
\end{figure}

We now present the streamline for both equations. For Stokes, the streamline of the {ROC approximations of the} FDM solution for $(\text{Re},\nu) = (1000,4)$ with $n= 2, 10, 15$ {and their corresponding error $E_q(n)$} are in Figure~\ref{fig:stokes:streamline}. {Those quantities for the} Navier-Stokes with {two sets of parameter values} and $n = 2, 10, 40$ are {in Figure ~\ref{fig:Nstokes:streamline}}. We see that for both equations, even with a bare minimum of two basis, the ROC streamline is {quite} accurate. As $n$ increases, the error quickly decreases until reaching machine accuracy with 15 basis for both equations when Reynolds number is $1000${, part of} the training grid. {We test the more challenging} Navier-Stokes {equation on the} parameter value $(500,2)$ that is outside of the training grid. {We see that,} as $n$ increases from $2$ to $40$, the largest error decreases from $10^{-3}$ to $10^{-5}$. 
To investigate the distribution of selected collocation points, streamlines {corresponding to} different Reynolds numbers together with the selected collocation points are shown in Figure~\ref{fig:navierstokes:streamline3}. As $n$ increases, more collocation points are identified. It is interesting to note that most points concentrate on the upper-right corner and some points are scattered around the upper-left corner. More {interestingly}, 
{points are only selected on the edges of MAC cells, that is, the collocation points are all from the velocities. This is because the reduced basis are divergence free, thus, the residual of continuity equation is negligible.}
To {showcase} this phenomena more intuitively, we show two examples in the third and fourth columns in Figure~\ref{fig:navierstokes:streamline3}.

\subsection{Time-dependent Navier-Stokes equations}

Finally, we consider the nonlinear time-dependent Navier-Stokes equations described in Section \ref{sec:lid}. {We vary the two parameters one at a time.} That is, $\bmu = (\text{Re},\nu) \in \Omega_p =[10, 500] \times \{0\}$ or $\{500 \} \times [0,2]$. The high fidelity solutions {are obtained} with $n_x = n_y =100, \tau =0.002, T=35$ {adopting} backward Euler {as the time evolution} scheme. The training set of the ROC algorithm is {trained by uniformly sampling} the parameter domain $21$ {times}. {The testing set consists of $5$ uniformly sampled points on} $[10+2.2, 500-2.2] \times \{0\}$, or $\{500\} \times [0+0.22,2-0.22]$, which does not intersect with the training set. For the control parameter, we set $p_{adap}=40\%, n_{adap}^{max}=300$ (one percent of the degrees of freedom in the underlying FDM scheme). Under this problem setting, the offline process {generates} $70$ ($90$) basis and $344$ ($835$) collocation points .

\begin{figure}[thbp]
\centering
\includegraphics[width=0.31\textwidth]{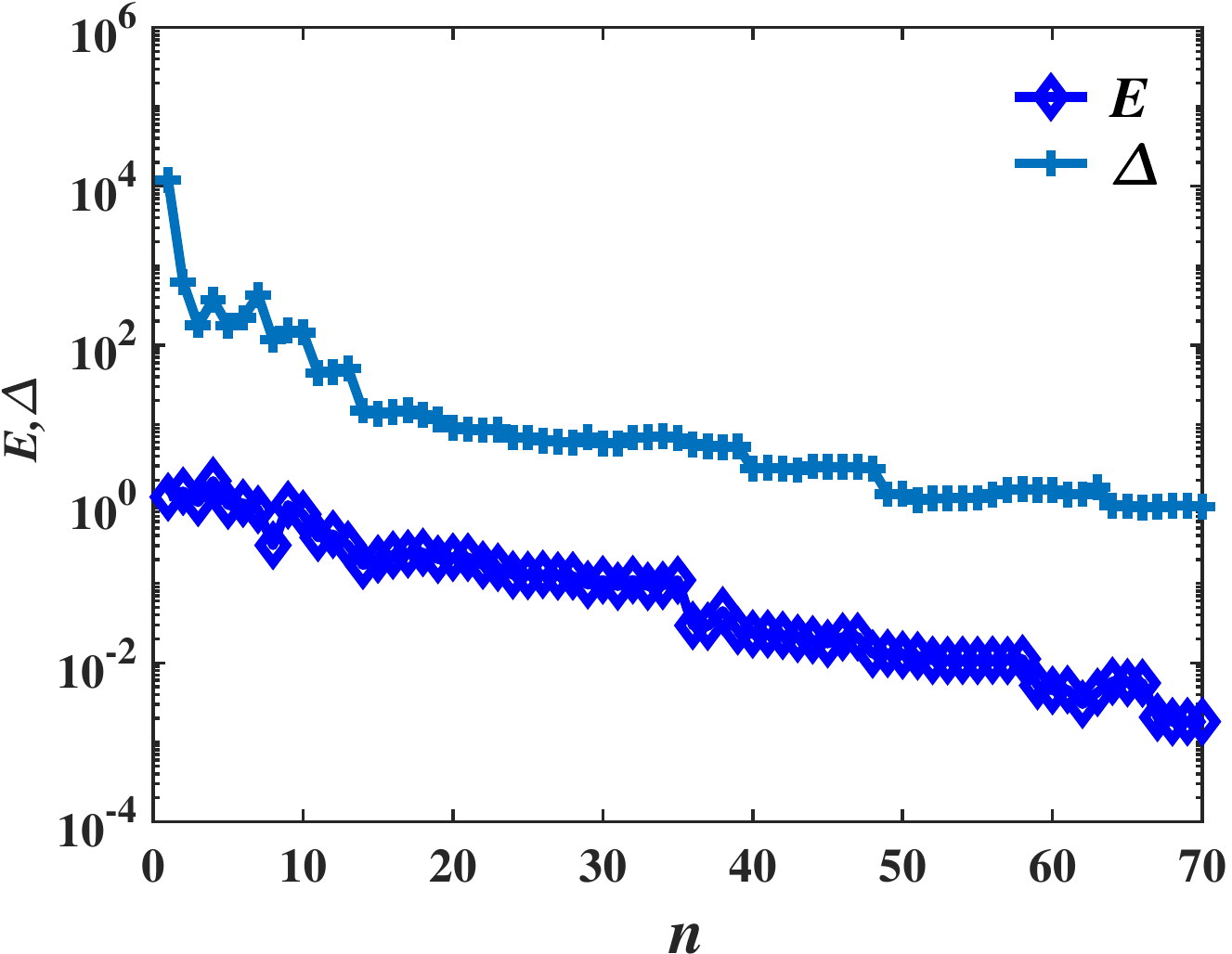}
\includegraphics[width=0.31\textwidth]{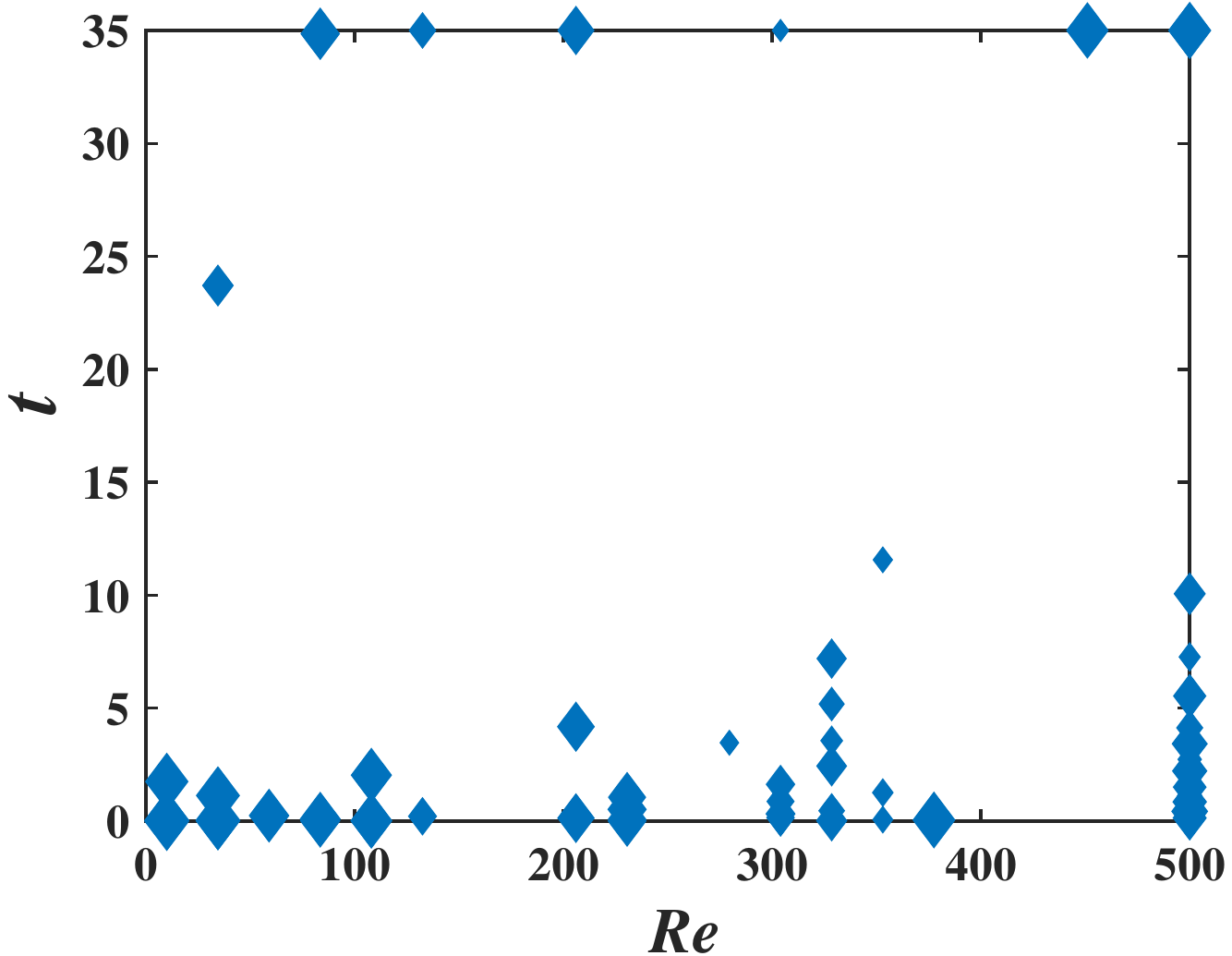}
\includegraphics[width=0.34\textwidth]{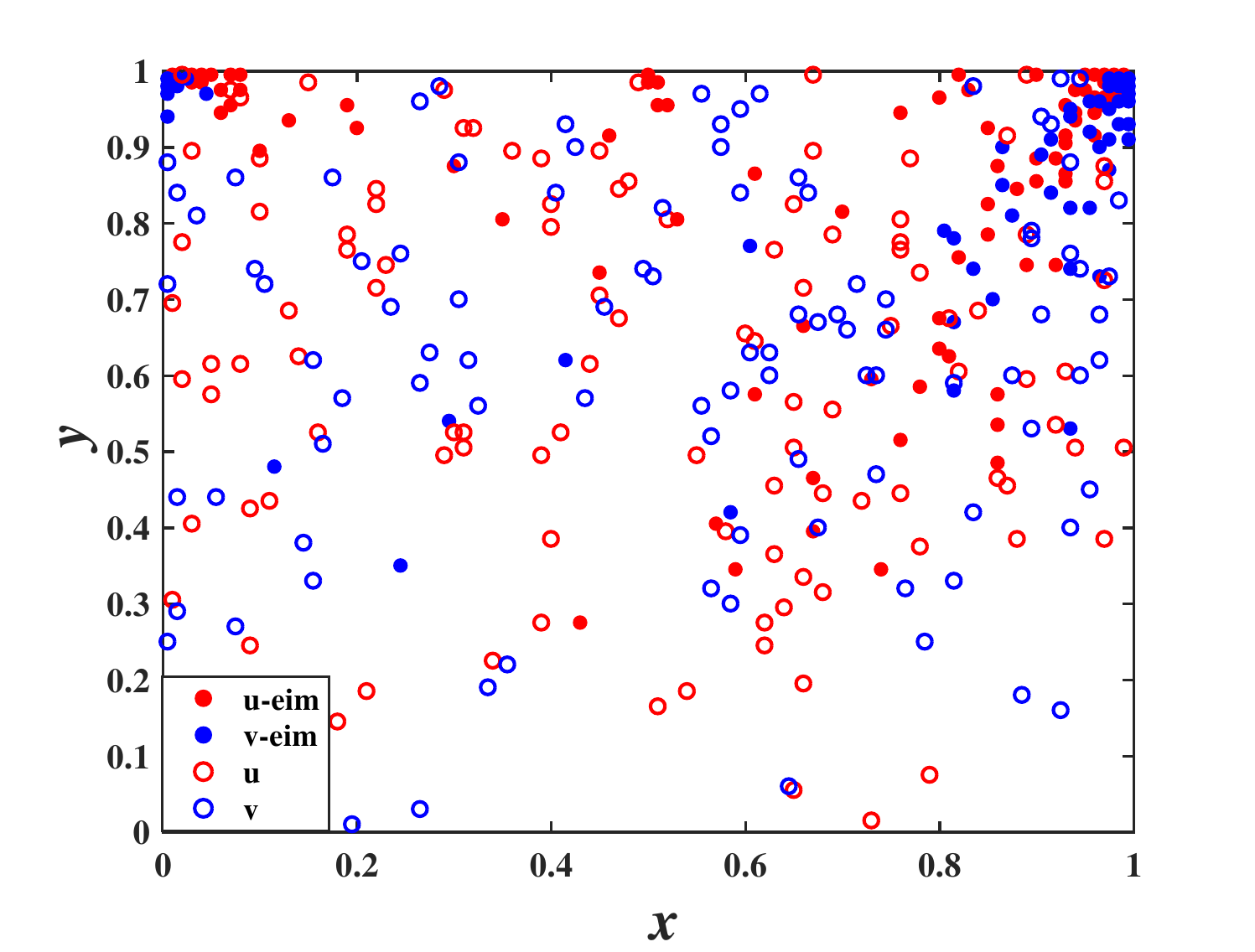}\\
\includegraphics[width=0.31\textwidth]{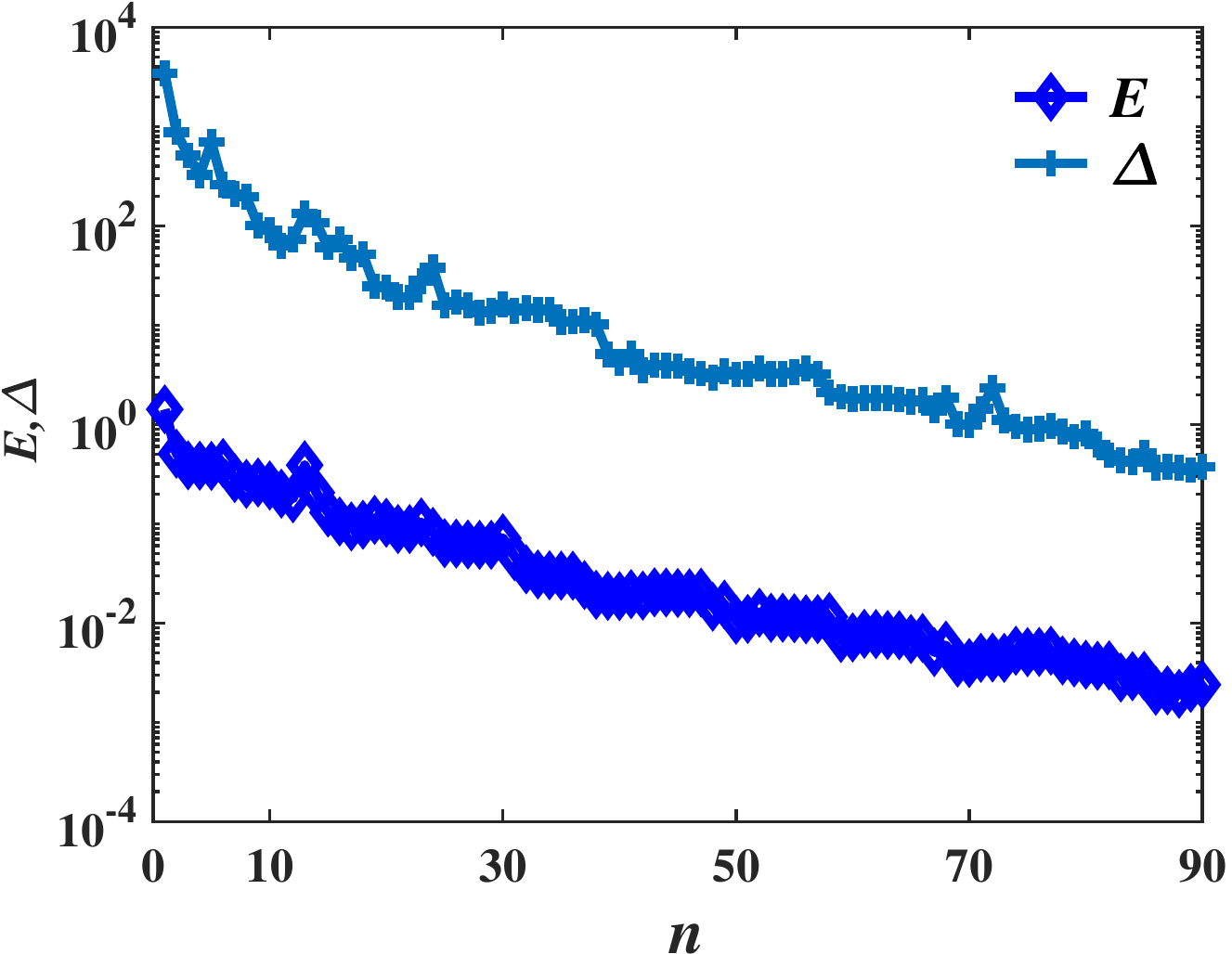}
\includegraphics[width=0.31\textwidth]{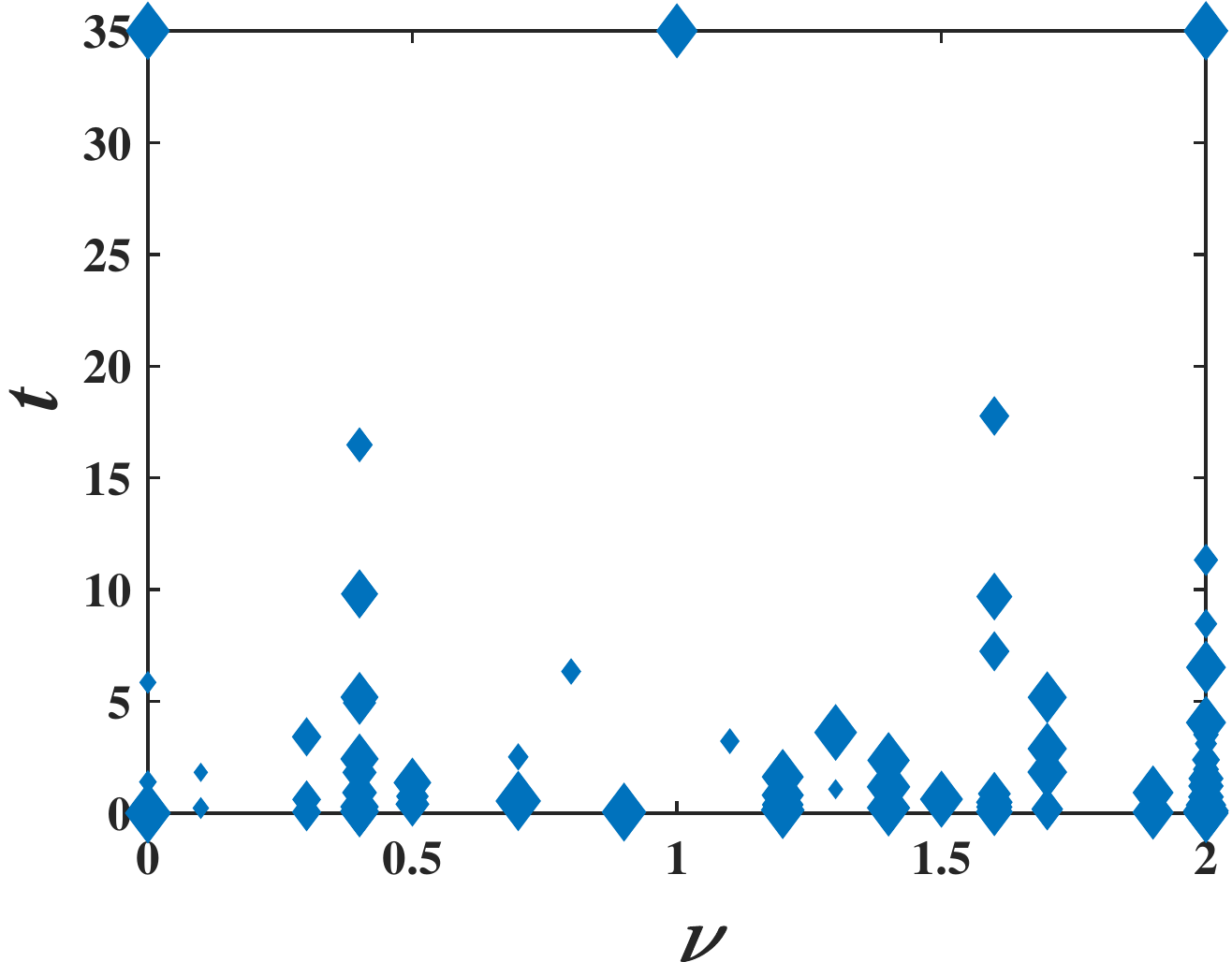}
\includegraphics[width=0.34\textwidth]{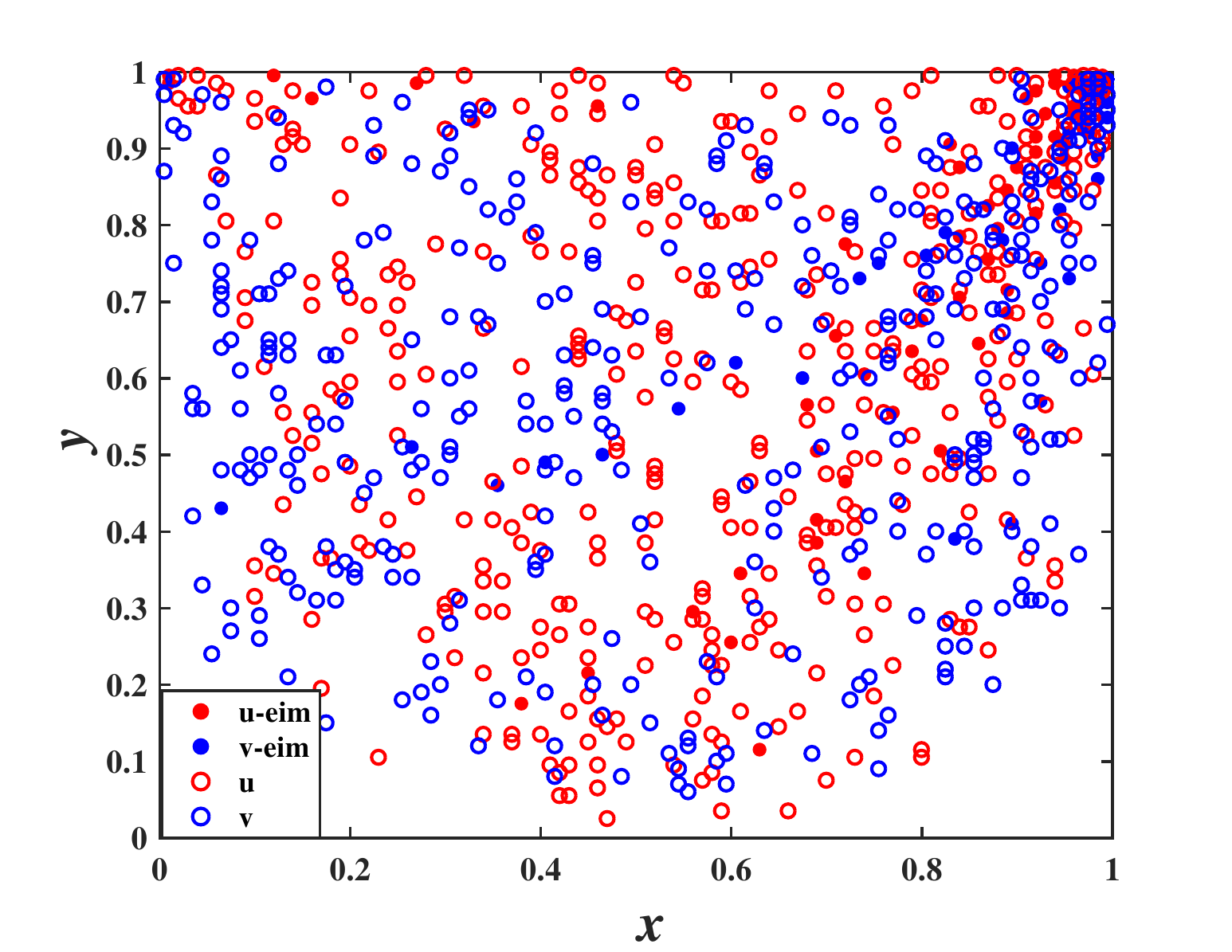}
\caption{History of convergence for the relative error and error estimator (left), the selected parameter-time pairs (middle), and selected collocation points (right) by adaptive ROC algorithm. Top is for the training set $\Omega_p =[10, 500] \times \{0\}$, and bottom $\{500 \} \times [0,2]$.}
\label{fig:lidDriven:error} 
\end{figure}

The error and error estimator curves are shown in Figure \ref{fig:lidDriven:error} left. 
Compared with the Stokes and steady-state Navier-Stokes equations, {these curves decrease slower} due to the nonlinearity and time evolution. {We note that this is consistent with the literature. For example,} 
authors of \cite{stabile2018finite} proposed a stabilized Galerkin projection based reduced basis method of unsteady incompressible NS equations with moderate Reynolds number. 
The relative error in $L^2$ norm at different time nodes for the velocity is at the $10^{-2}\sim 10^{-3}$ level, {while that of pressure is around} $10^{-2}$. 
{The chosen parameter-time pairs} are displayed in Figure \ref{fig:lidDriven:error} center.
Those with larger symbols are selected earlier.
{We show the} selected collocation points in Figure \ref{fig:lidDriven:error} {right which continue to not include any}  center of the MAC cells. Though much more collocation points are selected compared to the steady state NS and Stokes equations, collocation points in the upper left and upper right corners are still {more likely to be chosen than those from other parts of the domain. To show the effect of the adaptive enrichment,} we denote the collocation points  from the adaptive {enrichment with circles while the rest with filled diamonds}. The distribution of {the latter} is similar with that in the steady state NS equations.

\begin{figure}[thbp]
\centering
\includegraphics[width=0.32\textwidth]{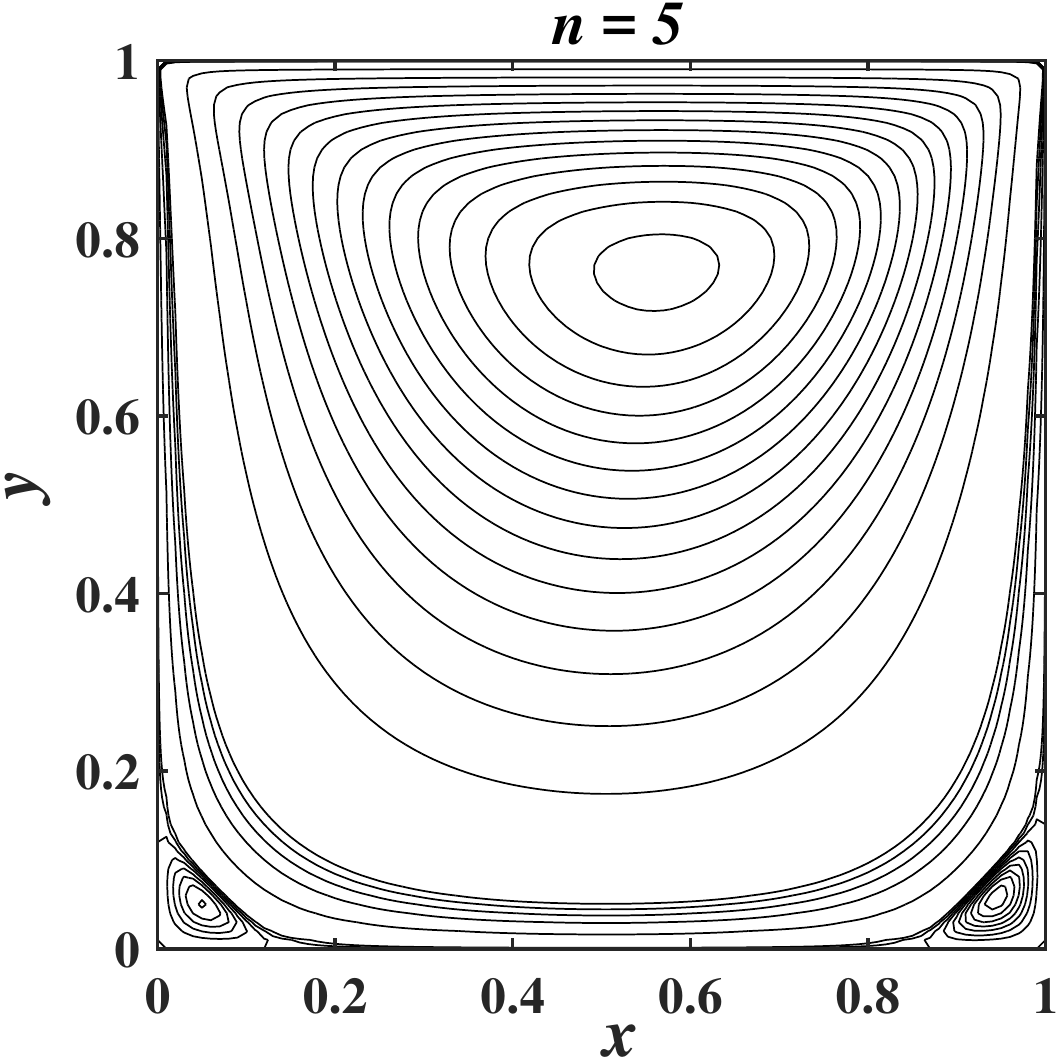}
\includegraphics[width=0.32\textwidth]{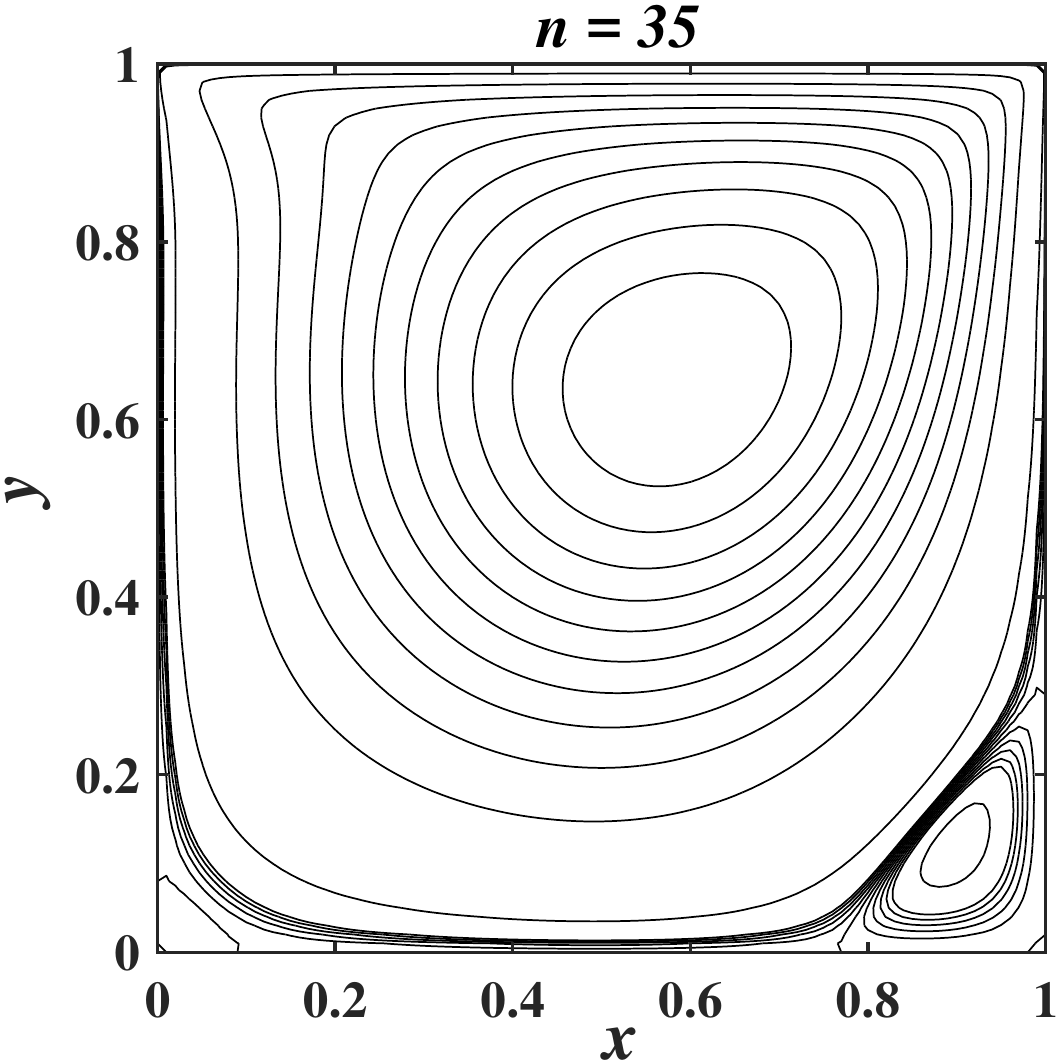}
\includegraphics[width=0.32\textwidth]{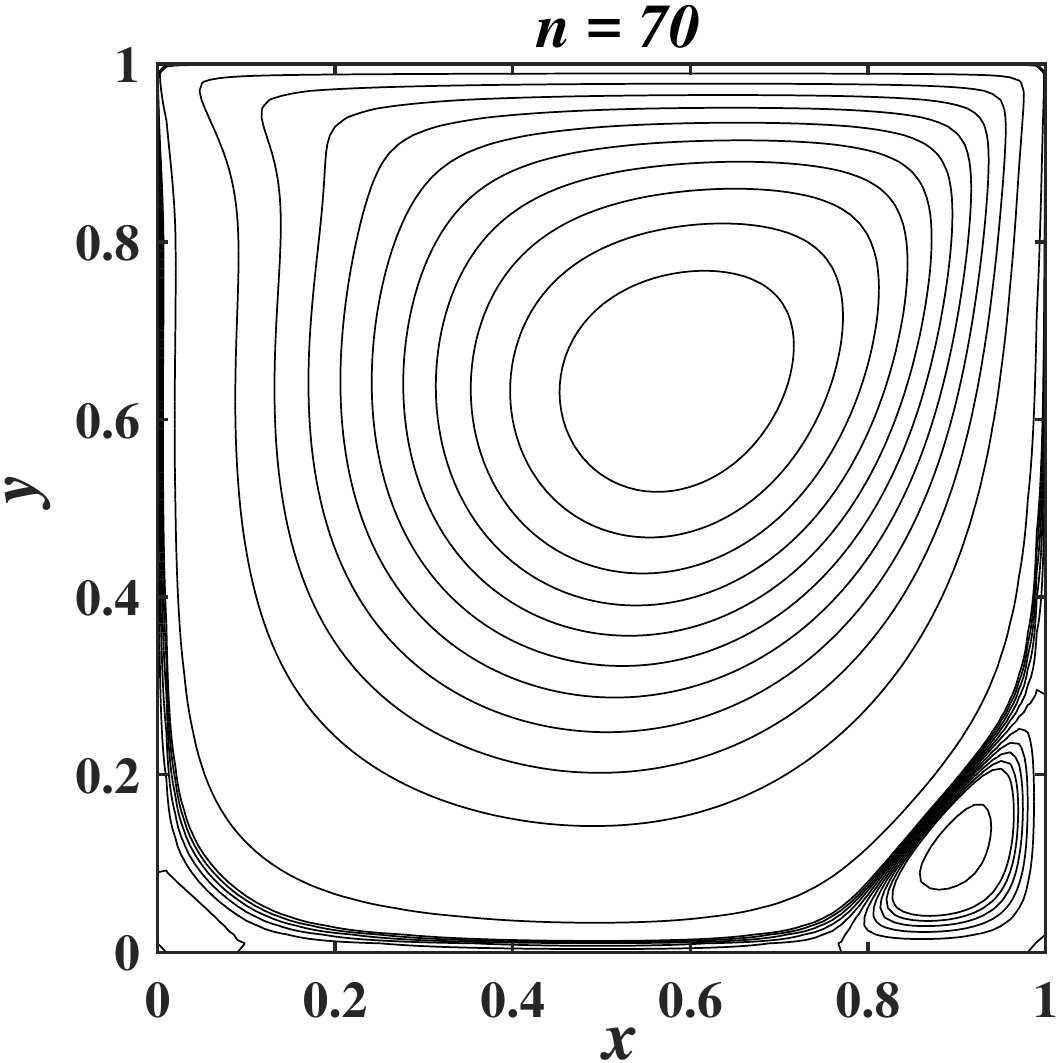}\\
\includegraphics[width=0.32\textwidth]{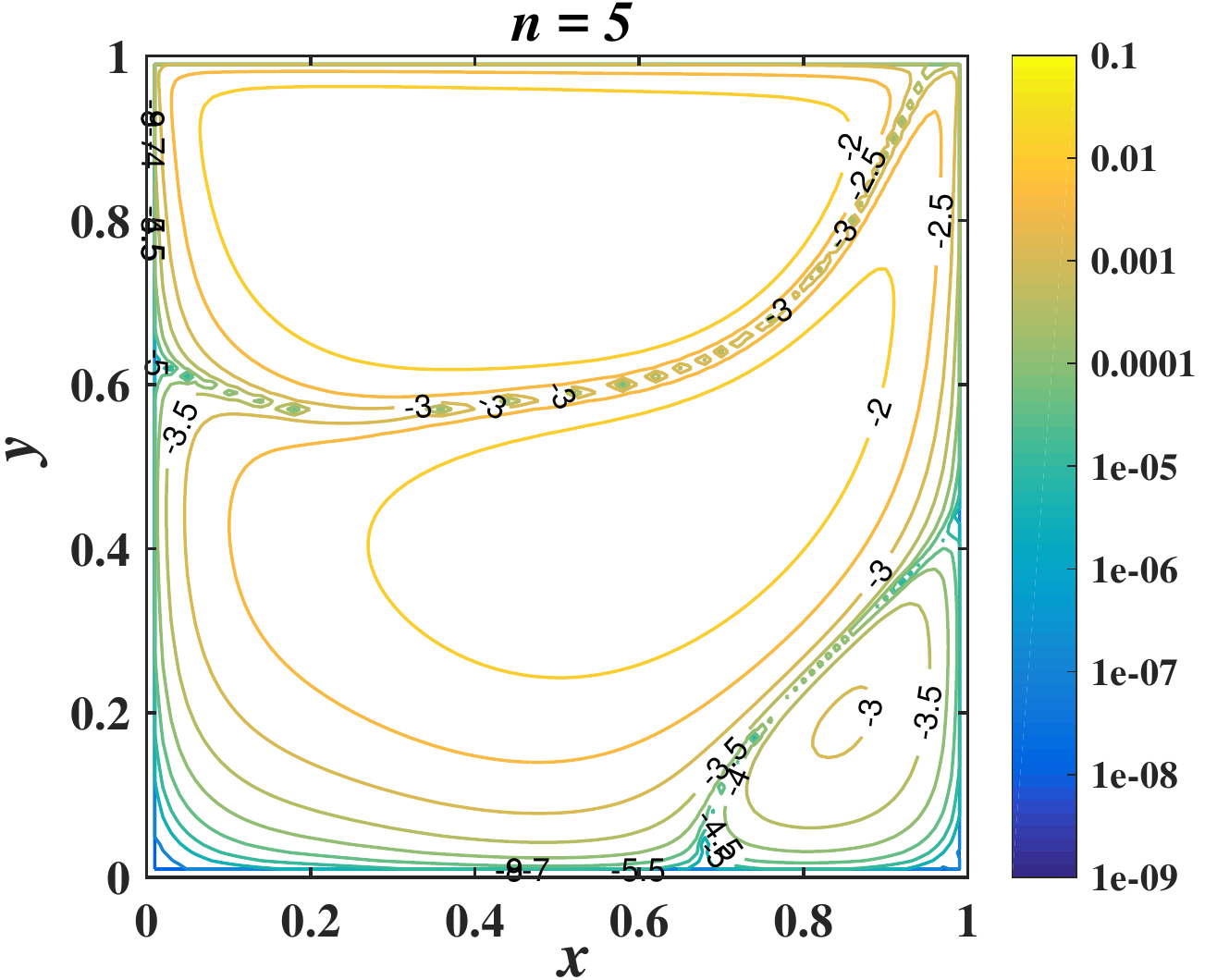}
\includegraphics[width=0.32\textwidth]{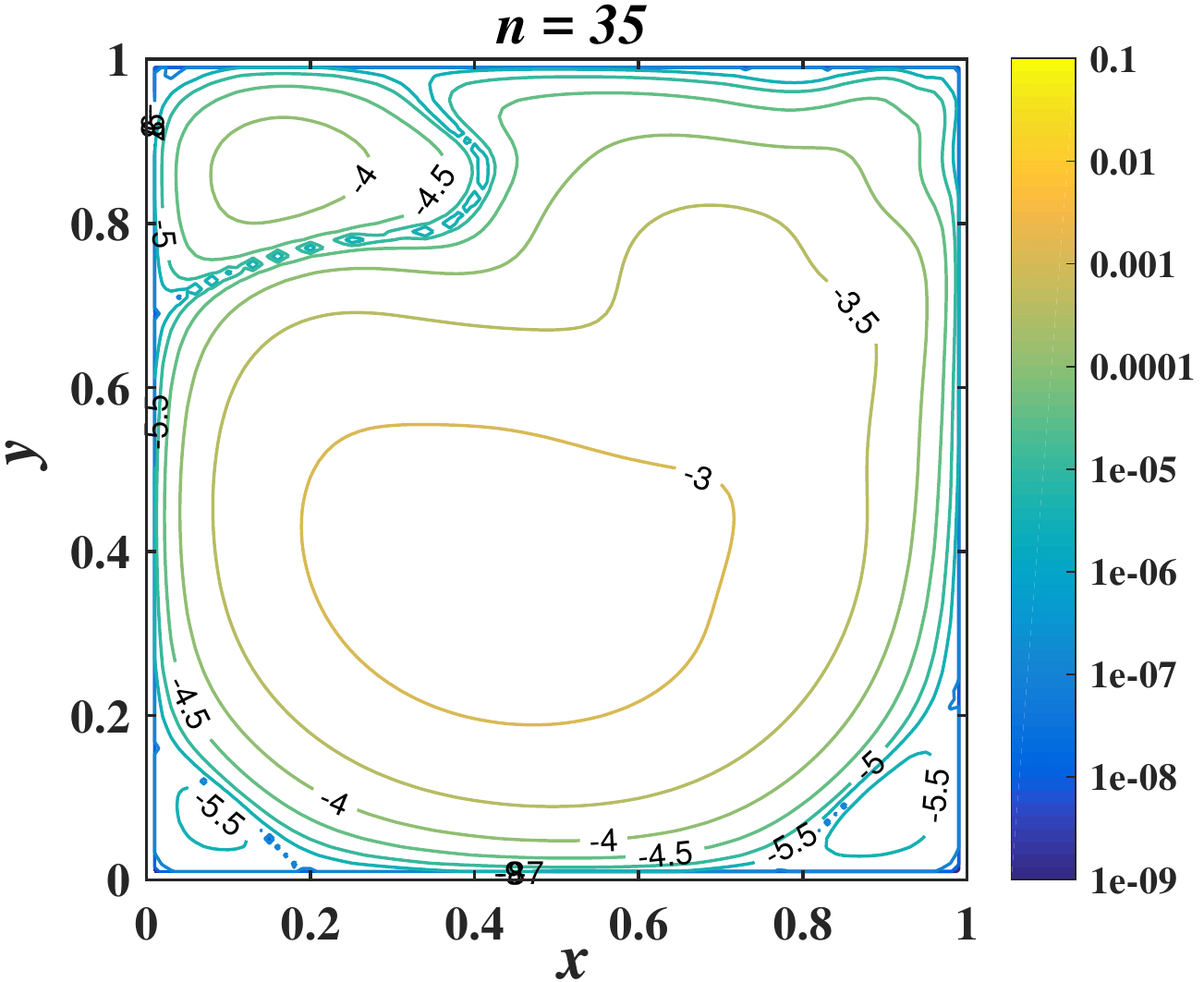}
\includegraphics[width=0.32\textwidth]{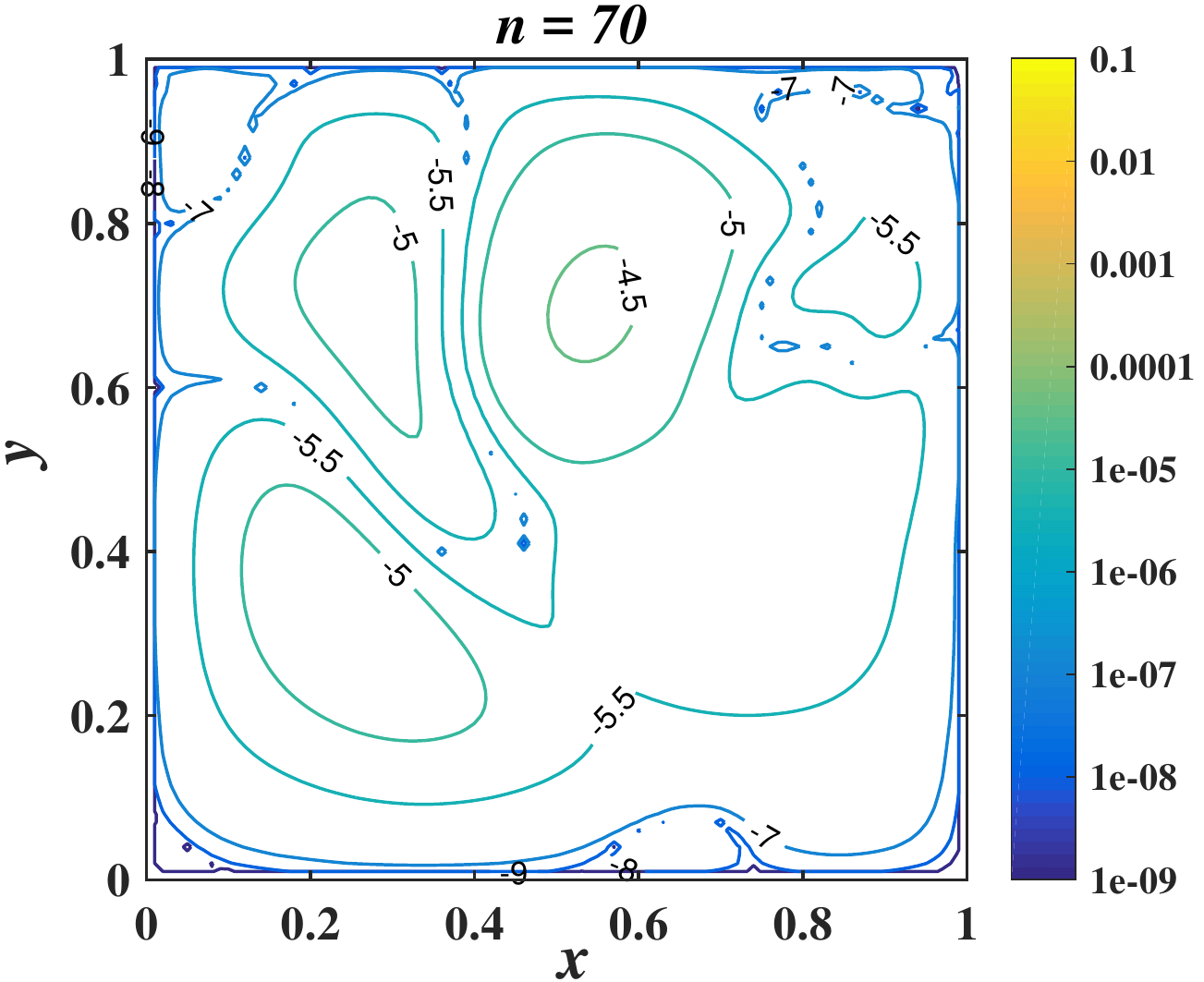}\\
\includegraphics[width=0.32\textwidth]{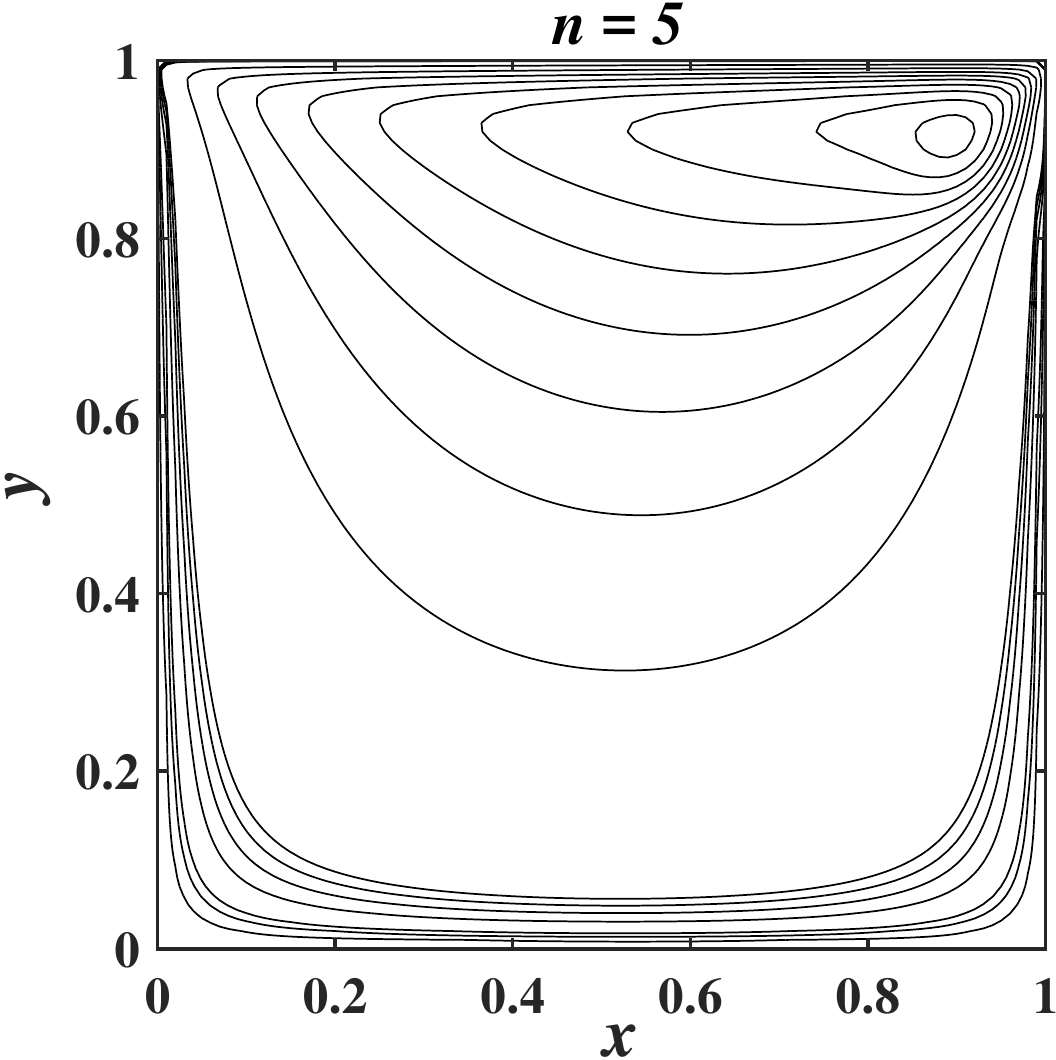}
\includegraphics[width=0.32\textwidth]{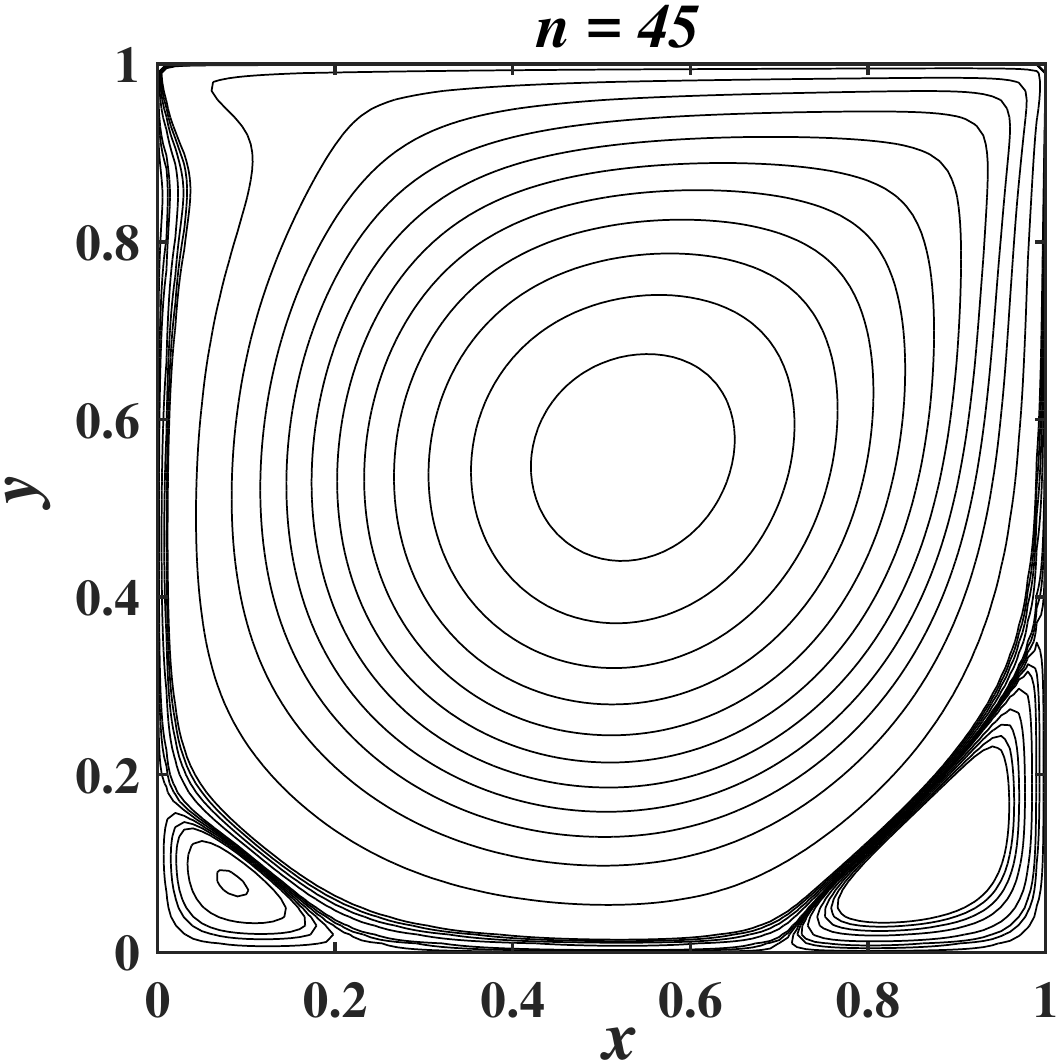}
\includegraphics[width=0.32\textwidth]{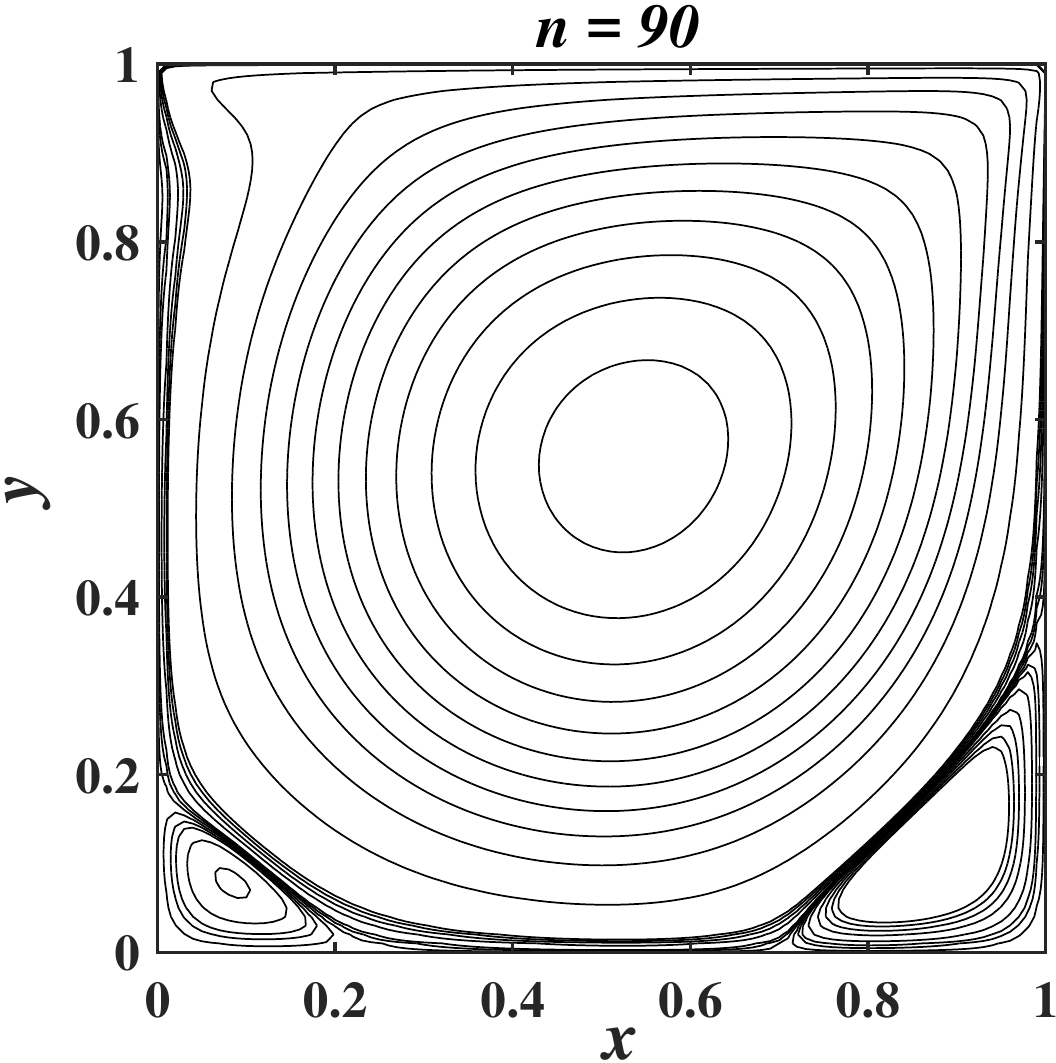}\\
\includegraphics[width=0.32\textwidth]{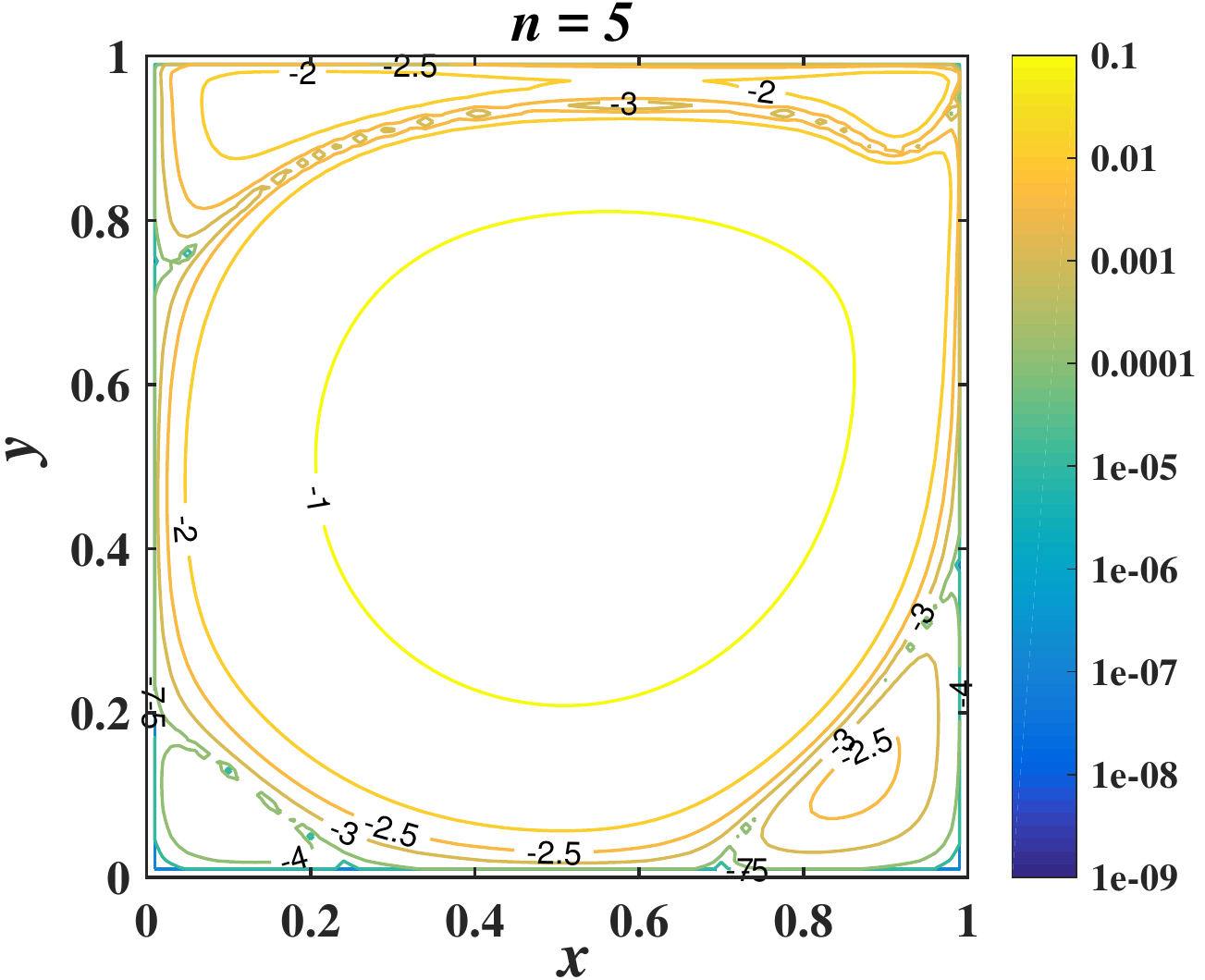}
\includegraphics[width=0.32\textwidth]{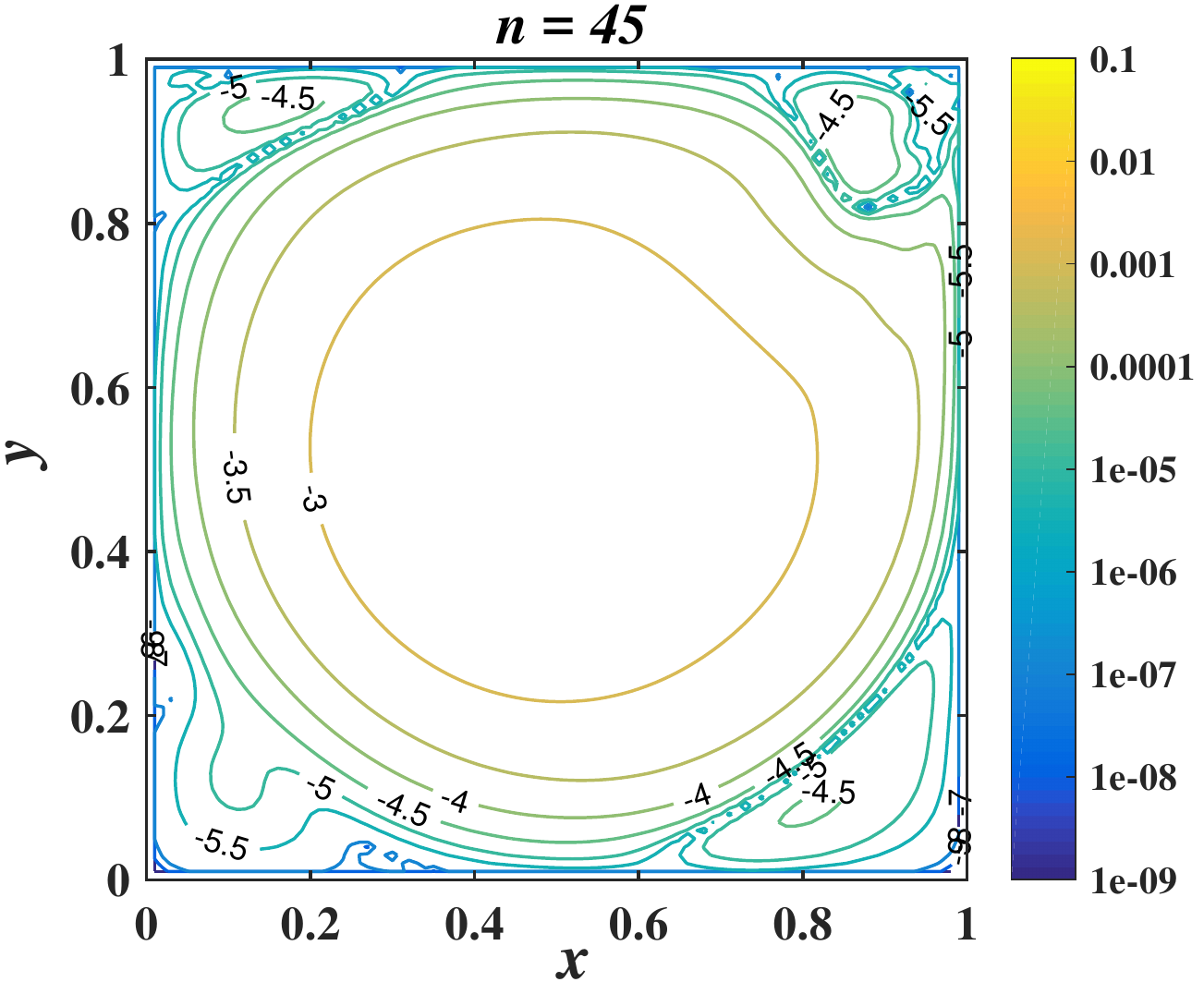}
\includegraphics[width=0.32\textwidth]{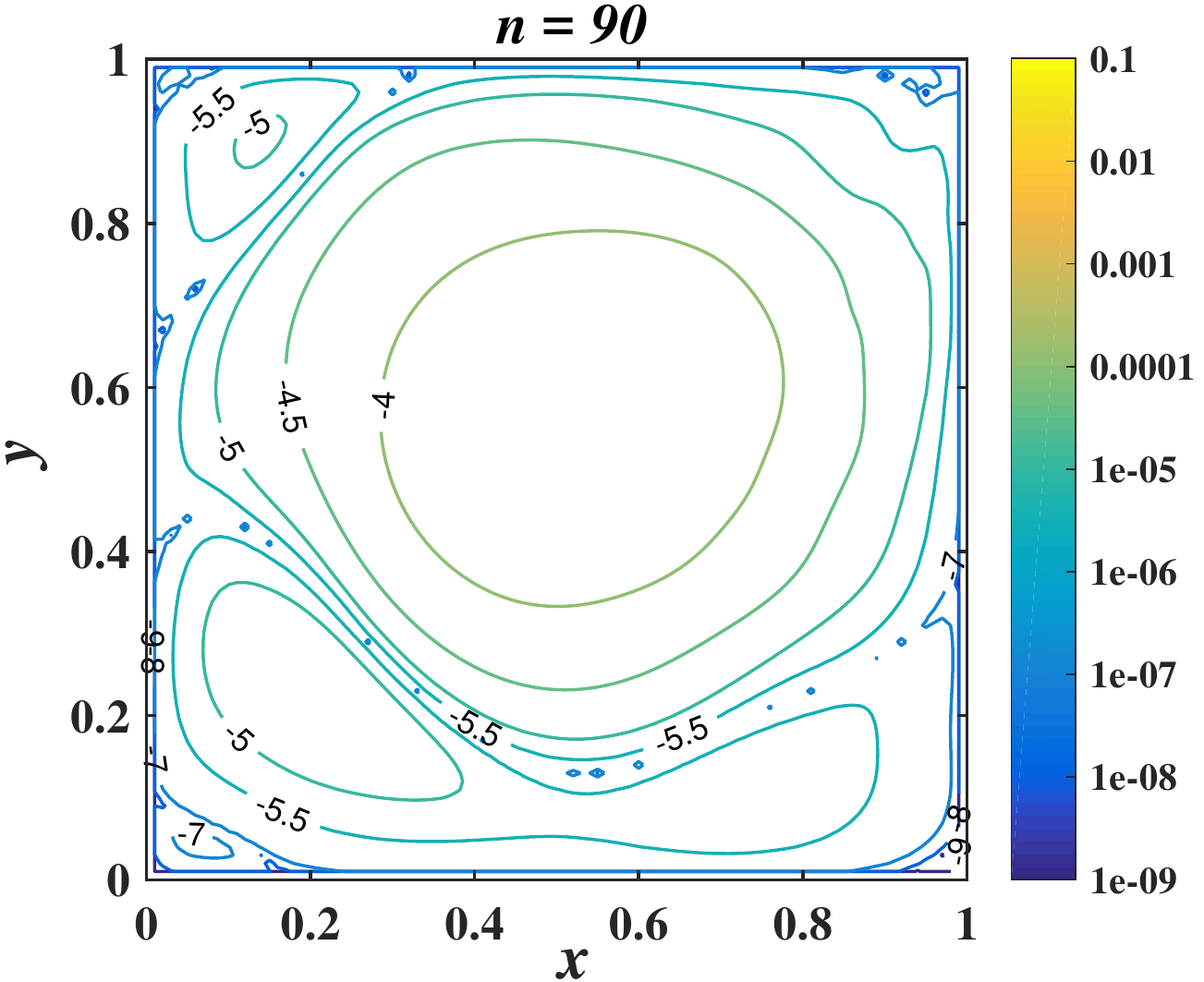}\\
\caption{Streamlines of ROC approximations with the denoted RB dimensions $n$, and the absolute streamline errors $E_q(n)$ in logarithmic scales for time dependent Navier-Stokes equation at $T=35$ (Rows 1,2 $(\text{Re}=255, \nu=0)$, rows 3,4 $(\text{Re}=500, \nu=1.78)$).}
\label{fig:lidDriven:streamline} 
\end{figure}

\begin{figure}[thbp]
\centering
\includegraphics[width=0.32\textwidth]{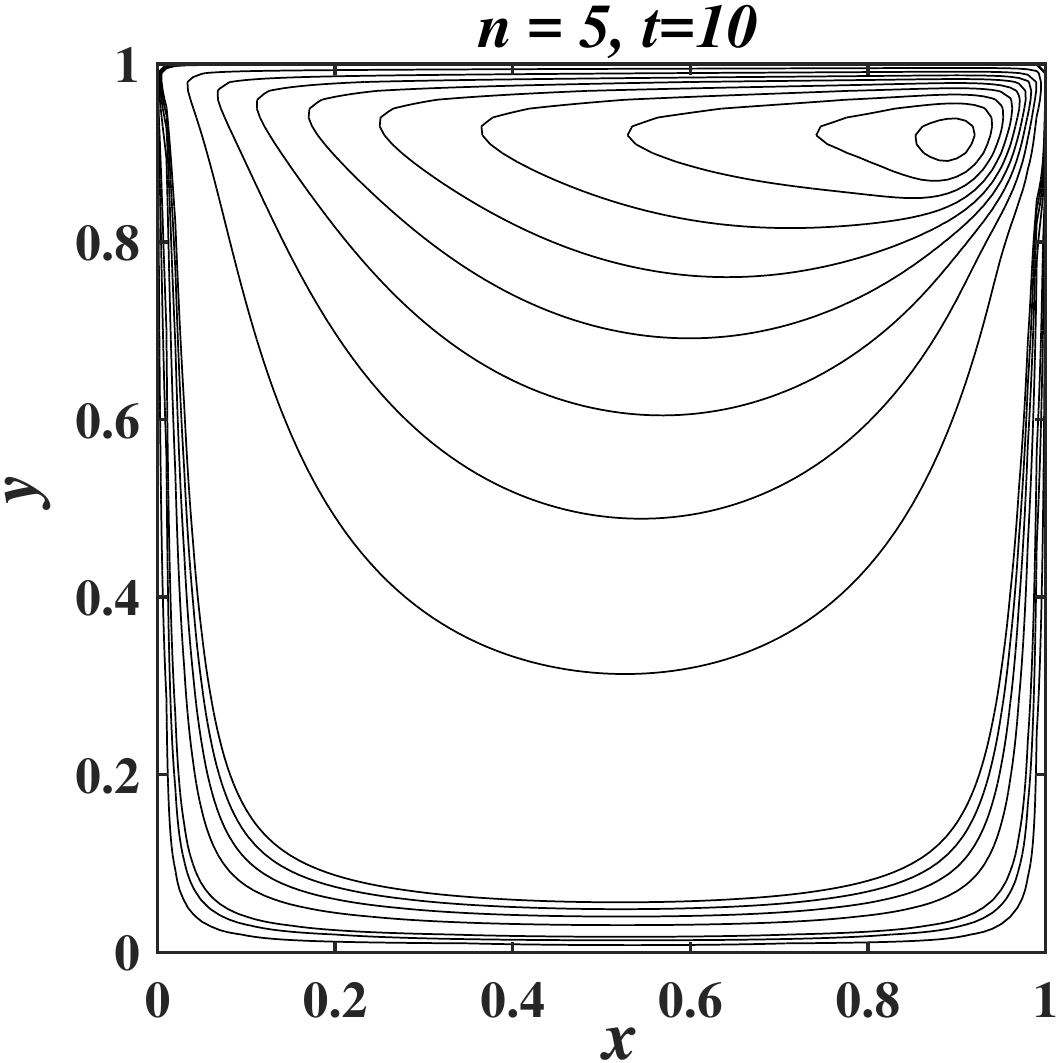}
\includegraphics[width=0.32\textwidth]{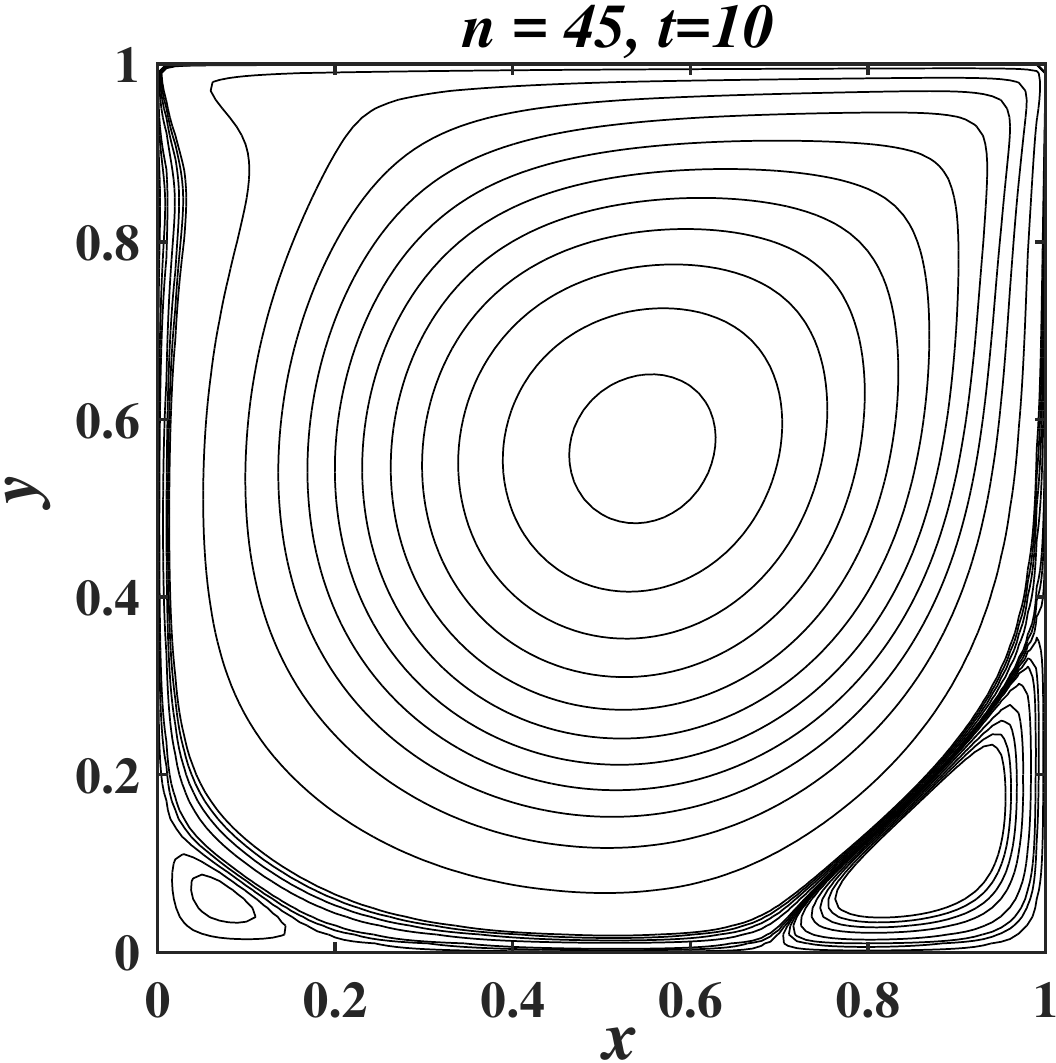}
\includegraphics[width=0.32\textwidth]{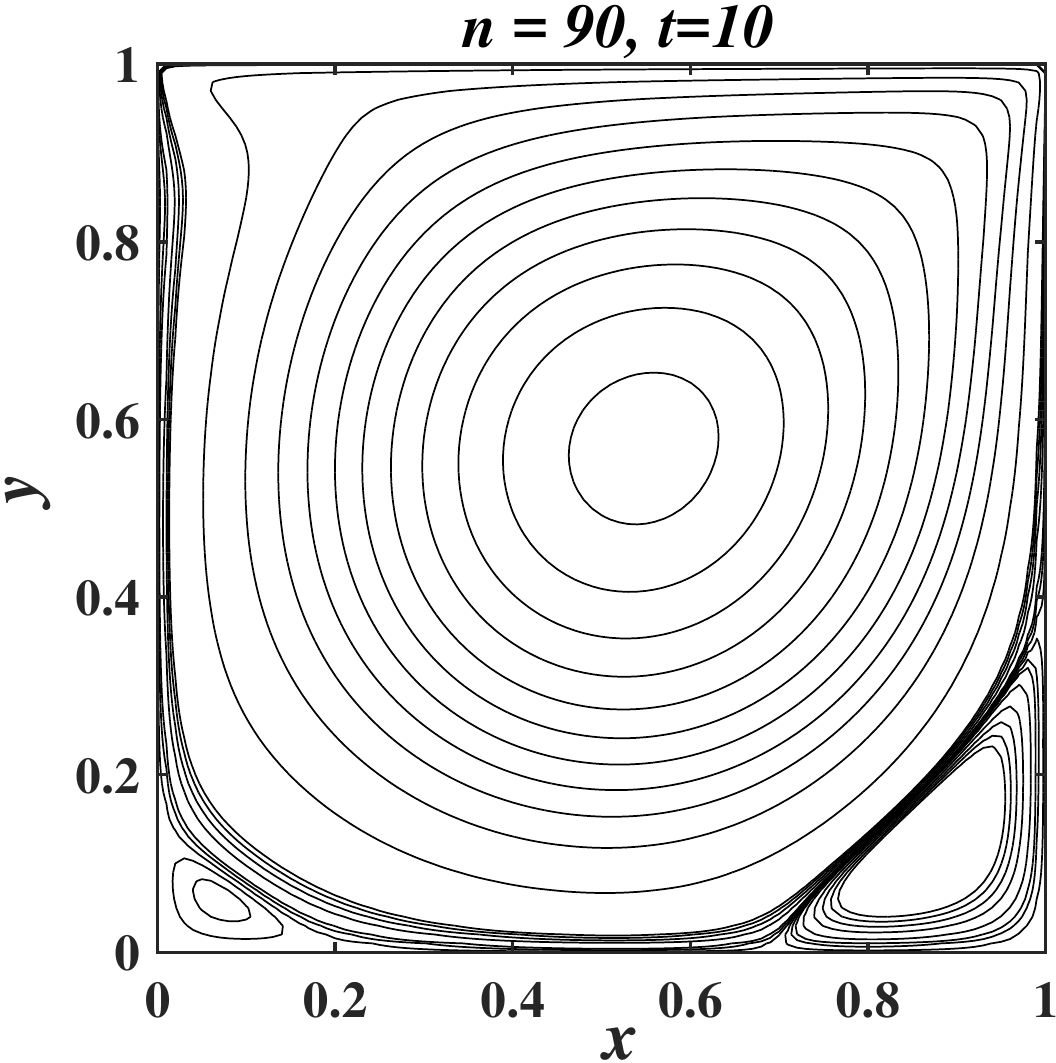}\\
\includegraphics[width=0.32\textwidth]{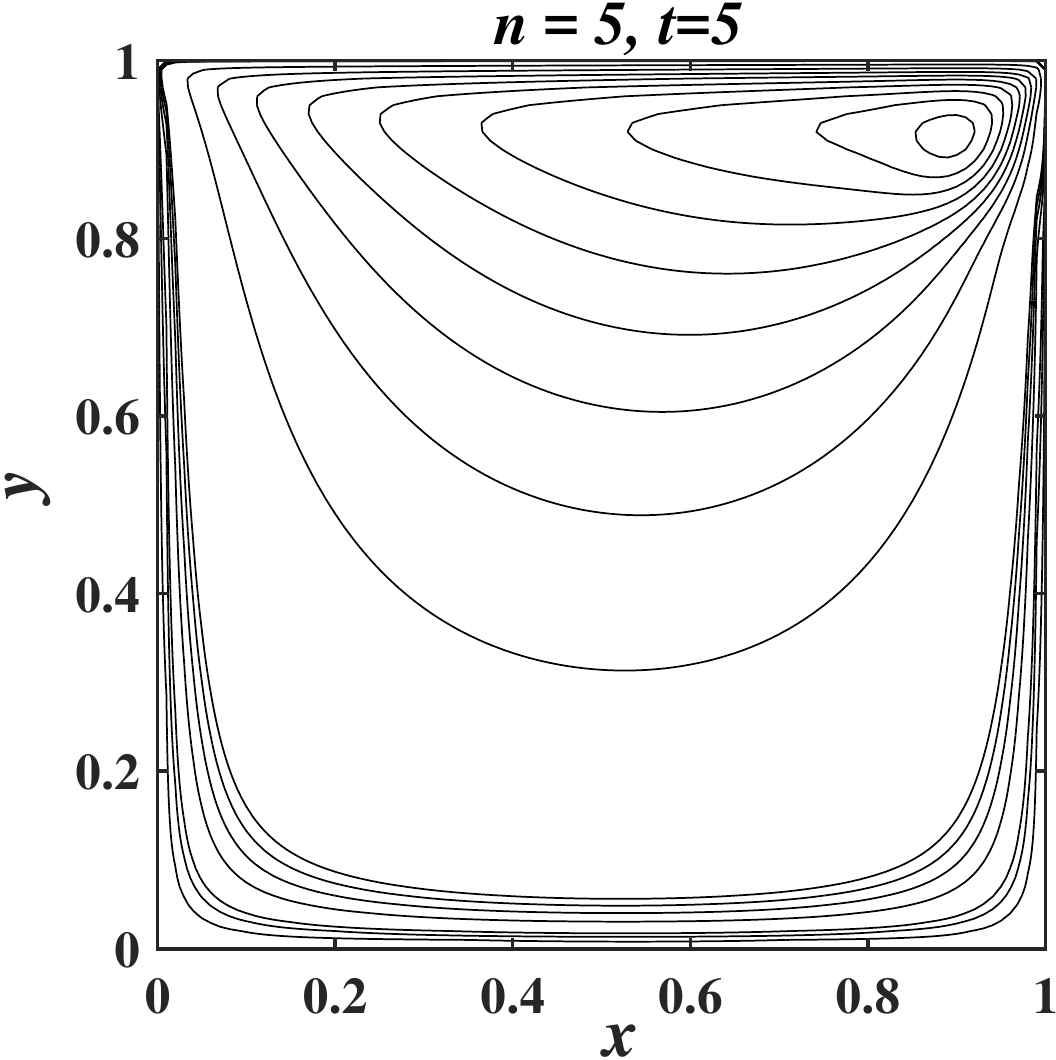}
\includegraphics[width=0.32\textwidth]{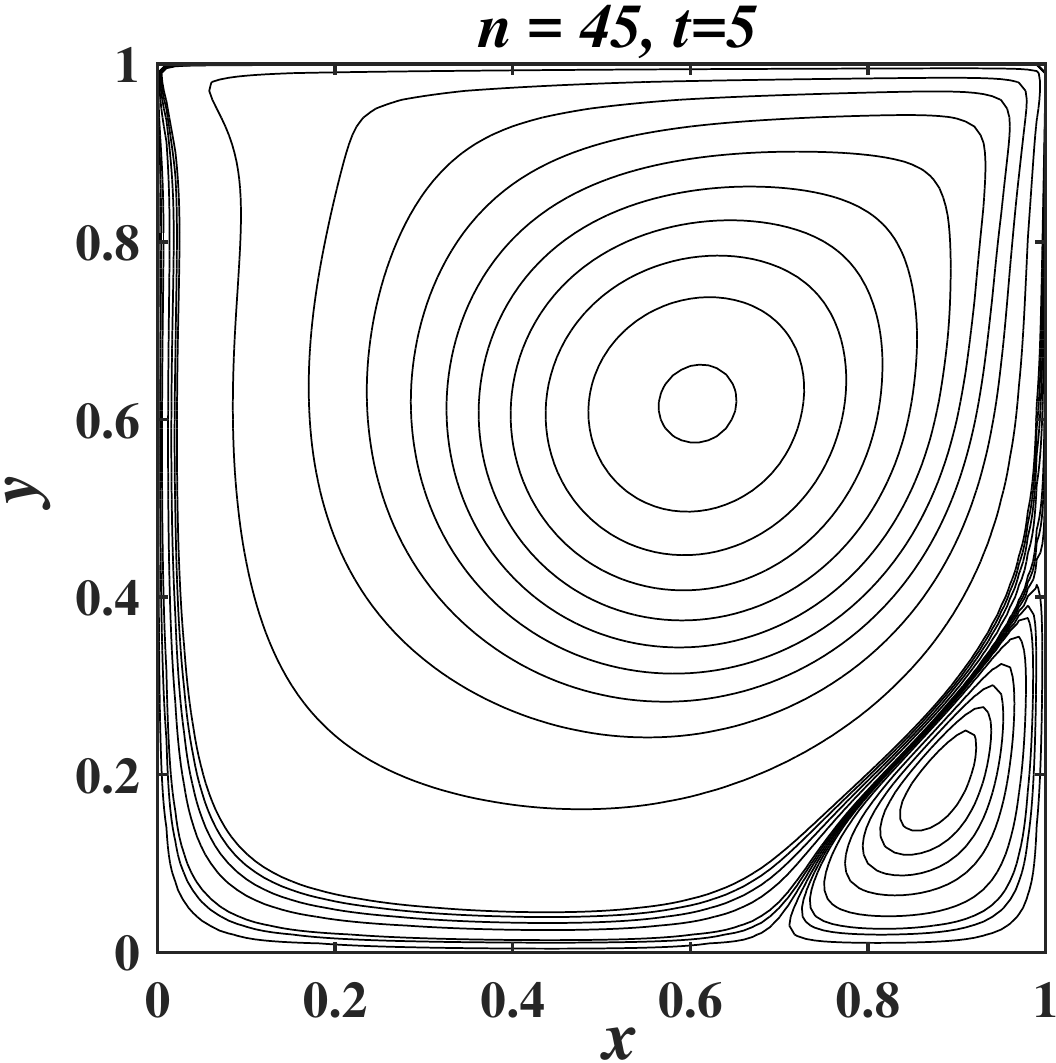}
\includegraphics[width=0.32\textwidth]{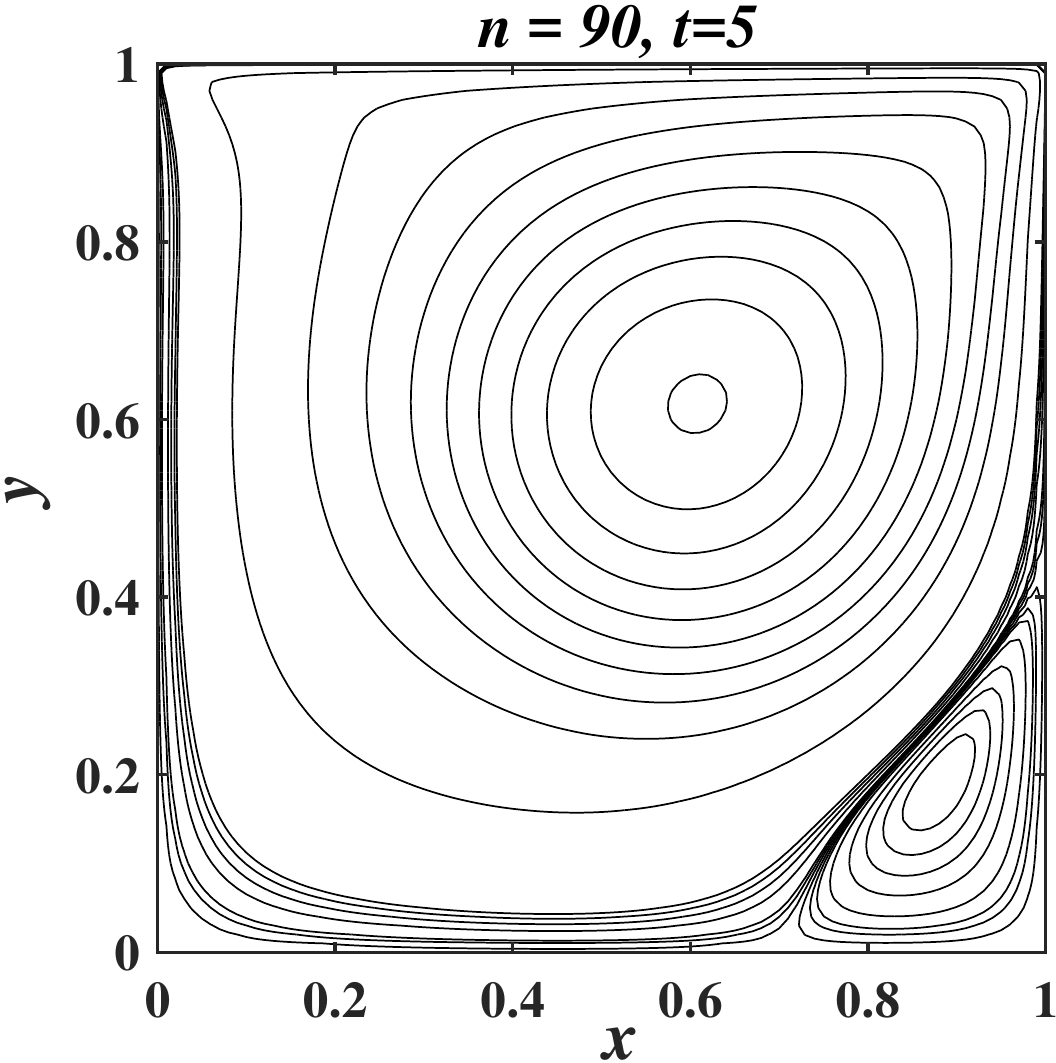}\\
\includegraphics[width=0.32\textwidth]{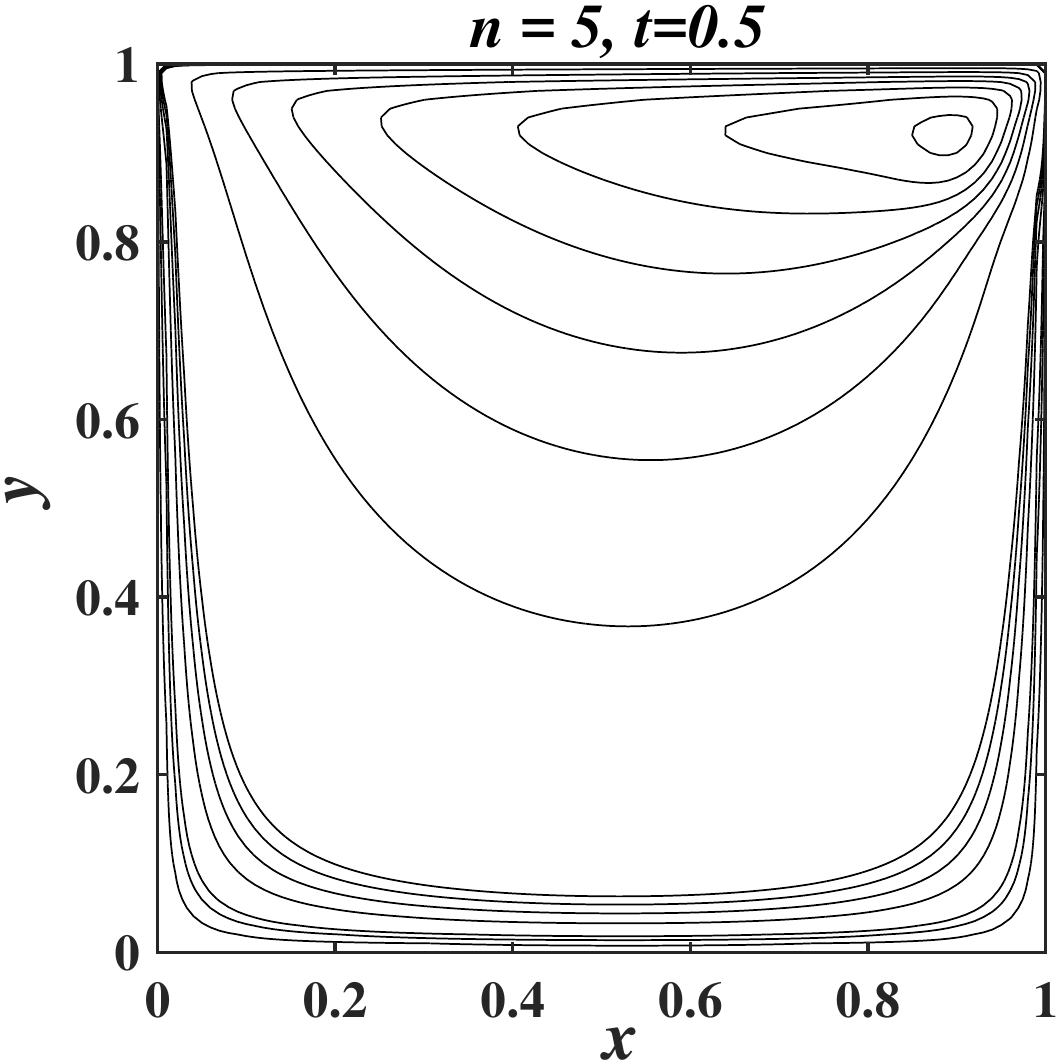}
\includegraphics[width=0.32\textwidth]{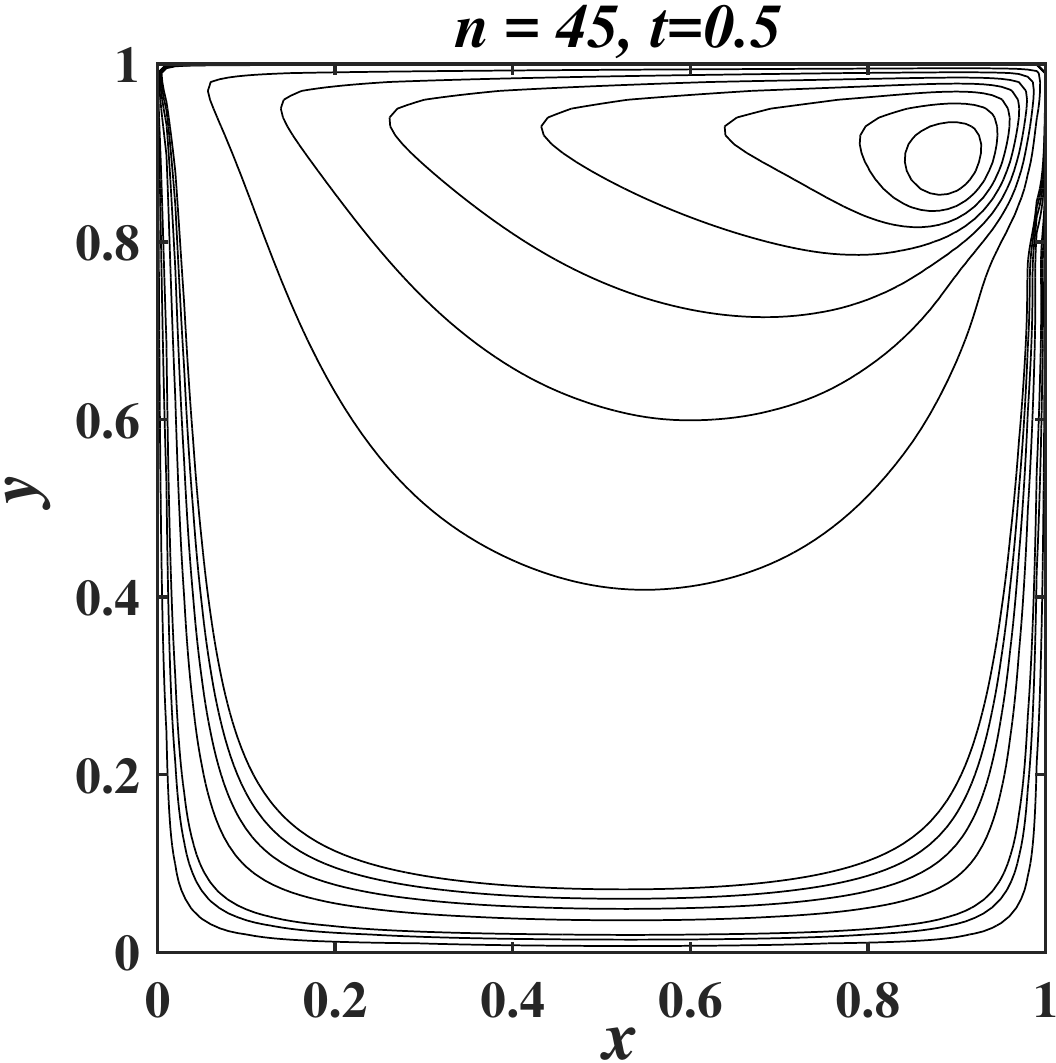}
\includegraphics[width=0.32\textwidth]{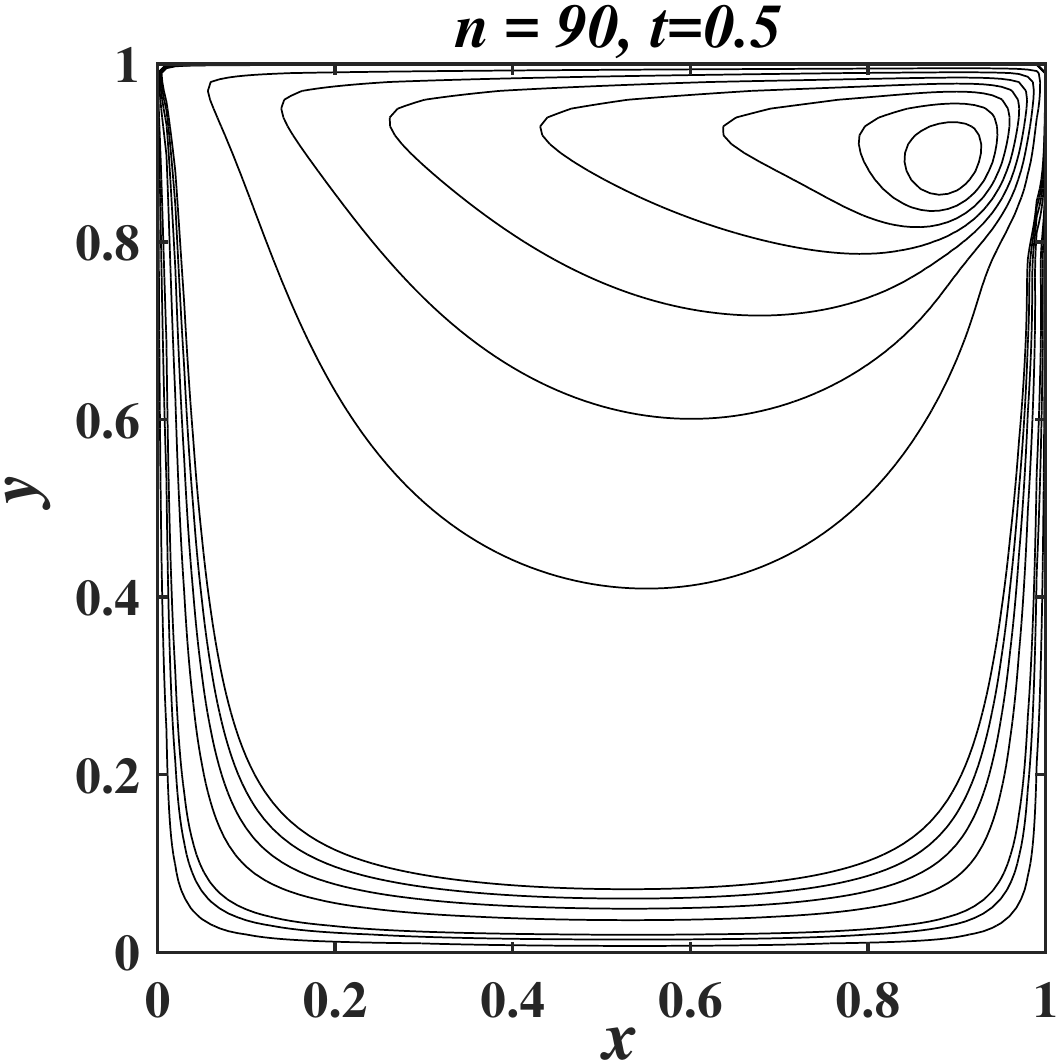}
\caption{Streamlines of ROC approximations with the denoted RB dimensions $n$ for time dependent Navier-Stokes equation at $t=0.5, 5, 10$ with $(\text{Re}=500, \nu=1.78)$.}
\label{fig:lidDriven:streamline2} 
\end{figure}

Streamlines and its ROC approximation errors at $(\text{Re}=255, \nu=0)$ ($(\text{Re}=500, \nu=1.78)$) with different basis numbers $n=5, 35, 70$ ($n=5, 45, 90$) at $T=35$ are shown in Figure~\ref{fig:lidDriven:streamline}. As $n$ increases from $5$ to $70$ ($90$), the largest error decreases from $10^{-2}$ to $10^{-4}$ and this is consistent with the decreasing rate of the error curve in Figure \ref{fig:lidDriven:error} left. Also, we investigate the streamlines at different time nodes $t=0.5, 5, 10$ with different basis numbers $n$ in Figure~\ref{fig:lidDriven:streamline2}, which indicates that our ROC algorithm can give an accurate approximation for the whole {time evolution history as} the basis number $n$ is increased.

\section{Conclusion}
\label{sec:conclusion}
This paper proposes a novel  reduced over-collocation method {with adaptive enrichment} for efficiently solving parametrized fluid systems. EIM/GEIM hyper reduction techniques produce two sets of collocation points {whose combined cardinality is much lower than the size of the full model, while} an adaptive reduced residual gives the greedy algorithm {an error indicator}. The resulting adaptive R2-ROC is online efficient,  and {stabilizes} the time dependent Navier-Stokes equations {with our adaptive enrichment strategy judiciously adding points on the physical domain for the reduced solver to collocate on}.

\bibliographystyle{abbrv}
\bibliography{rbmbib}

\begin{thebibliography}{10}

\bibitem{astrid2004fast}
P.~Astrid.
\newblock Fast reduced order modeling technique for large scale {LTV} systems.
\newblock In {\em Proceedings of the 2004 American Control Conference},
  volume~1, pages 762--767. IEEE, 2004.

\bibitem{Astrid2008MissingPoint}
P.~{Astrid}, S.~{Weiland}, K.~{Willcox}, and T.~{Backx}.
\newblock Missing point estimation in models described by proper orthogonal
  decomposition.
\newblock {\em IEEE Transactions on Automatic Control}, 53(10):2237--2251,
  2008.

\bibitem{ballarin2015supremizer}
F.~Ballarin, A.~Manzoni, A.~Quarteroni, and G.~Rozza.
\newblock Supremizer stabilization of {POD--Galerkin} approximation of
  parametrized steady incompressible {Navier--Stokes} equations.
\newblock {\em Int. J. Numer. Methods Eng.}, 102(5):1136--1161, 2015.

\bibitem{Barrault2004}
M.~Barrault, Y.~Maday, N.~C. Nguyen, and A.~T. Patera.
\newblock An `empirical interpolation' method: {A}pplication to efficient
  reduced-basis discretization of partial differential equations.
\newblock {\em C. R. Math.}, 339(9):667--672, 2004.

\bibitem{BenaceurEhrlacherErnMeunier2018}
A.~Benaceur, V.~Ehrlacher, A.~Ern, and S.~Meunier.
\newblock A progressive reduced basis/empirical interpolation method for
  nonlinear parabolic problems.
\newblock {\em SIAM J. Sci. Comput.}, 40(5):A2930--A2955, 2018.

\bibitem{BennerGugercinWillcox2015}
P.~Benner, S.~Gugercin, and K.~Willcox.
\newblock A survey of projection-based model reduction methods for parametric
  dynamical systems.
\newblock {\em SIAM Review}, 57(4):483--531, 2015.

\bibitem{BinevCohenDahmenDevorePetrovaWojtaszczyk}
P.~Binev, A.~Cohen, W.~Dahmen, R.~Devore, G.~Petrova, and P.~Wojtaszczyk.
\newblock {Convergence Rates for Greedy Algorithms in Reduced Basis Methods}.
\newblock {\em SIAM J. MATH. ANAL}, pages 1457--1472, 2011.

\bibitem{bos2004accelerating}
R.~Bos, X.~Bombois, and P.~Van~den Hof.
\newblock Accelerating large-scale non-linear models for monitoring and control
  using spatial and temporal correlations.
\newblock In {\em Proceedings of the 2004 American Control Conference},
  volume~4, pages 3705--3710. IEEE, 2004.

\bibitem{CarlbergBouMoslehFarhat2011}
K.~Carlberg, C.~Bou-Mosleh, and C.~Farhat.
\newblock Efficient non-linear model reduction via a least-squares
  {P}etrov-{G}alerkin projection and compressive tensor approximations.
\newblock {\em Int. J. Numer. Methods Eng.}, 86(2):155--181, 2011.

\bibitem{Casenave2014_M2AN}
F.~Casenave, A.~Ern, and T.~Leli\`evre.
\newblock Accurate and online-efficient evaluation of the {\it a posteriori}
  error bound in the reduced basis method.
\newblock {\em ESAIM: M2AN}, 48(1):207--229, 2014.

\bibitem{ChaturantabutSorensen2010}
S.~Chaturantabut and D.~C. Sorensen.
\newblock {Nonlinear model reduction via discrete empirical interpolation}.
\newblock {\em SIAM J. Sci. Comput.}, 32(5):2737--2764, 2010.

\bibitem{MAC_chenlong}
L.~Chen.
\newblock Finite difference method for {Stokes} equations: {MAC} scheme.
\newblock {\em University of California, Irvine.
  https://www.math.uci.edu/~chenlong/226/MACStokes.pdf}.

\bibitem{ChenGottlieb2013}
Y.~Chen and S.~Gottlieb.
\newblock Reduced collocation methods: {Reduced basis methods in the
  collocation framework}.
\newblock {\em J. Sci. Comput.}, 55(3):718--737, 2013.

\bibitem{ChenSigalJiMaday2021}
Y.~Chen, S.~Gottlieb, L.~Ji, and Y.~Maday.
\newblock An {EIM-degradation} free reduced basis method via over collocation
  and residual hyper reduction-based error estimation.
\newblock {\em J. Comput. Phys.}, 444:110545, 2021.

\bibitem{chen2012certified}
Y.~Chen, J.~S. Hesthaven, Y.~Maday, J.~Rodriguez, and X.~Zhu.
\newblock Certified reduced basis method for electromagnetic scattering and
  radar cross section estimation.
\newblock {\em Comput. Methods Appl. Mech. Eng.}, 233-236:92--108, 2012.

\bibitem{ChenJiNarayanXu2020}
Y.~Chen, L.~Ji, A.~Narayan, and Z.~Xu.
\newblock L1-based reduced over collocation and hyper reduction for steady
  state and time-dependent nonlinear equations.
\newblock {\em J. Sci. Compt.}, 87:10, 2021.

\bibitem{JiangChenNarayan2019}
Y.~Chen, J.~Jiang, and A.~Narayan.
\newblock {A robust error estimator and a residual-free error indicator for
  reduced basis methods}.
\newblock {\em Computers \& Mathematics with Applications}, 77:1963--1979,
  2019.

\bibitem{CohenDeVore2015}
A.~Cohen and R.~DeVore.
\newblock {Approximation of high-dimensional parametric {\{}PDE{\}}s}.
\newblock {\em Acta Numer.}, 24:1--159, 2015.

\bibitem{deparis2009reduced}
S.~Deparis and G.~Rozza.
\newblock Reduced basis method for multi-parameter-dependent steady
  {Navier--Stokes} equations: applications to natural convection in a cavity.
\newblock {\em J. Comput. Phys.}, 228(12):4359--4378, 2009.

\bibitem{diez2017generalized}
P.~D{\'\i}ez, S.~Zlotnik, and A.~Huerta.
\newblock Generalized parametric solutions in {Stokes} flow.
\newblock {\em Comput. Methods Appl. Mech. Eng.}, 326:223--240, 2017.

\bibitem{farhat2014dimensional}
C.~Farhat, P.~Avery, T.~Chapman, and J.~Cortial.
\newblock Dimensional reduction of nonlinear finite element dynamic models with
  finite rotations and energy-based mesh sampling and weighting for
  computational efficiency.
\newblock {\em Int. J. Numer. Methods Eng.}, 98(9):625--662, 2014.

\bibitem{fick2017reduced}
L.~Fick, Y.~Maday, A.~T. Patera, and T.~Taddei.
\newblock A reduced basis technique for long-time unsteady turbulent flows.
\newblock {\em arXiv preprint arXiv:1710.03569}, 2017.

\bibitem{grepl2007efficient}
M.~A. Grepl, Y.~Maday, N.~Nguyen, and A.~T. Patera.
\newblock Efficient reduced-basis treatment of nonaffine and nonlinear partial
  differential equations.
\newblock {\em ESAIM: M2AN}, 41(3):575--605, 2007.

\bibitem{grepl2005posteriori}
M.~A. Grepl and A.~T. Patera.
\newblock A posteriori error bounds for reduced-basis approximations of
  parametrized parabolic partial differential equations.
\newblock {\em ESAIM: M2AN}, 39(1):157--181, 2005.

\bibitem{Haasdonk2017Review}
B.~Haasdonk.
\newblock {\em Chapter 2: Reduced basis methods for parametrized {PDE}s--a
  tutorial introduction for stationary and instationary problems}, volume~15,
  pages 65--136.
\newblock SIAM Philadelphia, 2017.

\bibitem{harlow1965numerical}
F.~H. Harlow and J.~E. Welch.
\newblock Numerical calculation of time-dependent viscous incompressible flow
  of fluid with free surface.
\newblock {\em The Physics of Fluids}, 8(12):2182--2189, 1965.

\bibitem{HesthavenRozzaStammBook}
J.~S. Hesthaven, G.~Rozza, B.~Stamm, et~al.
\newblock {\em Certified reduced basis methods for parametrized partial
  differential equations}, volume 590.
\newblock Springer, 2016.

\bibitem{HKCHP}
D.~B.~P. Huynh, D.~J. Knezevic, Y.~Chen, J.~S. Hesthaven, and A.~T. Patera.
\newblock A natural-norm {Successive Constraint Method} for inf-sup lower
  bounds.
\newblock {\em Comput. Methods Appl. Mech. Eng.}, 199:1963--1975, 2010.

\bibitem{HuynhSCM}
D.~B.~P. Huynh, G.~Rozza, S.~Sen, and A.~T. Patera.
\newblock {A successive constraint linear optimization method for lower bounds
  of parametric coercivity and inf-sup stability constants}.
\newblock {\em C. R. Acad. Sci. Paris, S$\acute{\rm e}$rie I.}, 345:473--478,
  2007.

\bibitem{iliescu2014variational}
T.~Iliescu and Z.~Wang.
\newblock Variational multiscale proper orthogonal decomposition:
  {Navier--Stokes} equations.
\newblock {\em Numer. Methods Partial Differ. Equ.}, 30(2):641--663, 2014.

\bibitem{ito1998reduced}
K.~Ito and S.~S. Ravindran.
\newblock A reduced-order method for simulation and control of fluid flows.
\newblock {\em J. Comput. Phys.}, 143(2):403--425, 1998.

\bibitem{MadayMulaPateraYano2015}
Y.~Maday, O.~Mula, A.~Patera, and M.~Yano.
\newblock The {Generalized Empirical Interpolation Method}: Stability theory on
  {Hilbert} spaces with an application to the {Stokes} equation.
\newblock {\em Comput. Methods Appl. Mech. Eng.}, 287:310 -- 334, 2015.

\bibitem{MadayMulaTurinici_GEIMSIAM}
Y.~Maday, O.~Mula, and G.~Turinici.
\newblock Convergence analysis of the generalized empirical interpolation
  method.
\newblock {\em SIAM J. Numer. Anal.}, 54(3):1713--1731, 2016.

\bibitem{MagicPt_2009}
Y.~Maday, N.~C. Nguyen, A.~T. Patera, and S.~H. Pau.
\newblock A general multipurpose interpolation procedure: {The} magic points.
\newblock {\em Communications on Pure \& Applied Analysis}, 8(1):383, 2009.

\bibitem{mu2017simple}
L.~Mu and X.~Ye.
\newblock A simple finite element method for the {Stokes} equations.
\newblock {\em Advances in Computational Mathematics}, 43(6):1305--1324, 2017.

\bibitem{nicolaides1992analysis}
R.~A. Nicolaides.
\newblock Analysis and convergence of the {MAC} scheme. {I}. the linear
  problem.
\newblock {\em SIAM J. Numer. Anal.}, 29(6):1579--1591, 1992.

\bibitem{PeherstorferButnaruWillcoxBungartz2014}
B.~Peherstorfer, D.~Butnaru, K.~Willcox, and H.~Bungartz.
\newblock Localized discrete empirical interpolation method.
\newblock {\em SIAM J. Sci. Comput.}, 36(1):A168--A192, 2014.

\bibitem{Quarteroni2015}
A.~Quarteroni, A.~Manzoni, and F.~Negri.
\newblock {\em Reduced basis methods for partial differential equations: {A}n
  introduction}, volume~92 of {\em Springer Series in Computational
  Mathematics}.
\newblock Springer, 2015.

\bibitem{rocha2020adaptive}
I.~Rocha, F.~van~der Meer, and L.~J. Sluys.
\newblock An adaptive domain-based {POD/ECM} hyper-reduced modeling framework
  without offline training.
\newblock {\em Computer Methods in Applied Mechanics and Engineering},
  358:112650, 2020.

\bibitem{Rozza2008}
G.~Rozza, D.~B.~P. Huynh, and A.~T. Patera.
\newblock Reduced basis approximation and a posteriori error estimation for
  affinely parametrized elliptic coercive partial differential equations.
\newblock {\em Arch. Comput. Methods Eng.}, 15(3):229--275, 2008.

\bibitem{rozza2007stability}
G.~Rozza and K.~Veroy.
\newblock On the stability of the reduced basis method for {Stokes} equations
  in parametrized domains.
\newblock {\em Comput. Methods Appl. Mech. Eng.}, 196(7):1244--1260, 2007.

\bibitem{ryckelynck2005priori}
D.~Ryckelynck.
\newblock A priori hyperreduction method: {An} adaptive approach.
\newblock {\em J. Comput. Phys.}, 202(1):346--366, 2005.

\bibitem{ryckelynck2009hyper}
D.~Ryckelynck.
\newblock Hyper-reduction of mechanical models involving internal variables.
\newblock {\em Int. J. Numer. Methods Eng.}, 77(1):75--89, 2009.

\bibitem{stabile2018finite}
G.~Stabile and G.~Rozza.
\newblock Finite volume {POD-Galerkin} stabilised reduced order methods for the
  parametrised incompressible {Navier--Stokes} equations.
\newblock {\em Computers \& Fluids}, 173:273--284, 2018.

\bibitem{stabile2020efficient}
G.~Stabile, M.~Zancanaro, and G.~Rozza.
\newblock Efficient geometrical parametrization for finite-volume-based reduced
  order methods.
\newblock {\em International Journal for Numerical Methods in Engineering},
  121(12):2655--2682, 2020.

\bibitem{wang2012proper}
Z.~Wang, I.~Akhtar, J.~Borggaard, and T.~Iliescu.
\newblock Proper orthogonal decomposition closure models for turbulent flows: A
  numerical comparison.
\newblock {\em Comput. Methods Appl. Mech. Eng.}, 237:10--26, 2012.

\bibitem{yano2019lp}
M.~Yano and A.~T. Patera.
\newblock An {LP} empirical quadrature procedure for reduced basis treatment of
  parametrized nonlinear {PDEs}.
\newblock {\em Computer Methods in Applied Mechanics and Engineering},
  344:1104--1123, 2019.

\end{thebibliography}

\end{document}